\documentclass[3p, number, sort&compress, times, lefttitle, dvipsnames]{elsarticle}

\journal{Elsevier}

\usepackage[utf8]{inputenc}
\usepackage{gensymb}
\usepackage{hyphenat}
\usepackage{subcaption}
\usepackage{multirow}
\usepackage{graphicx}
\usepackage{tikz}
\usetikzlibrary{patterns}
\definecolor{pal29}{rgb}{0.890, 0.953, 0.984}
\definecolor{pal58}{rgb}{0.184, 0.310, 0.455}

\usepackage{amsmath,amsfonts}
\renewcommand{\Re}{\mathbb{R}}

\usepackage{amsthm}
\theoremstyle{plain}

\theoremstyle{definition}
\usepackage{bm}
\newcommand{\vm}[1]{\bm{#1}}
\newcommand{\vx}{\vm{x}}

\usepackage{mathtools}

\usepackage[colorlinks=true]{hyperref}

\newcommand{\revone}[1]{#1}
\newcommand{\revtwo}[1]{#1}

\newcommand{\tref}[1]{Table~\ref{#1}}
\newcommand{\fref}[1]{Fig.~\ref{#1}}
\newcommand{\sref}[1]{Section~\ref{#1}}

\begin{document}

\title{Scaled boundary cubature scheme for numerical integration over
\revone{planar regions with affine and curved boundaries}}

\author[1]{Eric B.\ Chin\corref{cor1}}
\ead{chin23@llnl.gov}

\author[2]{N.\ Sukumar}

\cortext[cor1]{Corresponding author}

\address[1]{Lawrence Livermore National Laboratory, 7000 East Avenue,
Livermore, CA 94550, USA}

\address[2]{Department of Civil and Environmental Engineering, 
University of California, Davis, CA 95616, USA}

\begin{abstract}
This paper introduces the scaled boundary cubature (SBC) scheme for accurate and
efficient integration of functions over polygons and two\hyp{}dimensional
regions bounded by parametric curves.  Over two\hyp{}dimensional domains, the
SBC method reduces integration over a region bounded by $m$ curves to
integration over $m$ regions (referred to as curved triangular regions), where
each region is bounded by two line segments and a curve.  With proper
(counterclockwise) orientation of the boundary curves, the scheme is applicable
to convex and nonconvex domains.  Additionally, for star\hyp{}convex domains, a
tensor\hyp{}product cubature rule with positive weights and integration points
in the interior of the domain is obtained. If the integrand is homogeneous, we
show that this new method reduces to the homogeneous numerical integration
scheme; however, the SBC scheme is more versatile since it is equally applicable
to both homogeneous and non\hyp{}homogeneous functions. This paper also
introduces several methods for smoothing integrands with point singularities and
near\hyp{}singularities.  When these methods are \revone{used}, highly efficient
integration of weakly singular functions is realized. The SBC method is
applied to a number of benchmark problems, which reveal its broad applicability
and superior performance (in terms of time to generate a rule and accuracy per
cubature point) when compared to existing methods for integration.
\end{abstract}


\maketitle

\section{Introduction}
\revone{In this paper, we propose} a new method for the accurate and efficient
numerical integration of functions over \revone{planar (two\hyp{}dimensional)
regions bounded by affine and/or curved edges, such as polygons or conics}.  The
classes of functions include those that are either continuously differentiable
in the interior of the domain, or are weakly \revone{singular} or nearly singular at \revone{a point within or near}
the domain.  \revone{Regions of integration may be convex or nonconvex and the
boundary of the region is described by parametric curves.}  The boundary curves
can be given in closed\hyp{}form or represented using Bernstein polynomials,
B\hyp{}splines, or \revone{non\hyp{}uniform rational B\hyp{}splines (NURBS)}.
For integration over bounded domains with curved boundaries, we show that the
homogeneous numerical integration
scheme~\cite{Lasserre:1998:ICP,Chin:2015:NIH,Chin:2019:MCI} is a reduction of a
tensor\hyp{}product cubature scheme that is obtained on using a scaled boundary
(SB) parametrization~\cite{Song:1997:TSB,Wolf:2001:TSB} of the domain.  We refer
to this integration method as the \textit{scaled boundary cubature} (SBC)
scheme.  \revone{While only two\hyp{}dimensional domains are considered in this
paper, we note the SBC scheme is also applicable to} three\hyp{}dimensional
domains that are bounded by parametric surfaces.

\subsection{Related work}
Methods for numerical integration over polygons and domains bounded by
parametric curves is currently an area of active research. An early, but still
popular method of cubature is partitioning the integration domain, and then
using well\hyp{}known rules of integration on the resulting (potentially curved)
triangular and quadrilateral
partitions~\cite{Moes:1999:AFE,Sevilla:2008:NEF,Fries:2017:HOM,Artioli:2020:ACP}.
However, the process of generating a partition can be time consuming for
complex geometries and cubature rules generated from this process are far from
optimal.  Optimized rules that retain polynomial precision can be
formulated through the process of moment
fitting~\cite{Mousavi:2010:GGQ,Artioli:2020:ACP}, which reduces the number of
cubature points iteratively by the least\hyp{}squares Newton method.  However,
this process is computationally expensive since it requires the solution of many
linear systems of equations, and furthermore, moment fitting requires an initial
cubature rule that is typically generated through partitioning.  For the special
case of a hyperrectangle intersecting an implicitly defined surface,
high\hyp{}order accurate integration schemes have also been devised.  The need
for this type of integration arises in isogeometric analysis (IGA) with trimmed
patches and with unfitted finite element methods.  One such integration method
by Saye~\cite{Saye:2015:HOQ} introduces a height function that approximates the
location of the interface to high degree.  Another method by Scholz and
J\"{u}ttler~\cite{Scholz:2019:NIT,Scholz:2020:HOQ} introduce higher order
correction terms based on transport theorems to reduce error in integration on
curved surfaces.

\smallskip
To avoid partitioning, some recent contributions have utilized generalized
Stokes's theorem to develop cubature rules that scale with the number of
edges/curves in the polygonal/curved domains.  For instance, in Sommariva and
Vianello~\cite{Sommariva:2007:PGC}, Gauss\hyp{}Legendre quadrature is used to
compute the antiderivatives in Green's formula (Stokes's theorem in two
dimensions), producing an exact cubature rule for integrating polynomials over a
convex or nonconvex polygon.  With this scheme, named Gauss\hyp{}Green (G-G)
cubature, cubature points can appear outside the polygonal domain of
integration, though the accuracy of the rule is not necessarily affected.  This
approach is also viable on curved domains, as demonstrated by Sommariva and
Vianello~\cite{Sommariva:2009:GGC} and more recently by Gunderman et
al.~\cite{Gunderman:2021:SMF}.  Another approach utilizing generalized Stokes's
theorem is the homogeneous numerical integration (HNI) method.  Building on work
by Lasserre~\cite{Lasserre:1998:ICP,Lasserre:1999:IHF}, the HNI
method~\cite{Chin:2015:NIH} combines generalized Stokes's theorem with Euler's
homogeneous function theorem to produce a cubature rule for integrating
polynomials and homogeneous functions over polytopes and regions bounded by
homogeneous curves.  More recently in Chin and
Sukumar~\cite{Chin:2019:MCI,Chin:2020:AEM}, this idea is extended to domains
bounded by arbitrary parametric curves such as polynomial and rational
B\'{e}zier and B\hyp{}spline curves.  With the HNI method, a simple formula is
realized for reducing the dimension of integration by one, resulting in a very
efficient cubature rule.  However, the restriction to homogeneous functions
limits its broader appeal, and in addition, the implementation of the scheme is
also more involved~\cite{Chin:2017:MCD,Chin:2019:MCI}.

\smallskip
Novel ideas in the computer simulation of physical phenomena and
computer\hyp{}aided manufacturing have motivated the need for improved numerical
integration over polygons and regions bounded by curves.  \revtwo{Galerkin
methods based on the SB mapping, such as the ones explored by Chen et
al.~\cite{Chen:2016:ANB} and Klinkel and Reichel~\cite{Klinkel:2019:AFE} would benefit
from integration rules formulated from the scaled boundary parametrization.}
More broadly, techniques such as the extended finite element method
(X-FEM)~\cite{Moes:1999:AFE,Sukumar:2015:EFE}, polygonal finite element
methods~\cite{Hormann:2017:GBC}, trimmed patches with
IGA~\cite{Hughes:2005:IAC,Marussig:2018:ICI}, the \revone{NURBS} enhanced finite
element method~\cite{Sevilla:2008:NEF}, unfitted finite element methods such as
embedded interface and fictitious domain
methods~\cite{Bishop:2003:RSA,Burman:2014:CDG,Schillinger:2015:TFC}, virtual
element
method~\cite{BeiraodaVeiga:2013:BPV,BeiraodaVeiga:2019:TVE,Benvenuti:2019:EVE},
and discontinuous Galerkin methods~\cite{Lew:2008:ADG,Cangiani:2014:HVD} \revone{need}
accurate, efficient cubature over polygons and regions bounded by
curves.  Additionally, in the field of contact mechanics and domain
decomposition, accurately computing contributions to the force vector in a
mortar method requires integration over polygons~\cite{Puso:2004:AMS}, and for
\revone{high\hyp{}order} mortar methods, integration over planar regions bounded by
curves~\cite{Hesch:2011:TTD} is needed.  Furthermore, for applications in
computer graphics and animation~\cite{Guendelman:2003:NRB} and
computer\hyp{}aided design~\cite{Krishnamurthy:2011:AGA}, computing areas,
centroids, and moments of inertia require accurate integration over polygons and
regions bounded by parametric curves.

\subsection{Contributions and outline}
The HNI method provides a simple formula for reducing integration over a
$d$\hyp{}dimensional domain to scaled integrals over the
$(d\!-\!1)$\hyp{}dimensional boundary facets of the domain.  Therefore, cubature
points lie entirely on the boundary and only $\mathcal{O}(mp)$ integration
points are needed to integrate a homogeneous polynomial of degree $(2p-1)$
over an $m$\hyp{}gon.  The proof of the method utilizes generalized Stokes's
theorem to relate the domain integral to an equivalent integral over the
domain's boundary.  However, instead of computing antiderivatives of the
integrand, the HNI method utilizes Euler's homogeneous function theorem to
relate the integrand to its derivative.  While not immediately apparent from its
derivation, the HNI method implicitly and exactly computes a
one\hyp{}dimensional integral, which when combined with a parametrization for
the boundary of the domain, parametrizes the entire domain of integration. In
this paper, we characterize this parametrization, and show that it matches the
domain mapping (SB parametrization) introduced in the scaled boundary finite
element method~\cite{Song:1997:TSB,Wolf:2001:TSB}.  In this framework, the HNI
method can be cast as a special case of the mapped integral, which is valid when
the integrand is homogeneous.

\smallskip
Beyond the additional insights that the SB parametrization affords into the HNI
method, we also demonstrate that integrating over the mapped space delivers an
effective cubature rule for polytopes and curved solids in its own right.  The
domain need not be convex or simply\hyp{}connected if the integrand is
continuously differentiable and computable outside the domain of integration.
The method relies on tensor\hyp{}product quadrature rules over the unit square,
and requires $\mathcal{O}(mp^2)$ integration points to exactly integrate a
$(2p\!-\!1)$\hyp{}degree bivariate polynomial over an $m$-gon.  Moreover, the
limitation of homogeneous integrands is removed---non\hyp{}homogeneous functions
and functions \revone{that are} not known explicitly can be integrated.  In addition, we introduce
a mapping of the domain that we refer to as the \textit{generalized} SB
parametrization and integral transformations that smooth singular integrands,
enabling accurate cubature of weakly singular functions with relatively few
cubature points.  The name ``generalized SB parametrization'' reflects its
similarity to the generalized Duffy transformation introduced by Mousavi and
Sukumar~\cite{Mousavi:2010:GDT}.  These transformations can be applied to
computational methods that require the integration of weakly singular
functions, such as the X-FEM or the boundary element method (BEM).

\smallskip
The remainder of this paper is organized as follows.  In~\sref{sec:hni}, the HNI
method is introduced, focusing specifically on how integrals over polygons and
curved solids are decomposed using this method.  In~\sref{sec:radial-xform}, we
introduce the SBC method and develop several simplifications for polygonal
domains.  Transformations appropriate for integration of weakly singular
functions are presented in~\sref{sec:singular}, which increase cubature accuracy
for a given number of cubature points.  Next, in~\sref{sec:examples}, we present
numerical integration examples that arise in several distinct application areas
to establish the accuracy and to showcase the capabilities of the integration
method.  These application areas include polygonal finite element methods, the
X-FEM, and the BEM.  The examples demonstrate integration of polynomials, smooth
nonpolynomial functions such as the test functions in
Franke~\cite{Franke:1979:ACC}, and weakly singular functions over convex and
nonconvex polygons and domains bounded by parametric curves.
In~\sref{sec:ex-tmvi}, we apply the SBC scheme to compute the $L_2$ error for
transfinite mean value interpolation (TMVI)~\cite{Ju:2005:MVC,Dyken:2009:TMV}
over a curved domain.  Finally, we conclude with \sref{sec:conclusion}, where we
recap the main developments in this paper and provide a few future directions of
research.

\section{Domain partitioning into curved triangles}\label{sec:hni} 
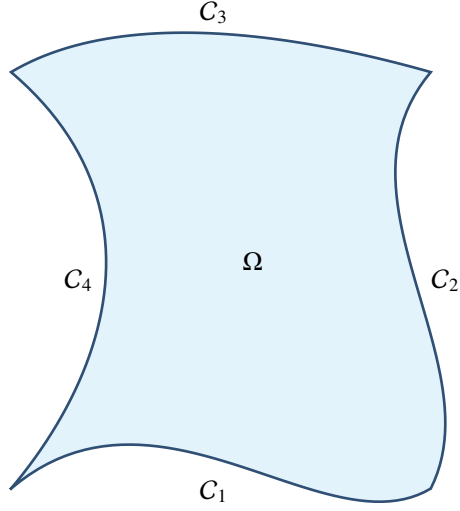
\begin{figure}
  \centering
  \begin{tikzpicture}
    \node at (0,0) {\includegraphics[scale=1]{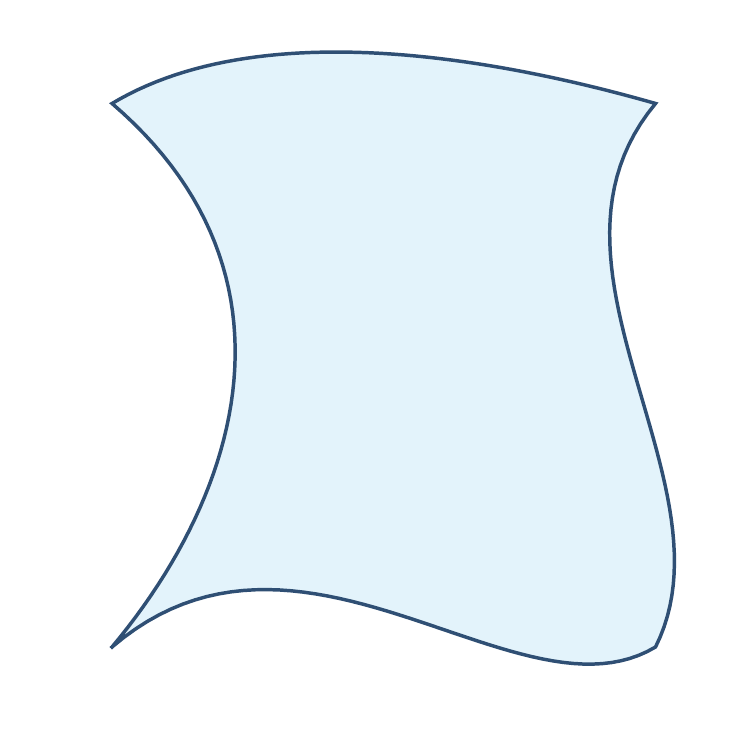}};
    \node at (0,-1.1in) {$\mathcal{C}_1$};
    \node at (1.2in,0) {$\mathcal{C}_2$};
    \node at (0,1.4in) {$\mathcal{C}_3$};
    \node at (-0.7in,0) {$\mathcal{C}_4$};
    \node at (0.2in,0.1in) {$\Omega$};
  \end{tikzpicture}
  \caption{A region bounded by cubic B\'{e}zier curves.}
  \label{fig:2d-curved-region}
\end{figure}

In $\Re^2$,
we define a closed region $\Omega \subset \Re^2$ with piecewise smooth,
orientable boundary $\partial \Omega$.  We assume that $\partial \Omega$ is
composed of $m$ parametric curves, $\mathcal{C}_i$ for $i = 1, \dotsc, m$,
\revtwo{where each curve occupies a unique subset of the boundary,} such that
\begin{equation}\label{eq:bdry-curves}
  \partial \Omega = \bigcup_{i=1}^m \mathcal{C}_i  \revtwo{\quad \text{ and } \quad
    \mathcal{C}_j \cap \mathcal{C}_k = \emptyset 
    \text{ for }j,k=1,\dotsc,m \text{ and } j \neq k}.
\end{equation}
An illustration of a region that conforms to this definition is provided
in~\fref{fig:2d-curved-region}.  All boundary curves are ordered sequentially as
the boundary is traversed in a counterclockwise direction.  We introduce a
coordinate $\vx_0$ that serves as the center of the partitioning, which leads to
a decomposition of $\Omega$ into partitions bounded by three curves: two line
segments and one curve.  We will refer to these partitions as \textit{curved
triangles} or by the symbol $\mathcal{T}$.

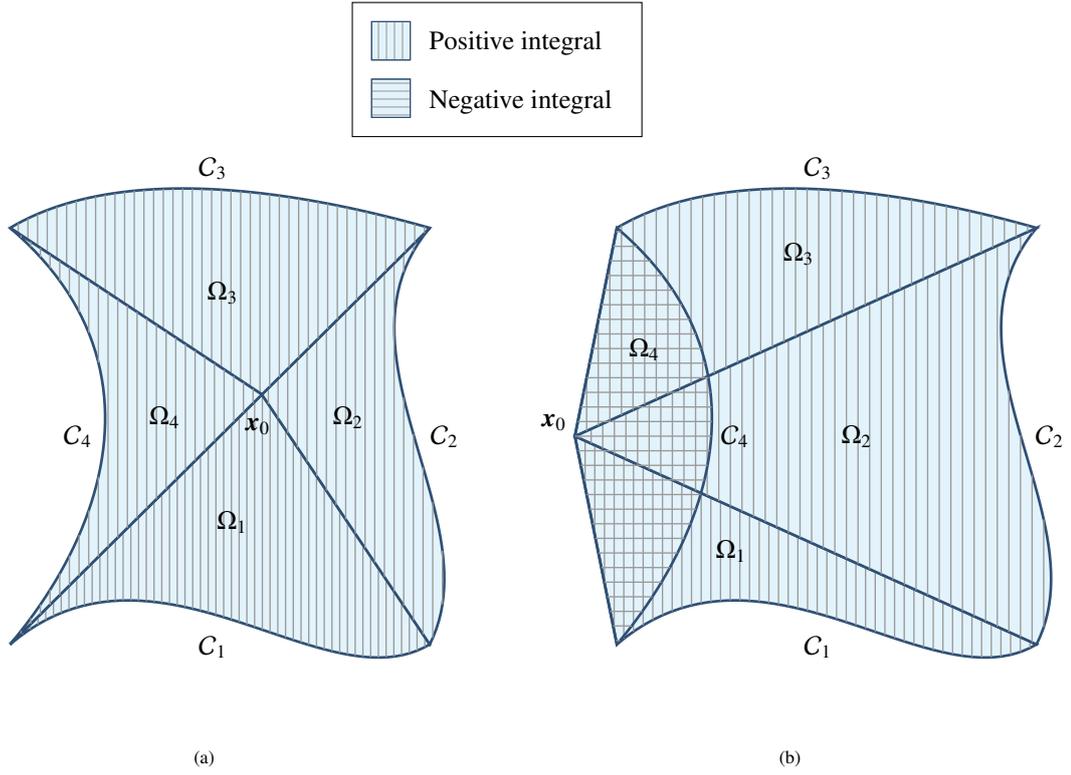
\begin{figure}[t]
  \centering
  \begin{tikzpicture}
    \draw[pal58, fill=pal29] (0,0) rectangle (0.2in,0.2in);
    \draw[pal58, pattern=vertical lines, pattern color=black!40!white] (0,0) rectangle (0.2in,0.2in);
    \node[anchor=west] at (0.25in,0.09in) {Positive integral};
    \draw[pal58, fill=pal29] (0,-0.1in) rectangle (0.2in,-0.3in);
    \draw[pal58, pattern=horizontal lines, pattern color=black!40!white] (0,-0.1in) rectangle (0.2in,-0.3in);
    \node[anchor=west] at (0.25in,-0.22in) {Negative integral};
    \draw (-0.1in,0.3in) rectangle (1.4in,-0.4in);
  \end{tikzpicture}\\
  \begin{subfigure}{3in}
    \begin{tikzpicture}
      \node at (0,0) {\includegraphics[scale=1]{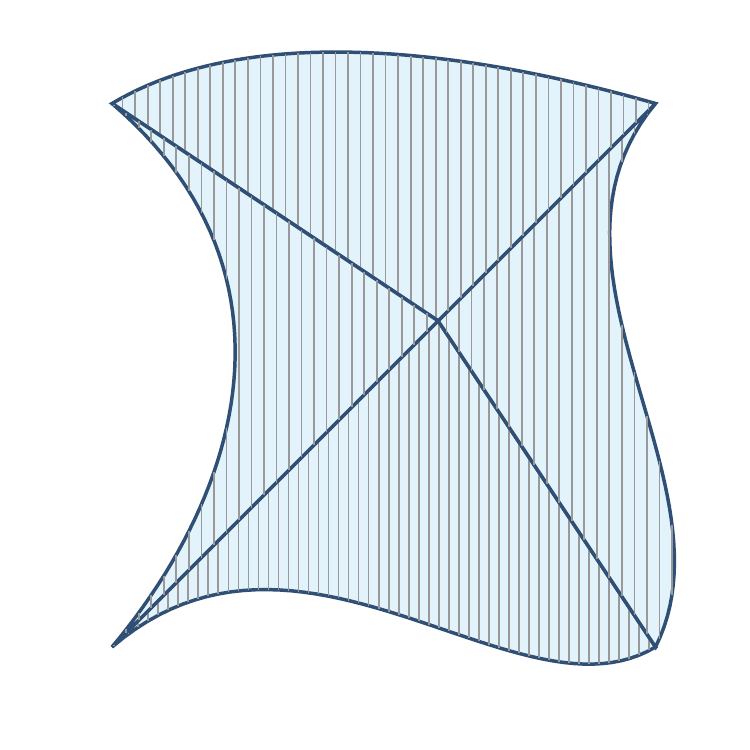}};
      \node at (0,-1.1in) {$\mathcal{C}_1$};
      \node at (1.2in,0) {$\mathcal{C}_2$};
      \node at (0,1.4in) {$\mathcal{C}_3$};
      \node at (-0.7in,0) {$\mathcal{C}_4$};
      \node at (0.1in,-0.45in) {$\Omega_1$};
      \node at (0.7in,0.1in) {$\Omega_2$};
      \node at (0.05in,0.75in) {$\Omega_3$};
      \node at (-0.25in,0.1in) {$\Omega_4$};
      \node at (0.23in,0.05in) {$\vx_0$};
    \end{tikzpicture}
    \caption{}\label{fig:2d-curved-region-star-convex}
  \end{subfigure}
  \begin{subfigure}{3in}
    \begin{tikzpicture}
      \node at (-0.1in,0) {\includegraphics[scale=1]{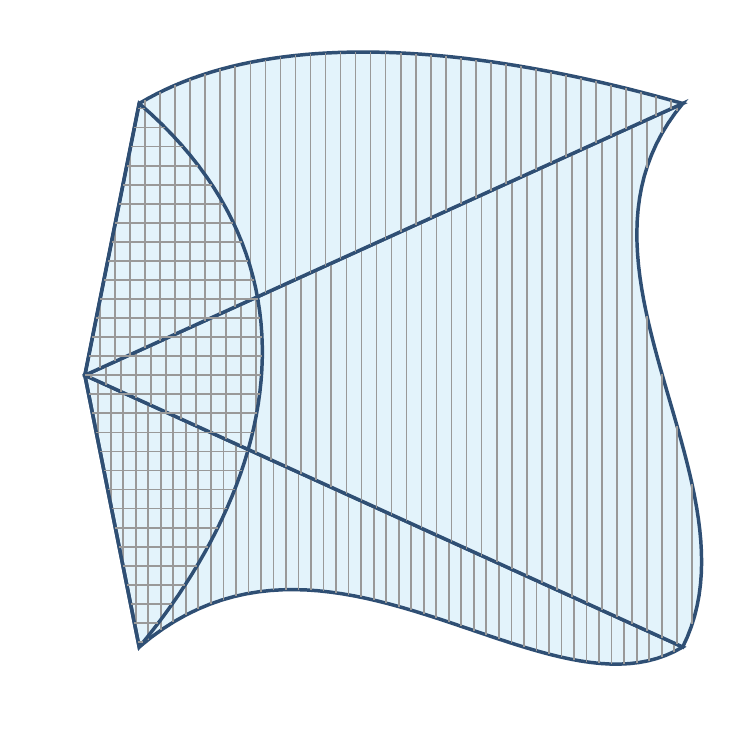}};
      \node at (0,-1.1in) {$\mathcal{C}_1$};
      \node at (1.2in,0) {$\mathcal{C}_2$};
      \node at (0,1.4in) {$\mathcal{C}_3$};
      \node at (-0.43in,0) {$\mathcal{C}_4$};
      \node at (-0.45in,-0.6in) {$\Omega_1$};
      \node at (0.2in,0) {$\Omega_2$};
      \node at (-0.1in,0.95in) {$\Omega_3$};
      \node at (-0.9in,0.45in) {$\Omega_4$};
      \node at (-1.37in,0.07in) {$\vx_0$};
    \end{tikzpicture}
    \caption{}\label{fig:2d-curved-region-no-star-convex}
  \end{subfigure}
  \caption{Two choices for decomposing the star\hyp{}convex domain $\Omega$
  in~\fref{fig:2d-curved-region}.  (a) This choice of $\vx_0$ results in
  positive regions of integration over all subdomains, (b) with this choice of
  $\vx_0$, the value of the integral over $\Omega_4$ must be negated from
  the sum of the integrals over the other subregions.}
  \label{fig:2d-curved-region-decomposed}
\end{figure}

\smallskip
If $\Omega$ is star\hyp{}convex, we can select a point $\vx_0 \in \Omega$ so
that the boundaries of the curved triangles do not intersect any of the boundary
curves, $\mathcal{C}_i$ ($i = 1, \dotsc, m$).  Or, more generally, we can choose
$\vx_0$ such that any line segment that begins at $\vx_0$ only intersects the
boundary at one point. This choice for $\vx_0$ with respect to the
star\hyp{}convexity of the domain is depicted
in~\fref{fig:2d-curved-region-star-convex}. However, even if $\Omega$ is not
star\hyp{}convex, or if $\vx_0$ is not selected with regards to the
star\hyp{}convexity of $\Omega$, accurate integration is still possible provided
a notion of positive and negative portions of a partition exists with the
integration rule.  A partitioning of $\Omega$ where positivity and negativity
must be accounted for is illustrated
in~\fref{fig:2d-curved-region-no-star-convex}.  Specifically, for
\fref{fig:2d-curved-region-no-star-convex}, the integral of a function $f$ over
the region is
\begin{equation}
  \int_\Omega f\,d\vx = \int_{\Omega_1} f\,d\vx + 
    \int_{\Omega_2} f\,d\vx +
    \int_{\Omega_3} f\,d\vx -
    \int_{\Omega_4} f\,d\vx ,
\end{equation}
where $\Omega_i$ is the curved triangle associated with $\mathcal{C}_i$. Note
the regions $\Omega_1$, $\Omega_2$, and $\Omega_3$ include $\Omega_4$, but
$\Omega_4$ is not in $\Omega$.  By subtracting the integral over $\Omega_4$, the
overall effect is to ignore the contribution of the integral outside $\Omega$.
To pursue this approach for the integration, the function $f$ must be computable
and continuously differentiable in $\Omega_1 \cup \ldots \cup \Omega_m$.  While
the example in~\fref{fig:2d-curved-region-no-star-convex} is applicable to
partitions whose integrals are entirely negative or positive, a single partition
could easily be defined that contains both positive and negative contributions
to the integral. This would occur if $\mathcal{C}_i$ intersects a line segment
starting at $\vx_0$ more than once.  Cubature rules with negative weights tend to
be avoided for strictly positive integrands, since they can compromise
accuracy due to round-off errors.  For these cases, we recommend further
subdivision of the integration domain until all subdivided integration domains
are star\hyp{}convex, coupled with careful selection of $\vx_0$.

\subsection{Curved triangles in the HNI method}\label{sec:hni-curved-triangles}
In this section, we introduce the HNI method and demonstrate it implicitly
divides $\Omega$, the domain of integration, into curved triangles, with the
origin serving as $\vx_0$.  Proof of the HNI method relies on Euler's
homogeneous function theorem and generalized Stokes's theorem.  A positively
homogeneous function $h(\vx) : V \rightarrow \Re$ with degree of homogeneity $q
\in \Re$ is defined such that
\begin{equation}\label{eq:homogeneous-def}
  h(\lambda \vx) = \lambda^q h(\vx) \quad \forall \vx \in V 
\end{equation}
and for any real constant $\lambda > 0$.  In~\eqref{eq:homogeneous-def}, the
domain $V = \Re^2$ when $h(\vx)$ is a polynomial and $V = \Re^2 \backslash
\{\vm{0} \}$ when $h(\vx) = \|\vx\|^\beta$ for $\beta < 0$.  For a continuously
differentiable, positively homogeneous scalar\hyp{}valued function $h := h(\vx)$
of degree $q$, Euler's homogeneous function theorem is:
\begin{equation}\label{eq:euler}
  \vx \cdot \nabla h = q h .
\end{equation}
For the vector field $\vx h$, we can write $\nabla \cdot (\vx h) = (\nabla \cdot
\vx) h + \vx \cdot \nabla h$, and on invoking Stokes's theorem (divergence
theorem or Green's theorem in the plane), we have
\begin{equation}\label{eq:stokes-thm}
  \int_\Omega ( \nabla \cdot \vx ) h \, d \vx
      + \int_\Omega \vx \cdot \nabla h \, d\vx
    = \int_\Omega \nabla \cdot (\vx h) \, d\vx
    = \int_{\partial \Omega} 
      ( \vx \cdot \vm{n} ) h \, ds ,
\end{equation}
where $\partial \Omega$ is the boundary of $\Omega$, $\vm{n}$ is the 
unit outward 
normal to $\mathcal{C}_i$, and $ds$ is the differential length on
$\mathcal{C}_i$.  Substituting $\nabla \cdot \vx = 2$ and Euler's homogeneous
function theorem~\eqref{eq:euler}, \eqref{eq:stokes-thm} becomes
\begin{equation}\label{eq:hni1}
  (2+q) \int_\Omega h \, d\vx
  = \int_\Omega \nabla \cdot (\vx h) \, d\vx
  = \int_{\partial \Omega} 
    ( \vx \cdot \vm{n} ) h \, ds .
\end{equation}

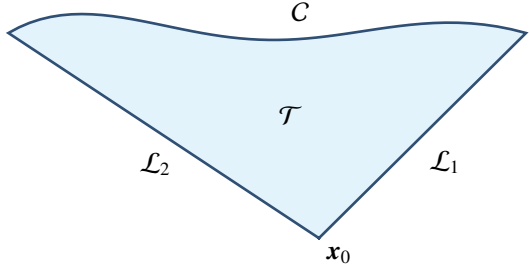
\begin{figure}
  \centering
  \begin{tikzpicture}
    \node at (0,0) {\includegraphics[scale=1]{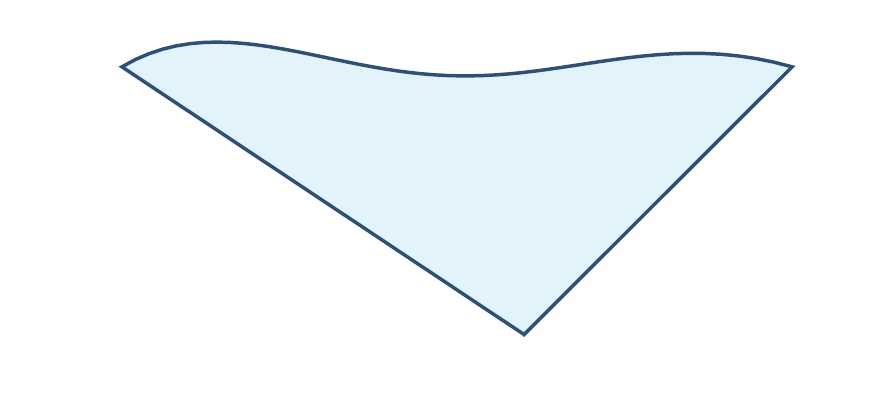}};
    \node at (0.25in,0.64in) {$\mathcal{C}$};
    \node at (1.0in,-0.14in) {$\mathcal{L}_1$};
    \node at (-0.5in,-0.15in) {$\mathcal{L}_2$};
    \node at (0.2in,0.1in) {$\mathcal{T}$};
    \node at (0.45in,-0.63in) {$\vx_0$};
  \end{tikzpicture}
  \caption{A curved triangle.}
  \label{fig:curved-triangle}
\end{figure}

\smallskip
Given a curved triangular domain $\mathcal{T}$, the boundary is the union of two
line segments: $\mathcal{L}_1$ and $\mathcal{L}_2$ and a curve $\mathcal{C}$
(see \fref{fig:curved-triangle}).  Thus, over $\mathcal{T}$, \eqref{eq:hni1} is
\begin{equation}\label{eq:hni-curvedtri1}
  (2+q) \int_\mathcal{T} h \, d\vx
  = \int_\mathcal{C} ( \vx \cdot \vm{n} ) h \, ds 
    + \sum_{i=1}^2 \int_{\mathcal{L}_i} ( \vx \cdot \vm{n} ) h \, ds .
\end{equation}
We note that line segments $\mathcal{L}_1$ and $\mathcal{L}_2$ pass through the
point $\vx_0$.  On choosing $\vx_0 = \vm{0}$, any point on $\mathcal{L}_1
\backslash \{\vm{0}\}$ or $\mathcal{L}_2 \backslash \{\vm{0}\}$ is then
perpendicular to $\vm{n}$, and therefore, $\vx \cdot \vm{n} = 0$ for both line
segments. As a result, \eqref{eq:hni-curvedtri1} reduces to
\begin{equation}\label{eq:hni-curvedtri2}
  (2+q) \int_\mathcal{T} h \, d\vx
  = \int_\mathcal{C} ( \vx \cdot \vm{n} ) h \, ds .
\end{equation}
Since a parametric description of $\mathcal{C}$ is presumed to be available, an
explicit expression for $\vm{n}$ is developed.  We now define
\begin{equation}\label{eq:curve-param}
  \mathcal{C} = \biggl\{\; \vx \in \Re^2 \;:\; \vx := \vm{c} (t) 
    \;\text{ for } t \in [0, 1] \;\biggr\} ,
\end{equation}
where we have introduced $\vm{c} (t)$ as the equation of the curve parametrized
by $t \in [0, 1]$.  Bounds for $t$ are limited between zero and one to simplify
the expressions that follow.  Curves that utilize different bounds of the
parameter $t$ can be scaled to match the definition used herein.  We require the
curve parameter $t$ to follow a counterclockwise orientation.  For example, on
the curved triangle in~\fref{fig:curved-triangle}, a counterclockwise
orientation of $t$ would result in $\vm{c} (0)$ and $\vm{c} (1)$ being points on
$\mathcal{L}_1$ and $\mathcal{L}_2$, respectively.  Now, we observe that the
unit normal vector is the unit tangent vector rotated through an angle of $-\pi
/ 2$:
\begin{equation}\label{eq:curve-norm}
  \vm{n} = \frac{\vm{c}^{\prime\perp}(t)}{\| \vm{c}^\prime (t) \|} ,
\end{equation}
where $\vm{c}^{\prime\perp}(t) = \text{rot } \vm{c}^\prime (t) := \bigl[
c^\prime_2(t) \ \ -c^\prime_1(t) \bigr]^T$ is $\vm{c}^\prime (t)$ rotated through
$-\pi /2$ radians. Combining these results with~\eqref{eq:hni-curvedtri2} yields
\begin{equation}\label{eq:hni-curvedtri}
  \int_\mathcal{T} h \, d\vx
  = \frac{1}{2+q} \int_\mathcal{C} h \, \frac{\vm{c}(t) \cdot \vm{c}^{\prime\perp}(t)}
    {\| \vm{c}^\prime (t) \|} \, ds ,
\end{equation}
where $\vm{c}(t) \cdot \vm{c}^{\prime\perp}(t) = c_1(t) c_2^\prime (t) - c_2(t)
c_1^\prime (t)$.  On using the HNI method, the integral of a homogeneous
function $h$ over a curved triangle $\mathcal{T}$ with $\vx_0 = \vm{0}$ is
reduced to a line integral over the bounding curve $\mathcal{C}$.  We remark the
quantity $\bigl( \vm{c}_i (t) \cdot \vm{c}^{\prime\perp} (t) \bigr) / \|
\vm{c}^\prime (t) \|$ is the signed distance of the tangent line of $\vm{c}_i
(t)$ to the point $\vx_0 = \vm{0}$.

\smallskip
The boundary of $\Omega$ is defined in~\eqref{eq:bdry-curves}. Partitioning the
boundary and using the parametric form of the curve given
in~\eqref{eq:curve-param}, \eqref{eq:hni1} becomes
\begin{equation}\label{eq:hni2}
  \int_\Omega h \, d\vx
  = \frac{1}{2+q} \sum_{i=1}^m 
    \int_{\mathcal{C}_i} h \, \frac{\vm{c}_i(t) \cdot \vm{c}_i^{\prime\perp}(t)}
    {\| \vm{c}_i^\prime (t) \|} \, ds ,
\end{equation}
where $\vm{c}_i(t)$ is the parametric form of the curve $\mathcal{C}_i$.
Using~\eqref{eq:hni-curvedtri}, \eqref{eq:hni2} can be equivalently expressed as
\begin{equation}\label{eq:hni-curve-simplified}
  \int_\Omega h (\vx) \, d\vx
    = \sum_{i=1}^m \int_{\Omega_i} h(\vx) d\vx ,
\end{equation}
which shows that the HNI method calculates integrals through a subdivision of
$\Omega$ into curved triangles.

\subsection{Scaled boundary parameter in the HNI method}
As~\eqref{eq:hni-curvedtri} makes clear, the HNI method converts integration
from a two\hyp{}dimensional curved triangle to a one\hyp{}dimensional integral
over the curve $\mathcal{C}$.  This implies the HNI method exactly computes the
integral over another parameter, provided $h(\vx)$ is homogeneous.  We introduce
the map $\vm{\psi} : [0, \rho(t)] \times [0, 1] \rightarrow \mathcal{T}^+
\supseteq \mathcal{T}$ that maps the radial distance from $\vx_0$, $r = \|\vx -
\vx_0\|$, and the curve parameter to the domain $\mathcal{T}^+$.  If
$\mathcal{T}$ is star\hyp{}convex with respect to $\vx_0$, then $\mathcal{T}^+ =
\mathcal{T}$.  With the HNI method, $\vx_0 = \vm{0}$, so the radial distance
reduces to $r = \| \vx \|$. Define $\rho(t) := \| \vm{c} (t) \|$ to be the
distance from the origin to $\mathcal{C}$ at parameter value $t$.  On closer
inspection of~\eqref{eq:hni-curvedtri}, the relationship between the HNI method
and the radial distance $r$ is revealed. In~\eqref{eq:hni-curvedtri}, the only
part of the integrand not related to the geometry of the curve is $h \bigl(
\vm{c}(t) \bigr)$.  So, $h \bigl(\vm{c}(t) \bigr) := (h \circ \vm{c})(t)$ must
be the result of a one\hyp{}dimensional integral times a function of $t$. Note
that,
\begin{equation}\label{eq:hni-r1}
  h \bigl( \vm{c}(t) \bigr) 
    = \frac{\rho^2(t)}{\rho^2(t)} h \bigl( \vm{c}(t) \bigr)
    = \frac{1}{\rho^2(t)} \int_0^{\rho(t)} \frac{\partial}{\partial r} 
      \left( r^2 (h \circ \vm{\psi})(r,t) \right) \, dr .
\end{equation}
On using the product rule, we obtain
\begin{equation}\label{eq:hni-r2}
  h \bigl( \vm{c}(t) \bigr)
    = \frac{1}{\rho^2(t)} \left( 2 \int_0^{\rho(t)} r (h \circ \vm{\psi})(r, t) \, dr + 
      \int_0^{\rho(t)} r^2 \frac{\partial (h \circ \vm{\psi})(r, t)}{\partial r} \, dr
      \right) .
\end{equation}
  
\begin{figure}
  \centering
  \begin{tikzpicture}
    \node at (0in,0in) {\includegraphics[scale=1]{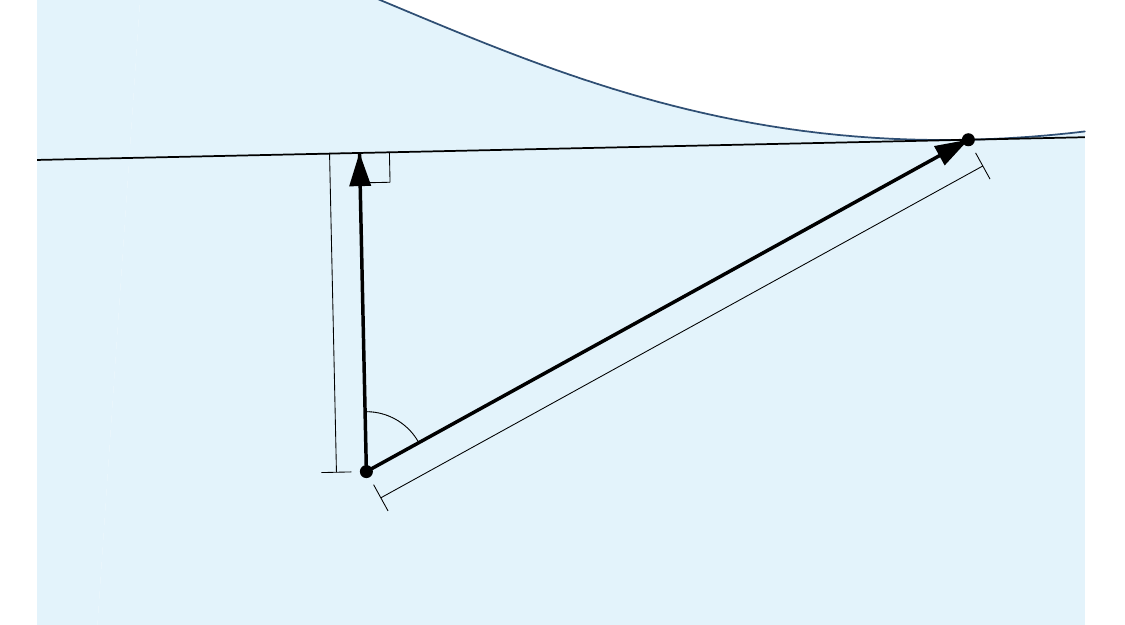}};
    \node at (1.75in,0.8in) {$\vm{c} (t)$};
    \node at (-0.63in,-0.37in) {$\alpha$};
    \node at (-0.87in,-0.75in) {$\vx_0$};
    \node at (0.65in,-0.15in) {$\rho(t)$};
    \node at (-1.475in,-0.03in) {$\dfrac{\bigl( \vm{c}(t) - \vx_0 \bigr) \cdot
      \vm{c}^{\prime\perp} (t)}{\| \vm{c}^\prime (t) \|}$};
  \end{tikzpicture}
  \caption{Interpreting quantities in~\eqref{eq:hni-radial} as the cosine of an
  angle $\alpha$.}
  \label{fig:angle-interp}
\end{figure}
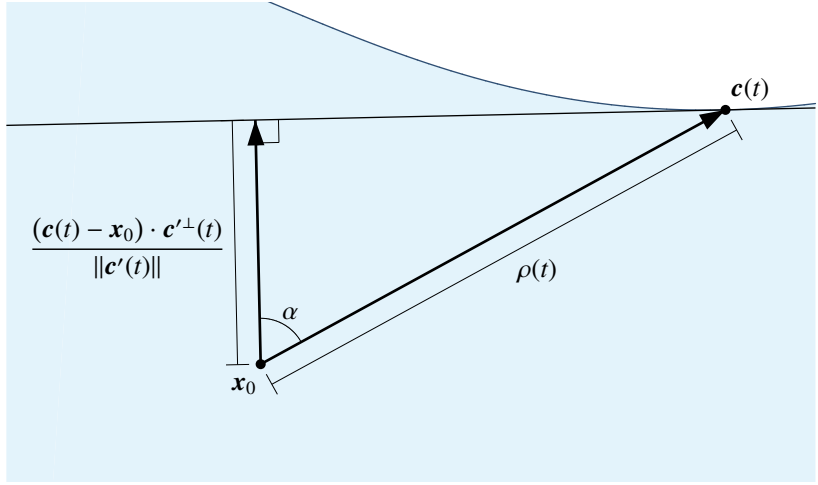

The map $\vm{\psi}$ is better understood through the lens of the polar
transformation:
\begin{subequations}\label{eq:polar-xform}
\begin{equation}
  x = r \cos \theta(t), \qquad y = r \sin \theta(t) ,
\end{equation}
and the mapping $\theta : (0, 1) \rightarrow (\theta_1, \theta_2)$,
\begin{equation}
  \theta = \arctan \left( \frac{c_2 (t)}{c_1(t)} \right) ,
\end{equation}
\end{subequations}
where $\vm{c}(t) = \bigl[ c_1 (t) \ \ c_2(t) \bigr]^T$.
From~\eqref{eq:polar-xform}, if $h(\vx)$ is a homogeneous function of degree
$q$, then $(h \circ \vm{\psi})(r, t)$ is $q$\hyp{}homogeneous in $r$ when $\vx_0
= \vm{0}$.  As a result, on applying Euler's homogeneous function theorem to the
second term on the right-hand side of~\eqref{eq:hni-r2}, we obtain
\begin{equation}\label{eq:hni-r3}
  \frac{1}{2+q} h \bigl( \vm{c}(t) \bigr)
    = \frac{1}{\rho^2(t)}\int_0^{\rho(t)} 
      r (h \circ \vm{\psi})(r, t) \, dr .
\end{equation}
Substituting~\eqref{eq:hni-r3} into~\eqref{eq:hni-curvedtri} results in
\begin{equation}\label{eq:hni-radial}
  \int_{\mathcal{T}} h (\vx) \, d\vx
    = \int_\mathcal{C} \int_0^{\rho(t)} 
      r (h \circ \vm{\psi})(r, t) \, dr \,
    \frac{\vm{c}_i(t) \cdot \vm{c}_i^{\prime\perp}(t)}
    {\rho^2 (t) \, \| \vm{c}_i^\prime (t) \|}
      \, ds .
\end{equation}
From~\eqref{eq:hni-radial}, the HNI method implicitly (and exactly) integrates a
parameter in the radial direction, provided $h(\vx)$ is a homogeneous function.
A physical interpretation of $\bigl( \vm{c}_i (t) \cdot \vm{c}^{\prime\perp} (t)
\bigr) / \bigl( \rho (t) \| \vm{c}^\prime (t) \| \bigr)$ as the cosine of an
angle $\alpha := \alpha(t)$ is illustrated in~\fref{fig:angle-interp}.
Introducing the angle $\alpha$ allows~\eqref{eq:hni-radial} to be written as
\begin{equation}\label{eq:hni-alpha}
  \int_{\mathcal{T}} h (\vx) \, d\vx
    = \int_\mathcal{C} \int_0^{\rho(t)} 
      (h \circ \vm{\psi})(r, t) \,
      \frac{r \cos \alpha}{\rho(t)} \, dr \, ds .
\end{equation}
The quantity $r \cos \alpha / \rho(t)$ can be interpreted as scaling of the
differential length $ds$.  Specifically, $r / \rho(t)$ provides a proportional
scaling of the length $ds$ based on the relative proximity of the radial
coordinate to the curve $\vm{c}(t)$.  Additionally, $\cos \alpha$ reduces $ds$
as the vector $\vm{c}(t) - \vx_0$ deviates in direction from the outward normal
vector at $\vm{c}(t)$.

\smallskip
To simplify~\eqref{eq:hni-radial}, we explicitly define
\begin{equation}\label{eq:r-to-xi}
  \xi = \frac{r}{\rho(t)}, \quad d \xi = \frac{dr}{\rho(t)} ,
\end{equation}
where $\xi \in [0, 1]$ is a scaled boundary (or a distorted polar) coordinate in
the radial direction. The point $\vx_0$ corresponds to $\xi = 0$ and along the
boundary curve $\vm{c} (t)$ we have $\xi = 1$.  Replacing~\eqref{eq:r-to-xi}
in~\eqref{eq:hni-radial} and integrating over the curve parameter $t$ leads to 
\begin{equation}\label{eq:hni-xi}
  \int_\mathcal{T} h (\vx) \, d\vx
    = \int_0^1 \int_0^1 
      (h \circ \vm{\psi})(\xi \rho(t),t)
      \, \xi \, d \xi \,
    \left[ \vm{c}_i(t) \cdot \vm{c}_i^{\prime\perp}(t) \right]
      \, d t .
\end{equation}
We define the mapped domain of integration $\Omega_0 = [0, 1]^2$ as the
parametric (reference) domain, where $(\xi, t) \in \Omega_0$ is a point in the
unit square.  Equation~\eqref{eq:hni-xi} is derived for a single curved
triangle, though the $(\xi, t)$\hyp{}parametrized integral can be applied to a
more general curved domain $\Omega$ using~\eqref{eq:hni-curve-simplified}, which
relies on the (implicit) partitioning of $\Omega$ into curved triangles. Lastly,
we mention that we can also derive~\eqref{eq:hni-xi} using the Poincar{\'e}
Lemma for vector fields~\cite{Terrell:2009:FTC} (see \ref{sec:poincare}).

\section{Scaled boundary transformation for integration}\label{sec:radial-xform}
The developments in~\sref{sec:hni} provide some insight into how the HNI method
works: by subdividing the domain into curved triangles and implicitly
integrating in the radial direction.  The integral in the radial direction is
manifested through the scaled boundary coordinate, $\xi \in [0,1]$.  However,
the derivation in~\sref{sec:hni} utilizes Euler's homogeneous function theorem,
which limits application to homogeneous functions.  Based on the insights
of \sref{sec:hni}, there must exist a transformation from the $(\xi,
t)$ parametrization to the $(x, y)$ parametrization that is not
reliant on the properties of homogeneous functions.  Herein, we show that this
mapping is in fact the SB parametrization.

\subsection{Integrating over curved \revone{planar domains} }\label{sec:sbc-curved-solids}
Mirroring the HNI method, we decompose the integral over a domain $\Omega$ into
the summation of integrals over its curved triangular partitions, $\Omega_i$ for
$i = 1, \dotsc, m$. Instead of requiring $\vx_0 = \vm{0}$, we now permit
arbitrary selection of $\vx_0$.  Judicious selection of $\vx_0$ on
star\hyp{}convex domains can result in positive cubature weights and/or improved
distribution of cubature points throughout the domain.  Recall the partition
$\Omega_i$ is bounded by two line segments and one parametric boundary curve,
$\mathcal{C}_i$. The parametric equation for $\mathcal{C}_i$ is given by
$\vm{c}_i (t)$, $t \in [0,1]$.  Next, we introduce the SB parametrization, which
transforms the domain from $\Omega^+_i \supseteq \Omega_i$ to $\Omega_0 = [0,
1]^2$.  If $\Omega_i$ is star\hyp{}convex with respect to $\vx_0$, then
$\Omega^+_i = \Omega_i$.  The SB parametrization $\vm{\varphi} : \Omega_0
\rightarrow \Omega^+_i$ is:
\begin{equation}\label{eq:sb-xform}
  \vx = \vm{\varphi}(\xi, t) = \vx_0 + \xi \bigl( \vm{c}_i(t) - \vx_0 \bigr) ,
\end{equation}
which matches the domain mapping introduced in Song~\cite{Song:1997:TSB}, though
we include the coordinate $\vx_0$ to allow shifting the center of the
parametrization away from the origin.

\smallskip
Let $J(\xi, t) = \det \bigl( \nabla \vm{\varphi} (\xi, t) \bigr)$ be the
Jacobian of the SB parametrization.  Over a curved triangle $\mathcal{T}$,
integration of a function $f : \mathcal{T} \rightarrow \Re$ is then
\begin{equation}\label{eq:2d-curved-subdomain-integral}
  \int_\mathcal{T} f(\vx)\,d\vx = \int_0^1
    \int_0^1 f (\vm{\varphi}) \, J(\xi, t) \, d\xi dt ,
\end{equation}
where $\vm{\varphi} := \vm{\varphi}(\xi, t)$.  Note that here the function $f$
is not required to be homogeneous.  The gradient of the parametrization is:
\begin{subequations}
\begin{equation}\label{eq:sb-grad}
  \nabla \vm{\varphi} (\xi, t) = \begin{bmatrix}
    \frac{\partial x_1}{\partial \xi} & \frac{\partial x_2}{\partial \xi}\\
    \frac{\partial x_1}{\partial t} & \frac{\partial x_2}{\partial t}\\
  \end{bmatrix},
\end{equation}
where
\begin{equation}
  \frac{\partial \vx}{\partial \xi} = \vm{c}(t) - \vx_0, \qquad
  \frac{\partial \vx}{\partial t} = \xi \vm{c}^\prime (t) .
\end{equation}
\end{subequations}
The Jacobian is then
\begin{equation}\label{eq:jacobian}
  J(\xi, t) = \det \bigl( \nabla \vm{\varphi} (\xi, t) \bigr) = \xi
    \det \, \Bigl[
      \vm{c}(t) - \vx_0 \quad \vm{c}_i^{\prime} (t)
    \Bigr]
    = \xi \left( \vm{c}(t) - \vx_0 \right) \cdot \vm{c}_i^{\prime\perp} (t) ,
\end{equation}
where $\vm{c}_i^{\prime\perp}(t)$ is defined in~\eqref{eq:curve-norm}.  The
vector $\vm{c}_i^{\prime\perp}(t)$ is in the outward normal direction from
$\mathcal{T}$ and $\bigl( \vm{c}(t) - \vx_0 \bigr) \cdot
\vm{c}_i^{\prime\perp}(t)$ is the signed distance of the tangent to $\vx_0$,
scaled by the inverse of the curve speed $\| \vm{c}^\prime (t) \|$ (see
\fref{fig:angle-interp}).  Equation~\eqref{eq:2d-curved-subdomain-integral} can
also be derived by combining the SB parametrization with the Poincar\'{e}
Lemma.  See~\ref{sec:poincare} for details.

\begin{figure}[t]
  \centering
  \begin{subfigure}{2in}
    \centering
    \includegraphics[scale=1]{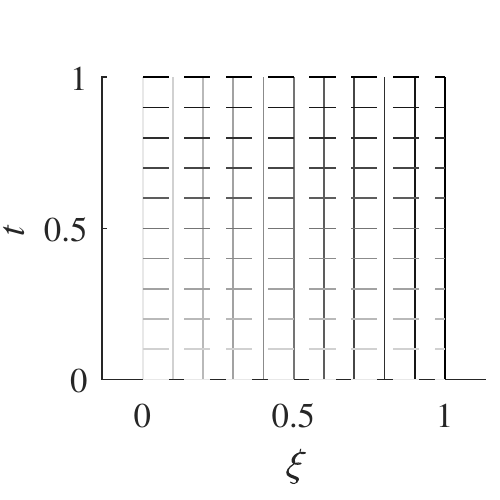}
    \caption{}\label{fig:xi-t-param}
  \end{subfigure}
  \begin{subfigure}{2in}
    \centering
    \begin{tikzpicture}
      \node at (0,0) {\includegraphics[scale=1]{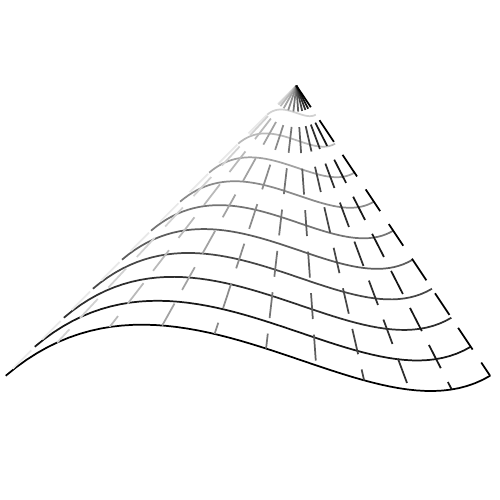}};
      \node at (0.3in,0.66in) {$\vx_0$};
    \end{tikzpicture}
    \caption{}\label{fig:xi-t-param-curvedtri}
  \end{subfigure}
  \begin{subfigure}{2in}
    \centering
    \includegraphics[scale=1]{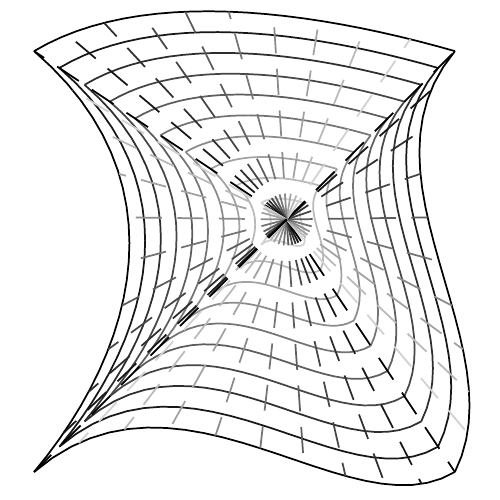}
    \caption{}\label{fig:xi-t-param-region}
  \end{subfigure}
  \caption{Isocontours of $\xi$ (solid lines) and $t$ (dashed lines) mapped from
    (a) the unit square to (b) a curved triangle.  (c) The region
    in~\fref{fig:2d-curved-region-star-convex} with $\xi$ and $t$ isocontours on
    its curved triangles.}
  \label{fig:region-param}
\end{figure}

\smallskip
An example of the SB parametrization over the region in
\fref{fig:2d-curved-region-star-convex} is illustrated in
\fref{fig:region-param}.  We note the map introduced in~\eqref{eq:sb-xform} is
not bijective since $J(0,t) = 0$, and $\xi = 0$ is a singular point.  In
\fref{fig:region-param}, $\xi = 0$ corresponds to the point $\vx_0$, which is
illustrated in \fref{fig:2d-curved-region-star-convex}. Lack of bijectivity has
ramifications when used within an isogeometric framework for the solution of
boundary\hyp{}value problems; for details and further connections, see Arioli et
al.~\cite{Arioli:2019:SBP} and the references therein.  However, when applied to
numerical integration, the SB parametrization does not pose any issues since the
edge in $\Omega_i$ that is associated with $\xi = 0$ has zero measure (area),
and hence bijectivity for $\xi > 0$ suffices.

\smallskip
To close this section, we mention the integral over $\Omega$ is simply the
summation of the integrals over each curved triangle:
\begin{equation}\label{eq:2d-curved-integral}
  \int_\Omega f(\vx)\,d\vx = \sum_{i=1}^m \int_0^1
    \int_0^1 f (\vm{\varphi}) \, \xi \, d\xi\,
    \Bigl[ \bigl( \vm{c}_i(t) - \vx_0 \bigr) \cdot 
    \vm{c}_i^{\prime\perp}(t)\Bigr]\,dt ,
\end{equation}
where the Jacobian of the transformation in~\eqref{eq:jacobian} is explicitly
introduced in~\eqref{eq:2d-curved-integral}.  Note when $\vx_0 = \vm{0}$, the
transformed integral in~\eqref{eq:2d-curved-integral} matches~\eqref{eq:hni-xi}.
Additionally, if the edges of $\Omega$ are properly ordered counterclockwise,
then the sign of the integral is correctly accounted for even if $\vx_0$ lies
outside $\Omega$ or if $\Omega$ is nonconvex.  If $\Omega_i$ is not
star\hyp{}convex, then the domain $\Omega^+_i \backslash \Omega_i$ is integrated
an even number of times, though multiplied by the scalar $-1$ half of those
times.  Finally, if $f(\vx)$ is a $m$-th degree polynomial and $\vm{c}_i(t)$ is
a $p$-th degree polynomial curve, then $f(\vm{\varphi})$ is a polynomial, and
furthermore, the integrand in~\eqref{eq:2d-curved-integral} is a polynomial of
degree $(m\!+\!1)$ in $\xi$ and degree $\bigl[(m\!+\!2)p\!-\!1\bigr]$ in $t$
that can be exactly integrated with an appropriate tensor\hyp{}product Gauss
rule.

\subsection{Integrating over polygons}\label{sec:polygons}
When the region $\Omega$ is a polygon, boundary curves $\mathcal{C}_i$ $(i = 1,
\dotsc, m)$ are reduced to line segments. One choice for parametrizing the curve
is
\begin{equation}\label{eq:line-segment-parametrization}
  \tilde{\vm{c}}_i(\tau) = \ell_i \vm{n}_i + \tau \vm{\tau}_i + \vx_0 ,
\end{equation}
where $\vm{n}_i$ is the unit outward normal vector to $\mathcal{C}_i$, $\ell_i$
is the signed distance from $\vx_0$ to the line on which the line segment
$\mathcal{C}_i$ lies, $\vm{\tau}_i = \vm{e}_3 \times \vm{n}_i$ is a unit vector
tangent to $\mathcal{C}_i$, and $(\tau_1)_i < \tau < (\tau_2)_i$ is a parameter.
The interval of $\tau$ is given by $(\tau_1)_i = \sqrt{\| (\vx_1)_i - \vx_0 \|^2
- \ell_i^2}$ and $(\tau_2)_i = \sqrt{\| (\vx_2)_i - \vx_0 \|^2 - \ell_i^2}$,
where $(\vx_1)_i$ and $(\vx_2)_i$ are the coordinates of the endpoints of
$\mathcal{C}_i$.  Rather than $t \in [0,1]$ as suggested
in~\sref{sec:hni-curved-triangles}, a different domain of definition for $t$ is
used because it simplifies integral transforms for weakly singular
functions that appear in~\sref{sec:singular-t}. To map the domain to the unit
square, the affine mapping
\begin{equation}\label{eq:tau}
  \tau = (\tau_1)_i \, (1-t) + (\tau_2)_i \, t, \quad
  d \tau = \left( (\tau_2)_i - (\tau_1)_i \right) dt = \Delta \tau_i \, dt ,
\end{equation}
is applied. In~\eqref{eq:tau}, $\Delta \tau_i = (\tau_2)_i - (\tau_1)_i$ is the
length of $\mathcal{C}_i$.  Given the boundary curve $\tilde{\vm{c}}_i (\tau)$,
\eqref{eq:2d-curved-subdomain-integral} simplifies to
\begin{equation}\label{eq:polygon-subdomain-integral}
  \int_{\Omega_i} f(\vx)\,d\vx
    = \ell_i \Delta \tau_i 
      \int_0^1 \int_0^1 f(\vm{\varphi}) \, \xi \, d\xi \, dt .
\end{equation}
Integration over the polygonal domain $\Omega$ is then given by
\begin{equation}\label{eq:polygon-integral}
  \int_\Omega f(\vx)\,d\vx 
    = \sum_{i=1}^m \ell_i \Delta \tau_i
      \int_0^1 \int_0^1 f(\vm{\varphi}) \, \xi \, d\xi \, dt .
\end{equation}
For polygonal domains, selecting $\vx_0$ as a vertex eliminates the need to
perform cubature on the two edges that contain the vertex. This is apparent
in~\eqref{eq:polygon-integral} since $\ell_i = 0$ for the two edges that pass
through $\vx_0$.  However, it is worth noting that the influence of $\vx_0$
being at a vertex on the number of cubature points reduces with increasing
numbers of edges in the polygon.  For domains that are bounded by both affine
and non\hyp{}affine curves, \eqref{eq:polygon-subdomain-integral}
and~\eqref{eq:2d-curved-subdomain-integral} can be utilized as needed to
simplify integration over the domain.

\smallskip
When $f$ is a $p$\hyp{}th degree polynomial, it is a sum of monomials that have
the form $x^\alpha y^\beta$, where $\alpha, \beta \in \mathbb{Z}^\geq$, $\alpha
+ \beta \leq p$, and $\mathbb{Z}^\geq$ is the set of integers greater than or
equal to zero.  The monomial can be expressed in $(\xi, t)$-coordinates through
the mapping of the domain introduced in~\eqref{eq:sb-xform}.  We assume
$\vm{e}_1 = \vm{n}$, $\vm{e}_2 = \vm{\tau}$, and $\vx_0 = \vm{0}$
in~\eqref{eq:line-segment-parametrization}, where $\vm{e}_1$ and $\vm{e}_2$ are
unit vectors in the $x$- and $y$-directions, respectively.  These values for
$\vm{n}$, $\vm{\tau}$, and $\vx_0$ can be obtained through a translational and
rotational mapping of the $(x,y)$\hyp{}coordinates.  When this mapping is
applied to $f(\vx)$, it does not change the polynomial nature of $f(\vx)$ or
polynomial degree of $f(\vx)$.  The transformed variables are then
\begin{equation}\label{eq:aligned-triangle-xform}
  x = \xi \ell_i, \quad y = \xi \tau.
\end{equation}
Inserting $f = x^\alpha y^\beta$ in~\eqref{eq:polygon-integral} and
substituting the transformed variables gives
\begin{equation}\label{eq:polygon-polynomial-integral}
  \int_\Omega f(\vx) \, d\vx 
    = \sum_{i=1}^m \ell_i^{\alpha + 1} \Delta \tau_i
      \int_0^1 \int_0^1 \xi^{\alpha + \beta + 1} 
      \left[ (\tau_1)_i \, (1-t) + (\tau_2)_i \, t \right]^\beta \, d\xi \, dt.
\end{equation}
Based on~\eqref{eq:polygon-polynomial-integral}, a cubature rule capable of
exactly integrating polynomials of degree $p + 1$ in the 
$\xi$ direction and
degree $p$ in the $t$ direction over each triangle $\Omega_i$ is capable of
exactly integrating a polynomial of degree $p$ over a polygon
using~\eqref{eq:polygon-integral}.

\smallskip
As will be explored in~\sref{sec:examples}, cubature over polygons using this
rule requires more integration points than a similar triangulation paired with
an optimized, symmetric cubature rule; however, there are several advantages to
the SBC approach.  First, the tensor\hyp{}product structure of the SBC rule
allows it to be defined with reduced memory requirements, which may be
beneficial in applications where memory is limited and data movement is
computationally expensive (such as in GPU programming). Also, the triangulation
of the domain is trivial with the SBC method and does not change with nonconvex
polygons. Finally, Gauss rules of arbitrary order can be easily generated to
compute~\eqref{eq:polygon-integral}, while high\hyp{}order symmetric rules
over triangles, if not available on\hyp{}hand, are time consuming to generate.

\smallskip
Finally, we point out that if $h (\vx)$ is a homogeneous function of degree $q$,
then on applying~\eqref{eq:hni-r3} and~\eqref{eq:r-to-xi}
to~\eqref{eq:polygon-integral}, we obtain
\begin{equation}\label{eq:polygon-hni}
  \int_\Omega h (\vx)\,d\vx = \frac{1}{2+q}
    \sum_{i=1}^m \ell_i
    \int_{\Omega_i} h \bigl( \vx \bigr) \, d s .
\end{equation}
Noting that $\ell_i = b_i / \|\vm{a}_i\|$, where $\vm{a}_i \cdot \vx = b_i$ is
the equation of the line on which $\mathcal{C}_i$ lies, \eqref{eq:polygon-hni}
matches the result in Chin et al.~\cite{Chin:2015:NIH}.  In Chin et
al.~\cite{Chin:2015:NIH}, integration over polynomials using the HNI method is
presented.  To integrate a $p$-th degree polynomial, the polynomial must first
be decomposed into $(p\!+\!1)$ homogeneous polynomials.  
Then, \eqref{eq:polygon-hni} is applied to each homogeneous term.  With the HNI
method, $(p\!+\!1)^2/2$ evaluations per edge are required to exactly integrate
$p$-th degree polynomials using a Gauss\hyp{}Legendre rule.  In comparison,
integrating a $p$-th degree polynomial using~\eqref{eq:polygon-integral}
requires $(p\!+\!2)(p\!+\!1)/4$ Gauss\hyp{}Legendre points per edge.

\section{Integration of weakly singular and nearly singular
functions}\label{sec:singular}
\revtwo{In this section, we introduce methods for handling integration of weakly
singular functions.  These methods modify the SB mapping from the unit square to
the curved triangle to remove or smooth singularities in the integrand.
\sref{sec:singular-xi} introduces an additional parameter in the SB mapping that
eliminates rational singularities in the $\xi$ direction.
\sref{sec:singular-t} introduces several mappings that smooth near
singularities in the $t$ direction of the SB map.}

\subsection{Cancelling weak singularities in the radial
direction}\label{sec:singular-xi}
\revtwo{Given a function $g(\vx) \in C^0(\Omega)$,  with a leading term that is 
a constant,} the integral
\begin{equation}\label{eq:singular-curved}
  I := \int_\Omega \frac{g(\vx)}{\| \vx - \vx_c \|^\beta} \, d\vx 
\end{equation}
for $0 < \beta < 2$ is said to be \textit{weakly singular} when $\Omega \subset
\Re^2$ and $\vx_c \in \revone{\bar \Omega}$.  \revtwo{For a weakly singular integral, $|I| <
\infty$.}  In~\eqref{eq:singular-curved}, $\vx_c$ serves as the location of a
point singularity.  When $\vx_c \in \revone{\bar \Omega}$ or when $\vx_c$ is near \revone{$\partial \Omega$},
standard polynomial\hyp{}precision cubature rules perform poorly when used to
integrate~\eqref{eq:singular-curved}.  To improve cubature convergence, special
integration techniques are usually employed.  One such technique is the Duffy
transformation~\cite{Duffy:1982:QOP}, which introduces the change of variables,
\begin{equation}
  x = u, \qquad y = xv = uv,
\end{equation}
which maps the standard triangle $T_0 := [0,1] \times [0,x]$ to the unit square
$\Omega_0 = [0,1]^2$.  This mapping eliminates singularities when $\beta = 1$
in~\eqref{eq:singular-curved}.  On $T_0$, the SB parametrization is
the same as the Duffy transformation (see~\eqref{eq:aligned-triangle-xform} and
note $\ell_i = 1$), though the SB parametrization is not limited to this
triangular domain.  Its ability to eliminate singularities when $\beta = 1$ can
be observed by applying the SB parametrization~\eqref{eq:sb-xform} with $\vx_0 =
\vx_c$ to~\eqref{eq:singular-curved}:
\begin{equation}\label{eq:duffy-sb-xform}
  \int_\Omega \frac{g(\vx)}{\| \vx - \vx_c \|}\,d\vx 
    = \sum_{i=1}^m \int_0^1
    \int_0^1 g (\vm{\varphi}) \, d\xi\,
    \frac{\bigl( \vm{c}_i(t) - \vx_c \bigr) \cdot 
    \vm{c}_i^{\prime\perp}(t)}
      {\| \vm{c}_i(t) - \vx_c \|}\,dt .
\end{equation}
In~\eqref{eq:duffy-sb-xform}, the singularity in the radial direction is
eliminated, thereby leaving a smooth integrand.

\begin{figure}[t]
  \centering
  \begin{subfigure}{2in}
    \centering
    \includegraphics[scale=1]{xi-t-param}
    \caption{}\label{fig:xi-t-param-a2}
  \end{subfigure}
  \begin{subfigure}{2in}
    \centering
    \includegraphics[scale=1]{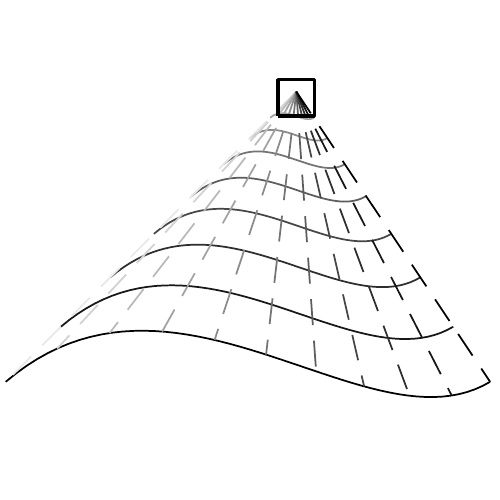}
    \caption{}\label{fig:xi-t-param-curvedtri-a2}
  \end{subfigure}
  \begin{subfigure}{2in}
    \centering
    \begin{tikzpicture}
      \node at (0,0) {\includegraphics[scale=1]{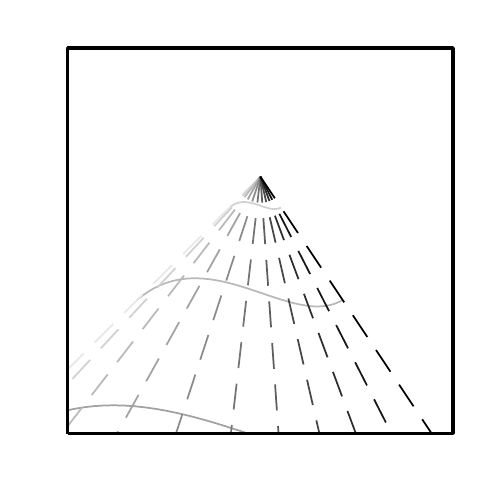}};
      \node at (0.12in,0.35in) {$\vx_0$};
    \end{tikzpicture}
    \caption{}\label{fig:xi-t-param-curvedtri-a2-zoom}
  \end{subfigure}
  \begin{subfigure}{2in}
    \centering
    \includegraphics[scale=1]{xi-t-param}
    \caption{}\label{fig:xi-t-param-a3}
  \end{subfigure}
  \begin{subfigure}{2in}
    \centering
    \includegraphics[scale=1]{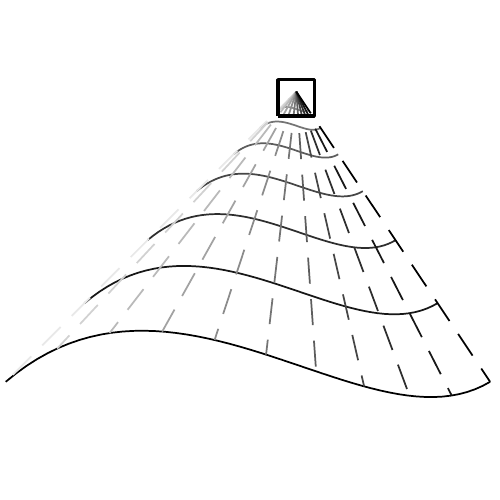}
    \caption{}\label{fig:xi-t-param-curvedtri-a3}
  \end{subfigure}
  \begin{subfigure}{2in}
    \centering
    \begin{tikzpicture}
      \node at (0,0) {\includegraphics[scale=1]{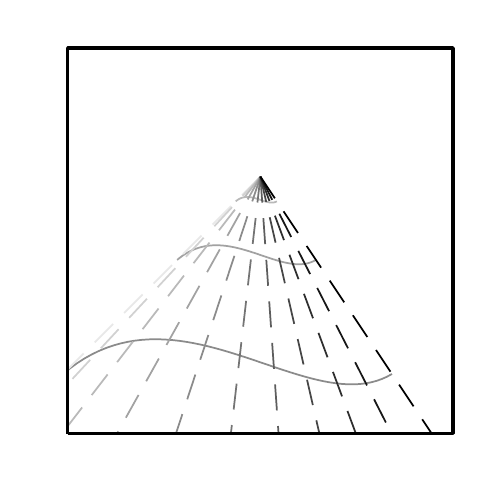}};
      \node at (0.12in,0.35in) {$\vx_0$};
    \end{tikzpicture}
    \caption{}\label{fig:xi-t-param-curvedtri-a3-zoom}
  \end{subfigure}
  \caption{Isocontours of $\xi$ (solid lines) and $t$ (dashed lines) mapped from
    (a), (d) the unit square to (b), (e) a curved triangle with the generalized
    SB parametrization.  (c), (f) \revone{Zoomed region illustrating more $\xi$
    isocontours near $\vx_0$ as $\alpha$ increases}.  In (a)--(c), $\alpha = 2$
    and in (d)--(f), $\alpha = 3$.  \revtwo{As $\alpha$ increases, the
    isocontours of $\xi$ shift toward $\vx_0$.}}
  \label{fig:region-gsb-param}
\end{figure}

\smallskip
More generally, when $0 < \beta < 2$, the generalized Duffy
transformation~\cite{Mousavi:2010:GDT} enables polynomial\hyp{}precision
integration to be recovered for $\Omega = T_0$ in~\eqref{eq:duffy-sb-xform}.
The transformation is:
\begin{equation}\label{eq:gd-xform}
  x = u^\alpha , \qquad y = xv = u^\alpha v.
\end{equation}
In~\eqref{eq:gd-xform}, $\alpha$ is selected such that $\alpha \in \mathbb{Z}_+$
and $\alpha(2 - \beta) \in \mathbb{Z}_+$.  Mimicking the generalized Duffy
transformation, a more general mapping on $\Omega_i$ is invoked:
\begin{equation}\label{eq:gsb-xform}
  \vx = \vm{\varphi}_\alpha (\xi, t) =
    \vx_0 + \xi^\alpha \bigl( \vm{c}_i (t) - \vx_0 \bigr) 
\end{equation}
for $\alpha > 0$.  We will refer to this as a generalized SB parametrization.
The isocontours of $\xi$ and $t$ with the generalized SB parametrization for
$\alpha = 2$ and $\alpha = 3$ are presented in~\fref{fig:region-gsb-param}.  As
$\alpha$ increases, the number of $\xi$\hyp{}isocontours in the vicinity of
$\vx_c$ increases.  In the context of a cubature rule, this has the practical
effect of placing more integration points closer to the singularity as $\alpha$
gets larger.  For the triangle $T_0$, the generalized SB transformation is
identical to the generalized Duffy transformation.  The map for $\vm{c}_i(t)$ in
$T_0$ is:
\begin{subequations}
\begin{equation}
  \vm{c}_i (t) = \bigl[ 1 \quad t \bigr]^T,
\end{equation}
so the generalized SB map is
\begin{equation}
  x = \xi^\alpha , \qquad y = xv = \xi^\alpha t ,
\end{equation}
\end{subequations}
which matches the generalized Duffy map in~\eqref{eq:gd-xform}.

\smallskip
The gradient of the generalized SB transformation is given
by~\eqref{eq:sb-grad}, where now
\begin{equation}
  \frac{\partial \vx}{\partial \xi} = \alpha \xi^{\alpha-1}
    \left( \vm{c}(t) - \vx_0 \right), \qquad
  \frac{\partial \vx}{\partial t} = \xi^\alpha \vm{c}^\prime (t) .
\end{equation}
The Jacobian of the mapping is:
\begin{equation}\label{eq:det-singular-xform}
  J(\xi, t) = \det \bigl( \nabla \vm{\varphi}_\alpha (\xi, t) \bigr) 
    = \alpha \xi^{2 \alpha - 1} \bigl( \vm{c}_i (t) - \vx_0 \bigr)
    \cdot \vm{c}_i^{\prime\perp}(t) .
\end{equation}
After applying the generalized SB transformation~\eqref{eq:gsb-xform}
to~\eqref{eq:singular-curved} and selecting $\vx_0 = \vx_c$, we obtain
\begin{equation}\label{eq:gsb-integral}
  \int_\Omega \frac{g(\vx)}{\| \vx - \vx_c \|^\beta}\,d\vx 
    = \alpha \sum_{i=1}^m 
    \int_0^1 \int_0^1 g (\vm{\varphi}_\alpha) \xi^{\alpha (2 - \beta) - 1} \, d\xi
    \frac{\bigl( \vm{c}_i (t) - \vx_c \bigr)
      \cdot \vm{c}_i^{\prime\perp}(t)}
      {\| \vm{c}_i (t) - \vx_c \|^\beta} \, d t ,
\end{equation}
where $\vm{\varphi}_\alpha := \vm{\varphi}_\alpha(\xi, t)$.  On selecting
$\alpha$ such that $\alpha (2 - \beta) \in \mathbb{Z}_+$, the
$r^{-\beta}$\hyp{}singularity in~\eqref{eq:singular-curved} is eliminated.
However, depending on the choice for $\alpha$ and the form of
$g(\vm{\varphi}_\alpha)$, integration convergence may still be poor.  For
example, if $g(\vx)$ is a polynomial function, selecting $\alpha \notin
\mathbb{Z}_+$ results in $g(\vm{\varphi}_\alpha)$ being non\hyp{}polynomial,
which potentially increases cubature error and reduces convergence rates when
using a Gauss rule in the $\xi$-direction. To ensure $g(\vm{\varphi}_\alpha)$
remains polynomial, we select the smallest $\alpha \in \mathbb{Z}_\beta$, where
\begin{equation}
  \mathbb{Z}_\beta = \Bigl\{\; \alpha \in \mathbb{Z}_+ \;:\; 
    \alpha (2 - \beta) \in \mathbb{Z}_+ \;\Bigr\} .
\end{equation}
The same criteria for selecting $\alpha$ is used in the generalized Duffy
transformation~\cite{Mousavi:2010:GDT}.

\smallskip
The generalized SB parametrization eliminates fractional singularities in an
integrand and recovers polynomial\hyp{}precision integration for the
non\hyp{}singular portion of an integrand.  However, the degree of the
transformed polynomials can be significantly higher than the untransformed
integrand.  For example, if $\beta = 1 / \gamma$ for $\gamma \in \mathbb{Z}_+$,
then $\alpha = \gamma$ is the best selection in $\mathbb{Z}_\beta$.  If $g(\vx)$
is a $p$\hyp{}th degree polynomial, the $\xi$ integrand after applying the
generalized SB mapping is a polynomial of degree $\gamma
(p\!+\!2\!-\!\beta)\!-\!1 = \mathcal{O}(\gamma p)$.  While the increase in
polynomial order may be a worthwhile tradeoff for improved handling of the
singularity for small values of $\gamma$, we note the polynomial order of
$g(\vm{\varphi}_\alpha)$ quickly increases as $\gamma$ grows.  Rather than
resorting to a high\hyp{}order Gauss\hyp{}Legendre rule paired with the
generalized SB transformation, we find that Gauss\hyp{}Jacobi
quadrature~\cite{Chernov:2012:ECG} in the 
$\xi$ direction provides a much
more efficient quadrature rule for integrands when $\gamma$ is large.

\smallskip
One important distinction between the generalized Duffy transformation and the
generalized SB parametrization is the applicability of the SB parametrization to
arbitrary convex and nonconvex polygons and regions bounded by curves.  This
allows the generalized SB parametrization to be applied directly to these
shapes, whereas the generalized Duffy transformation first requires a
partitioning of $\Omega$ and a mapping from each subdomain to $T_0$. This
mapping, if not chosen appropriately, can distort the singularity. In addition,
if the boundaries of the partitions are not affine, the mapping may be difficult
to compute.  Also, with the SB parametrization, integration over regions
where $\vx_c \notin \revone{\bar \Omega}$ is handled intrinsically (see \sref{sec:hni} for
details), which allows the transformation to be easily applied to integrate over
regions that contain nearby singularities. Integration of functions with
singularities very close to the region of integration
(\textit{nearly singular} integrands) have poor cubature convergence when
integrated with Gauss rules~\cite{Chin:2017:MCD}.

\subsection{Cancelling near\hyp{}singularities
            in the edge direction of polygons}\label{sec:singular-t}
When integrands contain singularities 
in the vicinity of $\revone{\partial \Omega}$, integrands in
the $t$ direction can be nearly singular.  While these integrands remain smooth
and differentiable, they are poorly approximated by polynomials and standard
Gauss quadrature rules do not \revone{lead to efficient evaluation of these
integrals}. When $\Omega$ is a polygon, the boundary curves $\mathcal{C}_i$ are
line segments, which simplifies integration using the SB parametrization.
Furthermore, this simplified representation allows the formulation of integral
transformations that eliminate from the integrand near\hyp{}singularities that
are proximal to $\mathcal{C}_i$. This idea was presented in Ma and
Kamiya~\cite{Ma:2002:DTN}, where it was shown to improve integration in the BEM.
More recently, Lv et al.~\cite{Lv:2019:ASD} demonstrated the applicability of
the transformation in combination with the generalized Duffy transformation.  In
this section, we show the transformation is also relevant here, and furthermore,
we introduce two more transformations that cancel other types of
near\hyp{}singularities.

\smallskip
When the generalized SB parametrization is applied to the integral
in~\eqref{eq:singular-curved}, the result is~\eqref{eq:gsb-integral}.  However,
when $\Omega$ is a polygon, the curve parametrization simplifies
to~\eqref{eq:line-segment-parametrization} with $\vx_0 = \vx_c$:
\begin{equation}\label{eq:line-segment-singular-xform}
  \tilde{\vm{c}}_i(\tau) = \ell_i \vm{n}_i + \tau \vm{\tau}_i + \vx_c .
\end{equation}
Given this parametrization, the distance to the singularity simplifies to
\begin{equation}\label{eq:line-segment-distance}
  r(\tau) = \|\tilde{\vm{c}}_i(\tau) - \vx_c\| = \sqrt{\ell_i^2 + \tau^2}.
\end{equation}
Inserting~\eqref{eq:line-segment-distance} into~\eqref{eq:gsb-integral} and
parametrizing the boundary curves by $\tau$ \revone{results in}
\begin{equation}\label{eq:polygon-gsb-xform}
  \int_\Omega \frac{g(\vx)}{\| \vx - \vx_c \|^\beta}\,d\vx 
    = \alpha \sum_{i=1}^m \ell_i \int_{(\tau_1)_i}^{(\tau_2)_i} \int_0^1 
    g (\vm{\varphi}_a) \, \xi^{\alpha (2 - \beta) - 1} \, d\xi
    \frac{1}{(\ell_i^2 + \tau^2)^{\beta/2}} \, d \tau .
\end{equation}

\smallskip
When $\beta = 1$ (and $\alpha = 1$), the SB parametrization is recovered and
$1/r$ singularities are removed in the $\xi$ direction.  To remove the
near\hyp{}singularity in the $\tau$ direction caused by the radial
dependence of the integrand, we define the transformation
\begin{equation}\label{eq:2d-transformation}
  d\tilde{\tau} = \dfrac{1}{\sqrt{\ell_i^2 + \tau^2}} \, d\tau .
\end{equation}
Integrating~\eqref{eq:2d-transformation} gives
\begin{subequations}
\begin{equation}
  \tilde{\tau} (\tau) = \ln \left( \tau + \sqrt{\ell_i^2 + \tau^2} \right) ,
\end{equation}
and solving for $t$ results in~\cite{Ma:2002:DTN,Lv:2019:ASD}
\begin{equation}\label{eq:2d-transformation-sub}
  \tau (\tilde{\tau}) 
    = \frac{1}{2} e^{-\tilde{\tau}} \left( e^{2 \tilde{\tau}} - \ell_i^2 \right).
\end{equation}
\end{subequations}
Applying the transform to~\eqref{eq:polygon-gsb-xform} (with $\beta = \alpha =
1$) \revone{yields} 
\begin{equation} \label{eq:polygon-1/r-xform}
  \int_\Omega \frac{g(\vx)}{\| \vx - \vx_c \|^\beta}\,d\vx 
    = \sum_{i=1}^m \ell_i \int_{\ln \left( (\tau_1)_i + 
    \sqrt{\ell_i^2 + (\tau_1)_i^2} \right)}^{\ln 
    \left( (\tau_2)_i + \sqrt{\ell_i^2 + (\tau_2)_i^2} \right)} \int_0^1 
    g (\vm{\varphi})) \, d\xi d \tilde{\tau} ,
\end{equation}
where the near\hyp{}singularity has been eliminated from the integral.  
Transforming integration to the unit square $\Omega_0$, we obtain:
\begin{equation} \label{eq:polygon-1/r-xform2}
  \int_\Omega \frac{g(\vx)}{\| \vx - \vx_c \|^\beta}\,d\vx 
    = \sum_{i=1}^m \ell_i \ln \left( \frac{(\tau_2)_i + 
    \sqrt{\ell_i^2 + (\tau_2)_i^2}}{(\tau_1)_i + 
    \sqrt{\ell_i^2 + (\tau_1)_i^2}} \right) 
    \int_0^1 \int_0^1 g (\vm{\varphi}) \, d\xi d t .
\end{equation}
The SB parametrization paired with the near\hyp{}singularity\hyp{}cancelling
transformation with $\beta = 1$ completely eliminates $1/r$ singularities from
the integrand.  However, even when the singularity is not of the type $1/r$, the
$\beta = 1$ integral transformation can still reduce errors in integration by
smoothing the near\hyp{}singularity.  This is explored further
in~\sref{sec:ex-singular-t}.

\smallskip
\revone{Following a procedure similar to that above, integral transformations
can be formulated for integrals with near-singularities that are of the
form~\eqref{eq:singular-curved} with $\beta = 2$ and $\beta = 3$.  While $| I |
= \infty$ in~\eqref{eq:singular-curved} when $\beta \geq 2$ and $\vx_c \in
\bar{\Omega}$, the transformations are designed to eliminate the
near\hyp{}singularity from $(\ell^2_i + \tau^2)^{-\beta/2}$
in~\eqref{eq:polygon-gsb-xform}, which remains bounded on all edges
$\mathcal{C}_i$ where $|T_i| \neq 0$ (that is, $\ell_i \neq 0$). These two
transformations, along with the transformation for $\beta = 1$, are listed
in~\tref{t:2d-integral-transformations}.} The performance of all three
transformations in integrating \revone{weak singularities} is examined
in~\sref{sec:ex-singular-t}.  Further examples illustrating the performance of
these transformations with the HNI method and singular integrands are available
in Chin~\cite{Chin:2019:NIH}.  Chin et al.~\cite{Chin:2015:NIH} presents an
alternative derivation for the integral transformation to cancel singularities
of the type $ 1 / r^2$ when $h(\vx)$ is a homogeneous function. This derivation
relies on parametrization of~\eqref{eq:polygon-hni} in terms of $\theta$.  The
transformation in Chin et al.~\cite{Chin:2015:NIH} is equivalent to the
transformation in~\tref{t:2d-integral-transformations}.  This is apparent when
noting the mapping to $\tilde{\tau}$ in~\tref{t:2d-integral-transformations} is
equivalent to the mapping to $\theta$, namely $\theta = \tan^{-1} (y/x)$, where
without any loss of generality, we choose $\vm{n} = \vm{e}_1$ and $\vm{t} =
\vm{e}_2$.

{\renewcommand{\arraystretch}{2.5}
\begin{table}[t]
  \centering
  \caption{Integral transformations in $\Re^2$ to cancel near\hyp{}singularities
    in~\eqref{eq:polygon-gsb-xform}.}%
  \label{t:2d-integral-transformations}
  \begin{tabular}{| c | c | c | c |}
    \hline
    $f(\vx)$              & $d \tilde{\tau}$
      & $\tilde{\tau}(\tau)$ 
      & $t(\tilde{\tau})$ \\
    \hline
    $\dfrac{g(\vx)}{r}$   & $\dfrac{d\tau}{\sqrt{\ell_i^2 + \tau^2}}$
      & $\ln \left( \tau + \sqrt{\ell_i^2 + \tau^2} \right)$
      & $\dfrac{1}{2} e^{-\tilde{\tau}} \left( e^{2 \tilde{\tau}} - \ell_i^2 \right)$ \\[1.5em]
    \hline
    $\dfrac{g(\vx)}{r^2}$ & $\dfrac{\ell_i d\tau}{\ell_i^2 + \tau^2}$
      & $\tan^{-1} \dfrac{\tau}{\ell_i}$
      & $\ell_i \tan \tilde{\tau}$ \\[0.75em]
    \hline
    $\dfrac{g(\vx)}{r^3}$ & $\dfrac{\ell_i^2 d\tau}{\bigl( \ell_i^2 + \tau^2 \bigr)^{3/2}}$
      & $\dfrac{\tau}{\sqrt{\ell_i^2 + \tau^2}}$
      & $\pm \dfrac{\ell_i \tilde{\tau}}{\sqrt{1 - \tilde{\tau}^2}}$ \\[1.5em]
    \hline
  \end{tabular}
\end{table}
}

\section{Numerical examples}\label{sec:examples}

\revtwo{This section evaluates the performance of the SBC method in numerically
integrating various functions over different types of domains. Integration of
polynomial and non\hyp{}polynomial functions over polygons is presented in
\sref{sec:ex-polygonal}.  We introduce comparisons to other methods of
integration based on both accuracy per cubature point and time to generate an
integration rule.  \sref{sec:ex-curved} focuses on curved domains, demonstrating
efficient integration of test functions over regions bounded by both polynomial
and rational curves.  In \sref{sec:ex-singular}, the examples center on the
integration of weakly singular functions using the methods introduced in
\sref{sec:singular}.  These methods are compared and a 
test case that arises in the X-FEM is studied. 
Finally, \sref{sec:ex-tmvi} highlights how the SBC method can
be used in conjunction with TMVI to compute $L_2$ error.  We also point out a
link between TMVI and homogeneous functions, which mirrors the connection
between the SBC and HNI methods.}

\subsection{Polygonal domains}\label{sec:ex-polygonal}

\subsubsection{Polynomial integration}
\begin{figure}[t]
  \centering
  \begin{subfigure}{3in}
    \centering
    \includegraphics[scale=1]{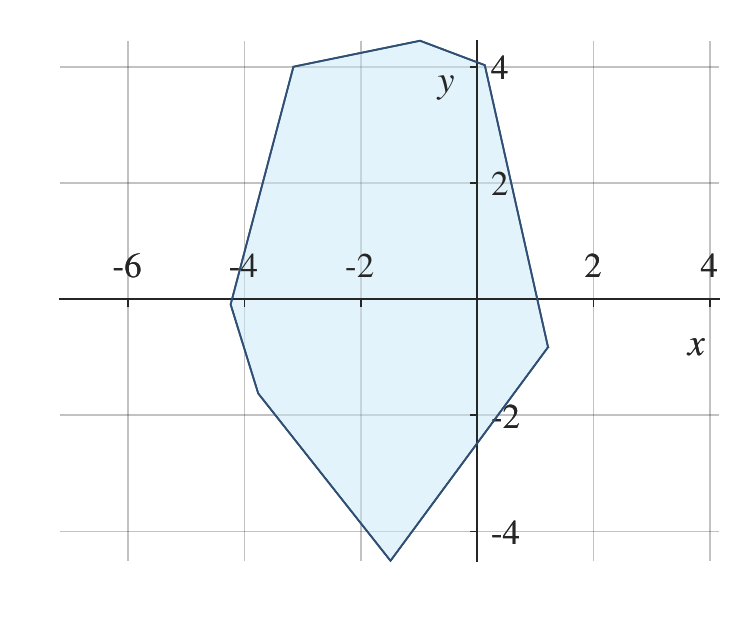}
    \caption{}\label{fig:polygon-convex-1}
  \end{subfigure}
  \begin{subfigure}{3in}
    \centering
    \includegraphics[scale=1]{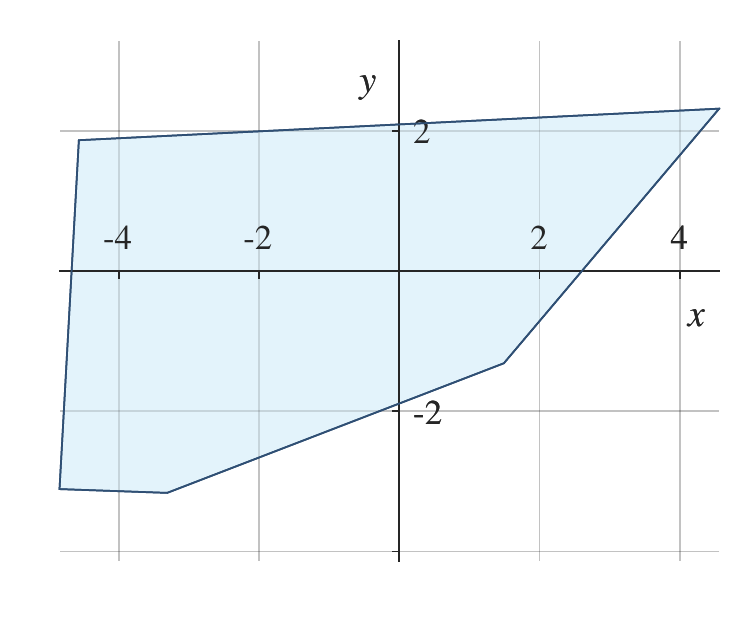}
    \caption{}\label{fig:polygon-convex-2}
  \end{subfigure}
  \begin{subfigure}{3in}
    \centering
    \includegraphics[scale=1]{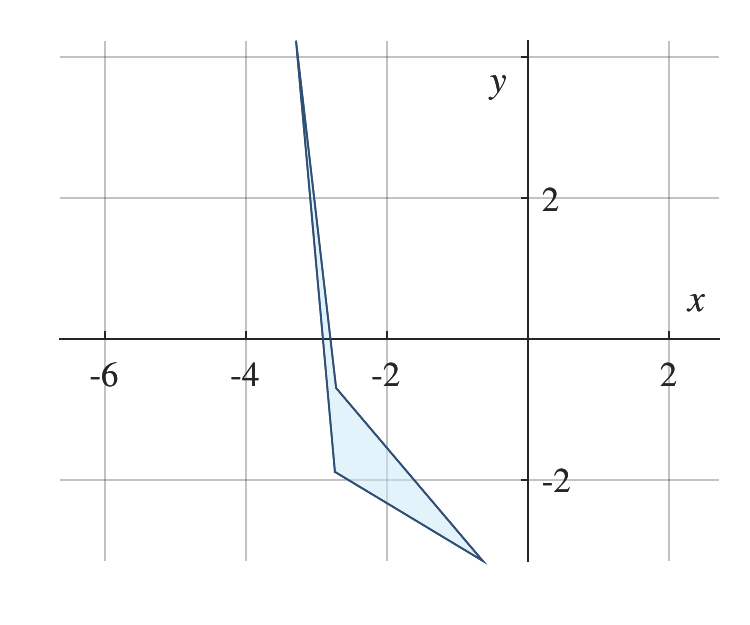}
    \caption{}\label{fig:polygon-nonconvex-1}
  \end{subfigure}
  \begin{subfigure}{3in}
    \centering
    \includegraphics[scale=1]{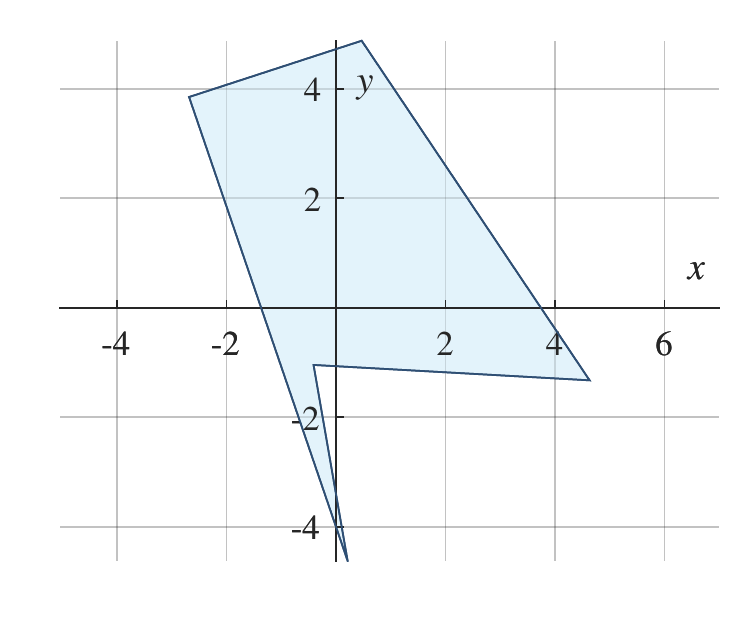}
    \caption{}\label{fig:polygon-nonconvex-2}
  \end{subfigure}
  \caption{Convex and nonconvex polygons over which polynomials in
  \tref{tab:ex-polynomials} are integrated.}
  \label{fig:ex-polynomials-polygons}
\end{figure}

As discussed in~\sref{sec:polygons}, the SBC scheme is capable of exactly
integrating polynomials of degree $p$ over polygons when provided a
polynomial\hyp{}precise quadrature rule that is capable of integrating
polynomials of degree $p+1$ in the $\xi$ direction and degree $p$ in the
$t$ direction.  This corresponds to a Gauss\hyp{}Legendre rule of greater
than or equal to $(p + 2)/2$ $\xi$\hyp{}points and $(p + 1)/2$ $t$\hyp{}points
for a $p$\hyp{}th degree polynomial integrand.  Beyond requiring a closed
polygonal region, no restrictions are placed on the shape of the polygon.
Therefore, integration using the SBC method is possible over both convex and
nonconvex polygons.  We verify this approach on several convex and nonconvex
polygons illustrated in~\fref{fig:ex-polynomials-polygons}.  The coordinates of
the vertices of these polygons are listed in Chin et al.~\cite{Chin:2015:NIH}.
Over each polygon, polynomials ranging from degree $p = 0$ to $p = 5$ are
integrated. These polynomials are listed in~\tref{tab:ex-polynomials}.  Exact
values of the integrals are computed using region integration with
\texttt{Integrate[]} in Mathematica 12.0.0.

{\renewcommand{\arraystretch}{1.4}
\begin{table}[t]
  \centering
  \caption{Polynomials integrated over the polygons illustrated in~\fref{fig:ex-polynomials-polygons}.}
  \label{tab:ex-polynomials}
  \begin{tabular}{|c|c|}
    \hline
    $p$ & $f(x,y)$ \\
    \hline
    $0$ & $1$ \\
    \hline
    $1$ & $x - 2y + 1$ \\
    \hline
    $2$ & $3x^2 + 4xy - 2y^2 - x + 2y - 3$ \\
    \hline
    $3$ & $4x^3 - 2x^2y - 3xy^2 + y^3 + 8x^2 - 4xy + 5y^2 - 6x - 4y + 7$ \\
    \hline
    $4$ & $-3x^4 - 5x^3y +2x^2y^2 - 9xy^3 + y^4 + 3x^3 - 2x^2y - xy^2 + 5y^3 + 4x^2 - 7xy - 6y^2 - 4x + 6y - 8$ \\
    \hline
    \multirow{2}{*}{$5$} & $10x^5 - 5x^4y - 7x^3y^2 + 6x^2y^3 + 3xy^4 + y^5 - x^4 + 2x^3y + 11x^2y^2 - 8xy^3 - 2y^4 - 3x^3 + 9x^2y$ \\[-0.05in]
        & $+ 8xy^2 - 10y^3 - 9x^2 - 6xy + 7y^2 + 5x - 4y + 4$ \\
    \hline
  \end{tabular}
\end{table}
}

\begin{figure}[t]
  \centering
  \begin{subfigure}{3in}
    \centering
    \includegraphics[scale=1]{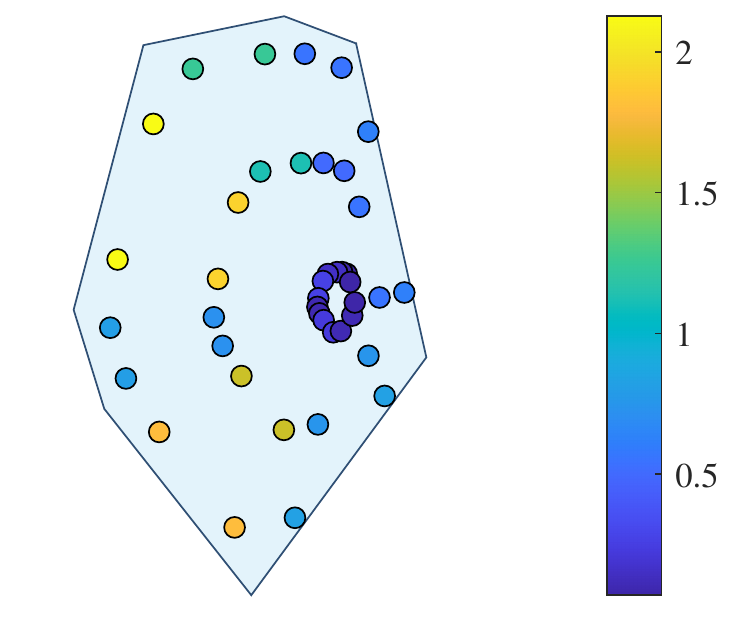}
    \caption{}\label{fig:ex-polygon-convex1-pts}
  \end{subfigure}
  \begin{subfigure}{3in}
    \centering
    \includegraphics[scale=1]{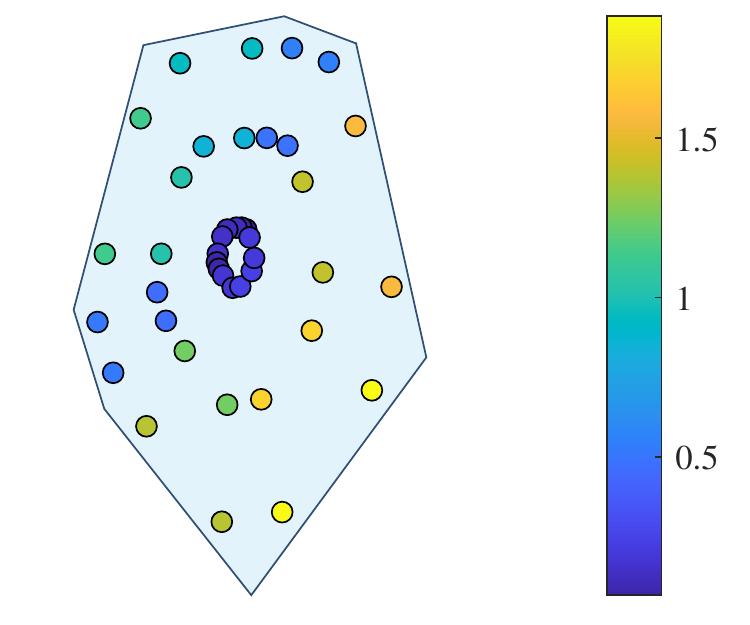}
    \caption{}\label{fig:ex-polygon-convex1-pts-shift}
  \end{subfigure}
  \begin{subfigure}{3in}
    \centering
    \includegraphics[scale=1]{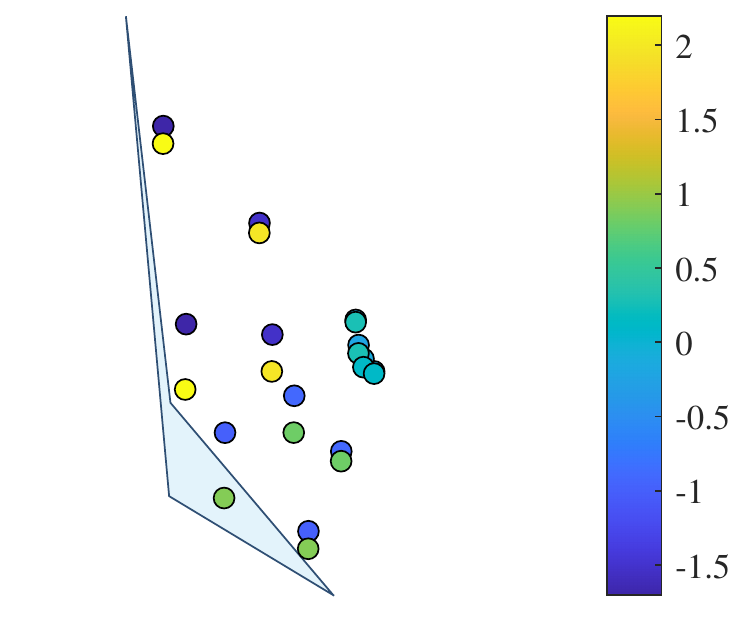}
    \caption{}\label{fig:ex-polygon-nonconvex1-pts}
  \end{subfigure}
  \begin{subfigure}{3in}
    \centering
    \includegraphics[scale=1]{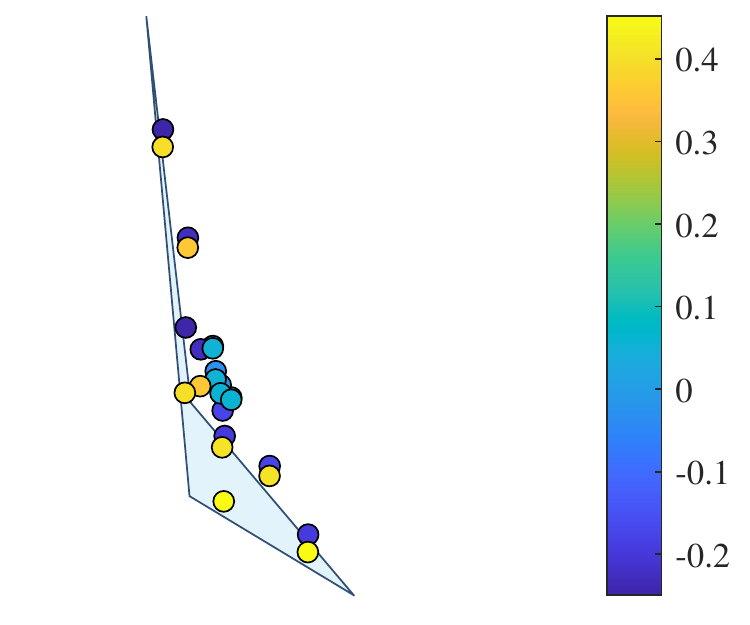}
    \caption{}\label{fig:ex-polygon-nonconvex1-pts-shift}
  \end{subfigure}
  \caption{Distribution of cubature points (closed circles) over polygons
  in~\fref{fig:polygon-convex-1} in (a) and (b) and
  in~\fref{fig:polygon-nonconvex-1} in (c) and (d). The point $\vx_0$ is the
  origin in (a) and (c) and $\vx_0$ is the average of the vertices in (b) and
  (d).  The shading of each point (and the scale on the right of each subfigure)
  corresponds to the cubature weight.}
  \label{fig:ex-polygon-pts}
\end{figure}

\smallskip
This study is performed for two cases: first with $\vx_0$ located at the vertex
average of the polygon and secondly with $\vx_0 = \vm{0}$.  While selecting
different locations for $\vx_0$ can affect the numerical accuracy of the
integral due to cancellation errors~\cite{Chin:2020:AEM}, integration with the
SB parametrization is capable of computing integrals with an arbitrary location
for $\vx_0$.  The distribution of cubature points and weights using SBC on the
polygons in~\fref{fig:polygon-convex-1} and~\fref{fig:polygon-nonconvex-1} for
$\vx_0$ at the origin and at the average of the polygonal vertices is depicted
in~\fref{fig:ex-polygon-pts}.  The cubature points on the polygons correspond to
three cubature points in the $\xi$ direction and two cubature points in the $t$ direction
per triangle.  When $\vx_0$ is chosen as the mean of the vertex
coordinates, cubature points are distributed more evenly across the convex
polygon in~\fref{fig:polygon-convex-1}.  For the nonconvex polygon
in~\fref{fig:polygon-nonconvex-1}, shifting $\vx_0$ to the average of the vertex
coordinates results in cubature points that are closer to the domain of
integration and the range of cubature weights is reduced.  Integration error
with the minimal number of Gauss\hyp{}Legendre points required for exact
integration is presented in~\tref{tab:ex-polynomials-error}.  The column
``Relative integration error (shifted $\vx_0$)'' refers to the case where
$\vx_0$ is the average of the polygonal vertex coordinates while the column
``Relative integration error ($\vx_0 = \vm{0}$)'' is for the case when $\vx_0$
is the origin. As \tref{tab:ex-polynomials-error} reveals, the majority of
integrals are computed with error that is close to machine precision
$\bigl($about $\mathcal{O}(10^{-15}) \bigr)$, regardless of the location of
$\vx_0$. For some integrals, error on the order of $\mathcal{O}(10^{-14})$ is
observed. Increased error in these integrals is due to cancellations caused by
subtraction of numbers with similar magnitude.  These results confirm the
predicted accuracy of integrating polynomials over polygons discussed
in~\sref{sec:polygons}.

{\renewcommand{\arraystretch}{1.4}
\begin{table}[t!]
  \centering
  \caption{Relative integration error for integrating polynomials over convex
  and nonconvex polygons using SBC.}
  \label{tab:ex-polynomials-error}
  \begin{tabular}{|c|c|c|c|p{1.38in}|p{1.38in}|}
    \hline
    $p$ & $\xi$\hyp{}points & $t$\hyp{}points 
      & Polygon in~\fref{fig:ex-polynomials-polygons} 
      & Relative integration error (shifted $\vx_0$) 
      & Relative integration error ($\vx_0 = \vm{0}$) \\
    \hline
    \multirow{4}{*}{$0$} & \multirow{4}{*}{$1$} & \multirow{4}{*}{$1$}
          & (a) & $2.2 \times 10^{-16}$ & $0$ \\
      \cline{4-6}
      & & & (b) & $2.2 \times 10^{-16}$ & $0$ \\
      \cline{4-6}
      & & & (c) & $1.2 \times 10^{-15}$ & $2.4 \times 10^{-15}$ \\
      \cline{4-6}
      & & & (d) & $0$ & $1.3 \times 10^{-16}$ \\
    \hline
    \multirow{4}{*}{$1$} & \multirow{4}{*}{$2$} & \multirow{4}{*}{$1$}
          & (a) & $0$ & $1.7 \times 10^{-16}$ \\
      \cline{4-6}
      & & & (b) & $6.4 \times 10^{-14}$ & $4.1 \times 10^{-15}$ \\
      \cline{4-6}
      & & & (c) & $1.5 \times 10^{-15}$ & $8.5 \times 10^{-15}$ \\
      \cline{4-6}
      & & & (d) & $6.4 \times 10^{-16}$ & $1.3 \times 10^{-16}$ \\
    \hline
    \multirow{4}{*}{$2$} & \multirow{4}{*}{$2$} & \multirow{4}{*}{$2$}
          & (a) & $2.8 \times 10^{-14}$ & $1.8 \times 10^{-14}$ \\
      \cline{4-6}
      & & & (b) & $0$ & $0$ \\
      \cline{4-6}
      & & & (c) & $7.1 \times 10^{-16}$ & $1.4 \times 10^{-16}$ \\
      \cline{4-6}
      & & & (d) & $5.2 \times 10^{-16}$ & $1.0 \times 10^{-15}$ \\
    \hline
    \multirow{4}{*}{$3$} & \multirow{4}{*}{$3$} & \multirow{4}{*}{$2$}
          & (a) & $2.8 \times 10^{-16}$ & $1.4 \times 10^{-16}$ \\
      \cline{4-6}
      & & & (b) & $3.5 \times 10^{-16}$ & $8.8 \times 10^{-16}$ \\
      \cline{4-6}
      & & & (c) & $1.1 \times 10^{-15}$ & $2.8 \times 10^{-15}$ \\
      \cline{4-6}
      & & & (d) & $2.2 \times 10^{-16}$ & $4.4 \times 10^{-16}$ \\
    \hline
    \multirow{4}{*}{$4$} & \multirow{4}{*}{$3$} & \multirow{4}{*}{$3$}
          & (a) & $4.7 \times 10^{-15}$ & $1.6 \times 10^{-15}$ \\
      \cline{4-6}
      & & & (b) & $0$ & $1.4 \times 10^{-16}$ \\
      \cline{4-6}
      & & & (c) & $2.8 \times 10^{-16}$ & $5.8 \times 10^{-15}$ \\
      \cline{4-6}
      & & & (d) & $4.0 \times 10^{-16}$ & $0$ \\
    \hline
    \multirow{4}{*}{$5$} & \multirow{4}{*}{$4$} & \multirow{4}{*}{$3$}
          & (a) & $0$ & $2.7 \times 10^{-16}$ \\
      \cline{4-6}
      & & & (b) & $2.1 \times 10^{-16}$ & $8.3 \times 10^{-16}$ \\
      \cline{4-6}
      & & & (c) & $1.6 \times 10^{-15}$ & $1.4 \times 10^{-15}$ \\
      \cline{4-6}
      & & & (d) & $6.3 \times 10^{-16}$ & $1.3 \times 10^{-15}$ \\
    \hline
  \end{tabular}
\end{table}
}

\subsubsection{Non\hyp{}polynomial integration}
\begin{figure}[t]
  \centering
  \begin{subfigure}{2in}
    \centering
    \includegraphics[scale=1]{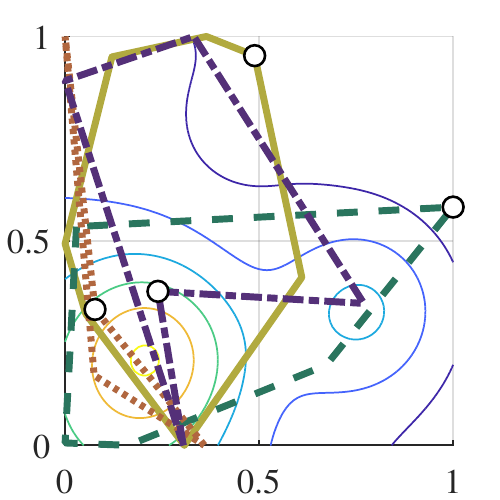}
    \caption{}\label{fig:ex-franke-function1}
  \end{subfigure}
  \begin{subfigure}{2in}
    \centering
    \includegraphics[scale=1]{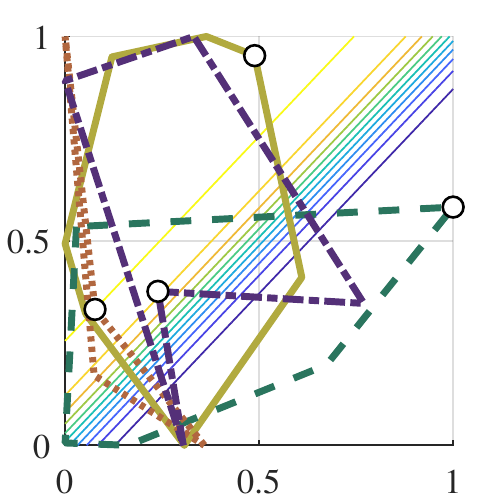}
    \caption{}\label{fig:ex-franke-function2}
  \end{subfigure}
  \begin{subfigure}{2in}
    \centering
    \includegraphics[scale=1]{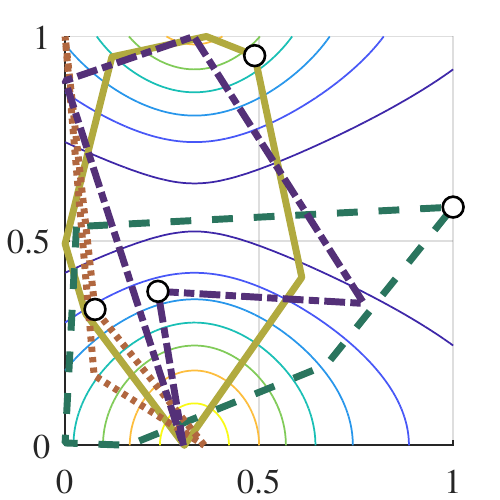}
    \caption{}\label{fig:ex-franke-function3}
  \end{subfigure}
  \caption{Contour maps of the three non\hyp{}polynomial functions
  in~\eqref{eq:franke-functions} that are integrated over the overlaid
  polygonal domains.  (a) $f_{F1} (\vx)$, (b) $f_{F2} (\vx)$, and (c) $f_{F3}
  (\vx)$. }
  \label{fig:ex-franke-functions}
\end{figure}

\begin{figure}[t!]
  \centering
  \begin{subfigure}{3in}
    \centering
    \includegraphics[scale=1]{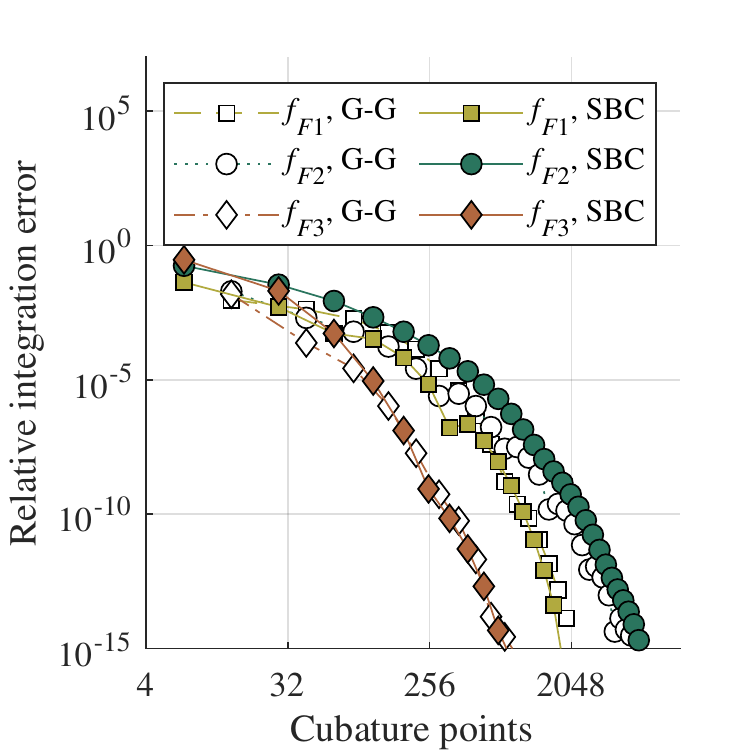}
    \caption{}\label{fig:ex-polygon-franke-conv-1}
  \end{subfigure}
  \begin{subfigure}{3in}
    \centering
    \includegraphics[scale=1]{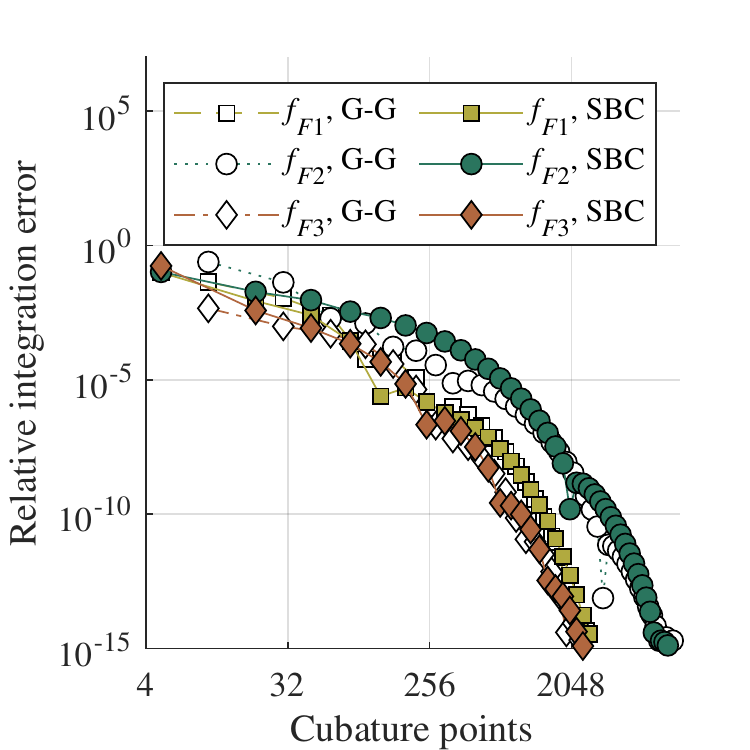}
    \caption{}\label{fig:ex-polygon-franke-conv-2}
  \end{subfigure}
  \begin{subfigure}{3in}
    \centering
    \includegraphics[scale=1]{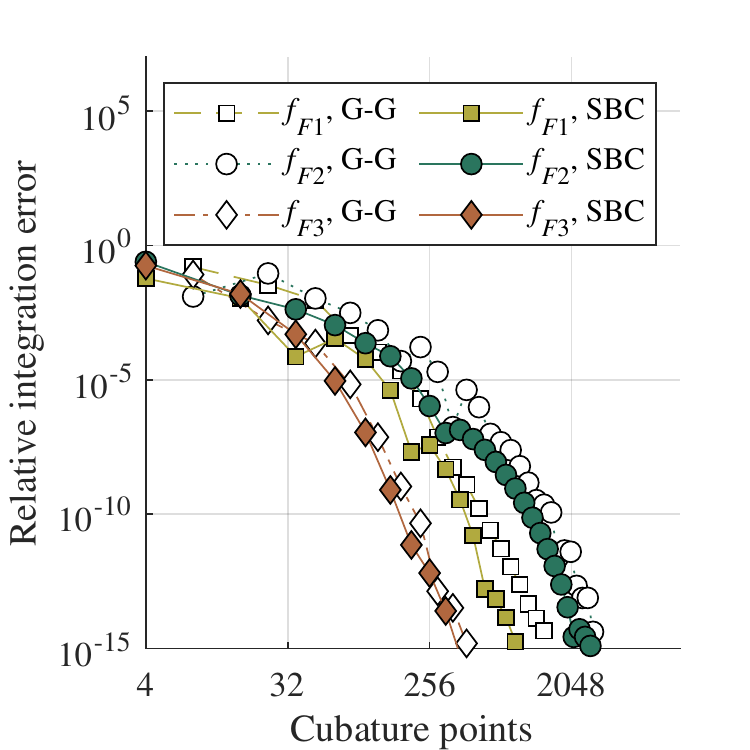}
    \caption{}\label{fig:ex-polygon-franke-conv-3}
  \end{subfigure}
  \begin{subfigure}{3in}
    \centering
    \includegraphics[scale=1]{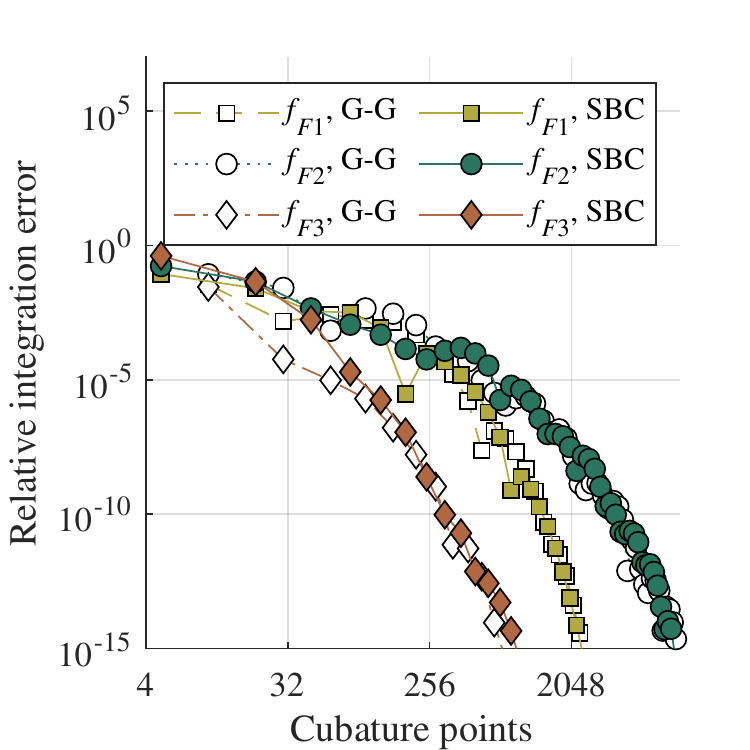}
    \caption{}\label{fig:ex-polygon-franke-conv-4}
  \end{subfigure}
  \caption{Relative integration error versus number of cubature points for the
  polygons in~\fref{fig:ex-polynomials-polygons} mapped to the unit square.
  Integration is performed using the SBC scheme with $\vx_0$ at the vertex
  average (filled markers) and with Gauss\hyp{}Green (G-G) cubature (open
  markers).  Polygon shape in (a) \fref{fig:polygon-convex-1}, (b)
  \fref{fig:polygon-convex-2}, (c) \fref{fig:polygon-nonconvex-1}, and (d)
  \fref{fig:polygon-nonconvex-2}. }
  \label{fig:ex-franke-conv}
\end{figure}

Franke~\cite{Franke:1979:ACC} used several non\hyp{}polynomial functions to
compare methods of interpolation for scattered data.  In this section, we
utilize these functions to validate the SBC method over polygons with
non\hyp{}polynomial integrands.  The functions we utilize are:
\begin{subequations}\label{eq:franke-functions}
\begin{align}
  f_{F1} (\vx) & = \frac{3}{4} \exp \left( 
      - \frac{(9 x - 2)^2 + (9 y - 2)^2}{4}
    \right) + \frac{3}{4} \exp \left(
      - \frac{(9 x + 1)^2}{49} - \frac{9 y+1}{10}
    \right) \nonumber\\
    & + \frac{1}{2} \exp \left(
      - \frac{(9 x - 7)^2 + (9 y - 3)^2}{4}
    \right) + \frac{1}{5} \exp \left(
      -(9 x - 4)^2 - (9 y - 7)^2
    \right) , \label{eq:franke-function}\\
  f_{F2} (\vx) & = \frac{1}{9} \left[ 
      \tanh \left( 9 y - 9 x \right) + 1
    \right] , \\
  f_{F3} (\vx) & = \frac{\frac{5}{4} 
      + \cos \left( \frac{27}{5} y \right)}
    {6 \left[ 1 + (3 x - 1)^2 \right]} .
\end{align}
\end{subequations}
We integrate these functions over the polygons
in~\fref{fig:ex-polynomials-polygons}.  In Franke~\cite{Franke:1979:ACC}, the
functions are designed to exhibit interesting features in the unit square,
$[0,1]^2$. However, the polygons in~\fref{fig:ex-polynomials-polygons} are not
contained in this domain.  Accordingly, we apply an affine translation and
scaling of each polygon such that the domain is the largest possible subset of
the unit square.  Contour maps of the three integrands are illustrated
in~\fref{fig:ex-franke-functions} with the scaled polygons overlaid. Integration
error versus the number of cubature points is revealed
in~\fref{fig:ex-franke-conv}.  For all functions and polygons, integration error
less than $\mathcal{O}(10^{-15})$ is obtained with fewer than $10,000$ cubature
points.  \revtwo{For the results in~\fref{fig:ex-franke-conv}, $\vx_0$ is the
average of the vertex coordinates of the polygons.}

\begin{figure}[t!]
  \centering
  \begin{subfigure}{3in}
    \centering
    \includegraphics[scale=1]{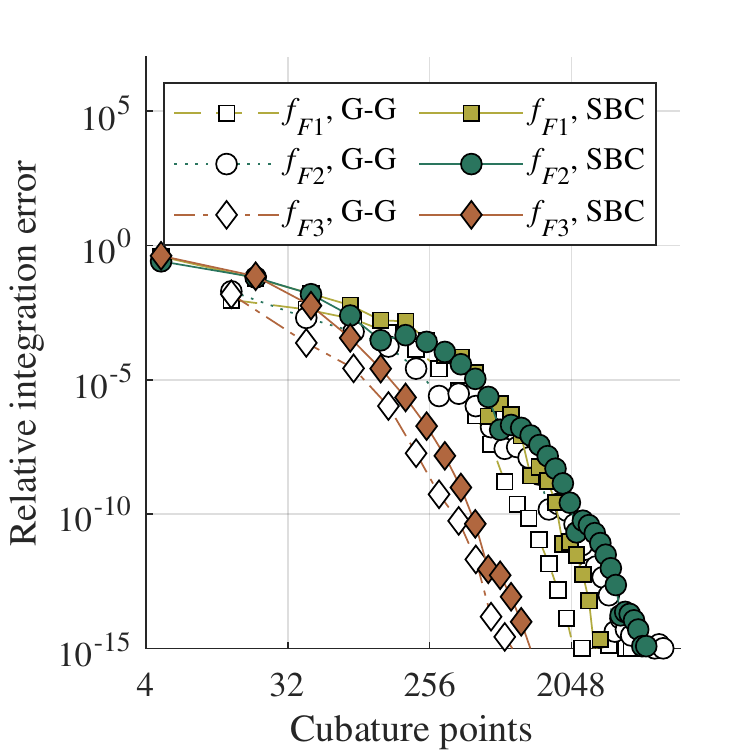}
    \caption{}\label{fig:ex-polygon-franke-conv-opt-1}
  \end{subfigure}
  \begin{subfigure}{3in}
    \centering
    \includegraphics[scale=1]{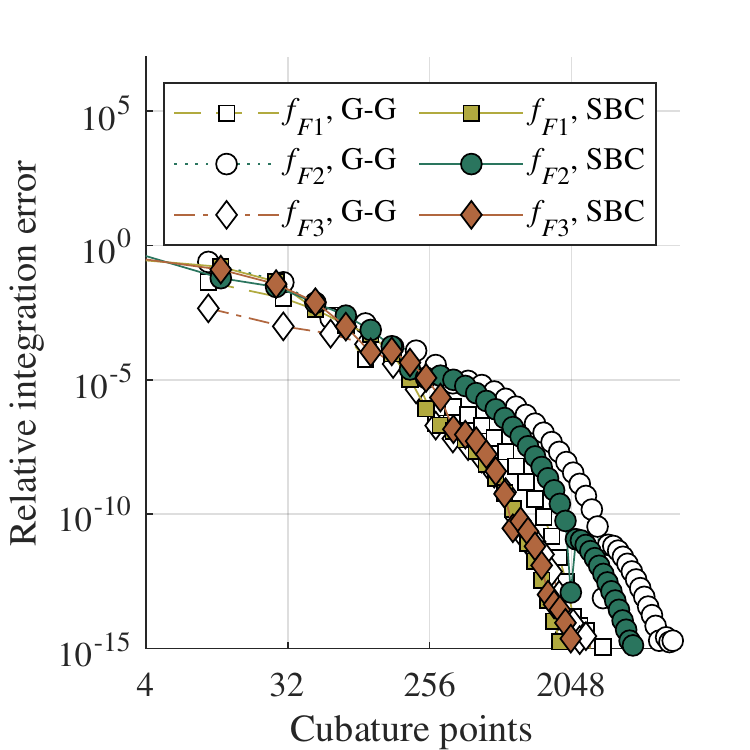}
    \caption{}\label{fig:ex-polygon-franke-conv-opt-2}
  \end{subfigure}
  \begin{subfigure}{3in}
    \centering
    \includegraphics[scale=1]{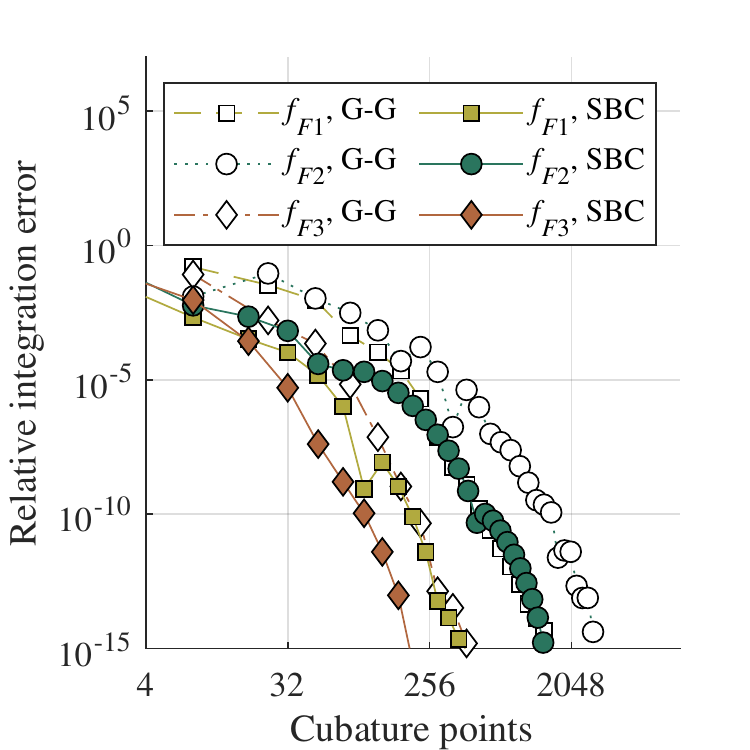}
    \caption{}\label{fig:ex-polygon-franke-conv-opt-3}
  \end{subfigure}
  \begin{subfigure}{3in}
    \centering
    \includegraphics[scale=1]{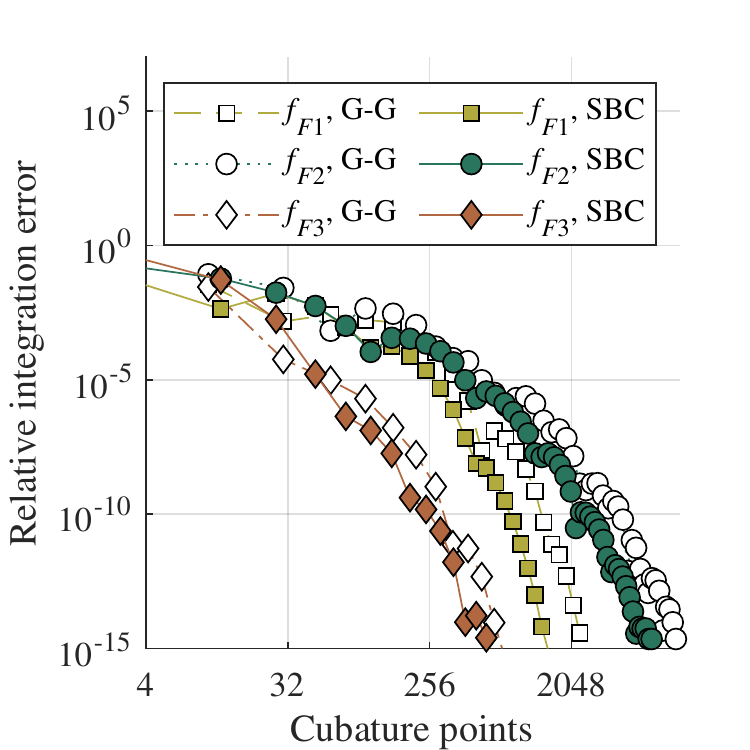}
    \caption{}\label{fig:ex-polygon-franke-conv-opt-4}
  \end{subfigure}
  \caption{Relative integration error versus number of cubature points for the
  polygons in~\fref{fig:ex-polynomials-polygons} mapped to the unit square.
  Integration is performed using SBC with $\vx_0$ at a vertex (filled markers)
  and with Gauss\hyp{}Green (G-G) cubature (open markers).  Polygon shape in (a)
  \fref{fig:polygon-convex-1}, (b) \fref{fig:polygon-convex-2}, (c)
  \fref{fig:polygon-nonconvex-1}, and (d) \fref{fig:polygon-nonconvex-2}. }
  \label{fig:ex-franke-conv-opt}
\end{figure}

\smallskip
For the purpose of comparison, we repeat the above tests using the G-G method of
Sommariva and Vianello~\cite{Sommariva:2007:PGC}.  These results are also
presented in~\fref{fig:ex-franke-conv}.  Integration accuracy with the G-G
scheme is on par with the SBC method.  However, the SBC scheme enables cubature
points to be placed entirely inside the domain of integration on more polygons
as compared to G-G integration, and furthermore, by placing $\vx_0$ at a vertex,
the number of SBC integration points is reduced without compromising the
polynomial precision of the cubature rule.  To explore these effects, the same
study in~\fref{fig:ex-franke-conv} is repeated with $\vx_0$ selected with
respect to star\hyp{}convexity and at a vertex of the polygon.  \revtwo{The
vertex selected as $\vx_0$ for each polygon is marked with an open circle in
\fref{fig:ex-franke-functions}.}  The results of this exercise are presented
in~\fref{fig:ex-franke-conv-opt}.  Over the two convex polygons tested, the
results are mixed.  Since the integrands are non\hyp{}polynomial, the accuracy
of the rule is strongly dependent on the location of the cubature points.  For
the polygon in~\fref{fig:polygon-convex-1}, placing $\vx_0$ at a vertex degrades
the performance of the SBC rule, whereas for the polygon
in~\fref{fig:polygon-convex-2}, integration error is reduced.  We note the
accuracy of the SBC rule is dependent on the vertex selected as $\vx_0$.  For
both nonconvex polygons, placing $\vx_0$ at a vertex selected with respect to
the star\hyp{}convexity of the domain results in improved cubature convergence.
For the quadrilateral domain in~\fref{fig:polygon-nonconvex-1}, the improvement
is most pronounced. For this domain, choosing $\vx_0$ to be located at a vertex
reduces the number of cubature points by a factor of one\hyp{}half without
reducing the polynomial precision of the rule.  We point out that for both
nonconvex polygons used, it is not possible to recover integration points that
are entirely within the domain using the G-G scheme.

\subsubsection{Timing to generate a cubature rule}
The time required to integrate a function to a desired accuracy given a cubature
rule is mainly dependent on two factors.
\begin{enumerate}
  \item The overhead time required to generate the integration points and
  weights.  This includes computing domain\hyp{}specific quantities, such as
  triangulations, mappings, and distances, and mapping a base integration rule
  (over an interval or over a triangle) of a given order to the domain of
  integration.
  \item The accuracy of the cubature method in integrating the given function.
  Less accurate methods require more cubature points to compensate, which can
  increase the time to generate the rule and require more function evaluations.
\end{enumerate}
In this example, we explore the first factor (overhead time) with three
benchmark cubature methods for polygons:  SBC, G-G integration, and a symmetric
Dunavant~\cite{Dunavant:1985:HDE} rule applied to a constrained Delaunay
triangulation (CDT) of the integration domain.  For a given cubature rule, the
computational costs are proportional to the number of cubature points.  More
cubature points require more function evaluations.  Timing to perform a function
evaluation is strongly dependent on the function, so we exclude this aspect from
our study in this section.  The other examples in this paper explore the economy
of the SBC method in terms of accuracy per cubature point for various types of
functions.

\smallskip
To ensure consistent timing computations, all examples are repeated 100 times,
and the times reported here are the average of the 100 executions.  All code is
run in MATLAB 9.8.0 on a Linux computer with an Intel Xeon E5-2695 v4 CPU and
128 GB RAM.  A MATLAB implementation of Gauss\hyp{}Green cubature is provided by
the \texttt{Polygauss} code from Sommariva.\footnote{Accessed from
\texttt{https://www.math.unipd.it/\~{}alvise/SOFTWARE/GAUSSCUB\_2013.zip}.} The
CDT is generated from the \texttt{delaunayTriangulation} class in MATLAB.
Implementations of the SBC method and the Dunavant rule are by the authors.  The
code for all methods (including \texttt{Polygauss}) is optimized through array
preallocation and by eliminating hot spots using the MATLAB Profiler, where
possible.

\begin{figure}
  \centering
  \includegraphics[scale=1]{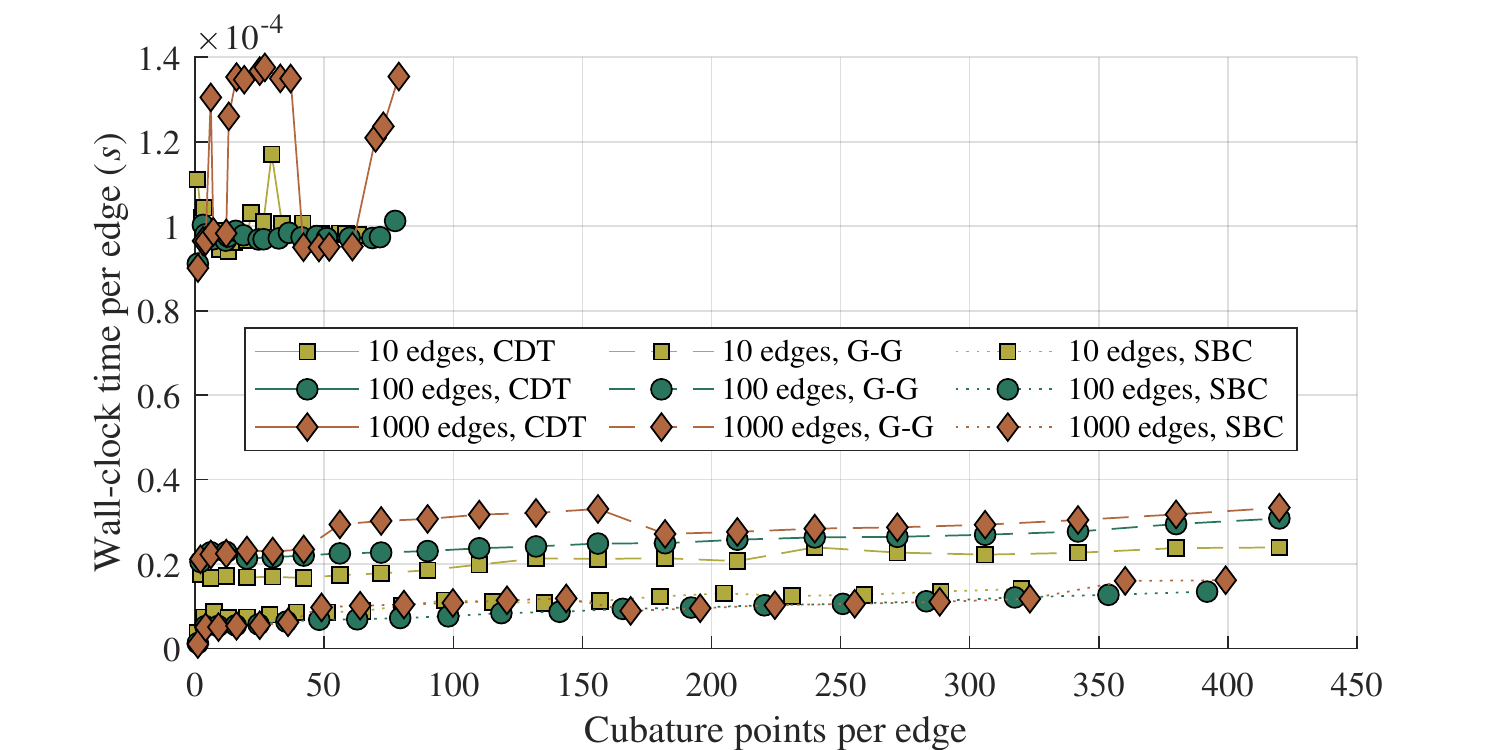}
  \caption{Time per edge to generate a cubature rule over polygons inscribed on
    a circle for the three benchmark methods. }
  \label{fig:ex-timing-overhead}
\end{figure}

\smallskip
To examine the effects of overhead time on polygons of varying complexity, we
compute the time to generate a rule versus the number of cubature points on
polygons with 10 sides, 100 sides, and 1000 sides using the three benchmark
methods.  The polygons are cyclic polygons generated by joining
evenly\hyp{}spaced points on a circle.  For the G-G and SBC methods, cubature
rules with different numbers of points per edge are generated by increasing the
number of Gauss points in each direction of the tensor\hyp{}product rule.  For
the CDT, the number of cubature points per edge is increased using a higher
order Dunavant rule.  \fref{fig:ex-timing-overhead} presents the numerical
results of this study.  The SBC rule is the fastest.  Rules of order $1 \times
1$ to $20 \times 20$ per edge are generated in less than $1.7 \times 10^{-5}$
seconds. This is approximately half the time required to generate the G-G rule
and about ten percent of the time required to generate the CDT/Dunavant rule.
For all three methods of integration, the time required to generate a rule per
edge is approximately constant and increasing the number of cubature points per
edge only slightly affects the time required to generate the integration rule.
These findings suggest the majority of the wall\hyp{}clock time is spent on
polygon\hyp{}specific tasks, such as generating a triangulation, mapping
domains, and computing distances.  Finally, we remark the CDT generation process
is eased by the simplicity of the polygonal domain in this example.  For
complicated polygonal domains (e.g.\ nonconvex or multiply\hyp{}connected
domains), the time to generate the CDT will increase.

\subsubsection{Wachspress coordinates}
\begin{figure}[t!]
  \centering
  \begin{subfigure}{3in}
    \centering
    \includegraphics[scale=1]{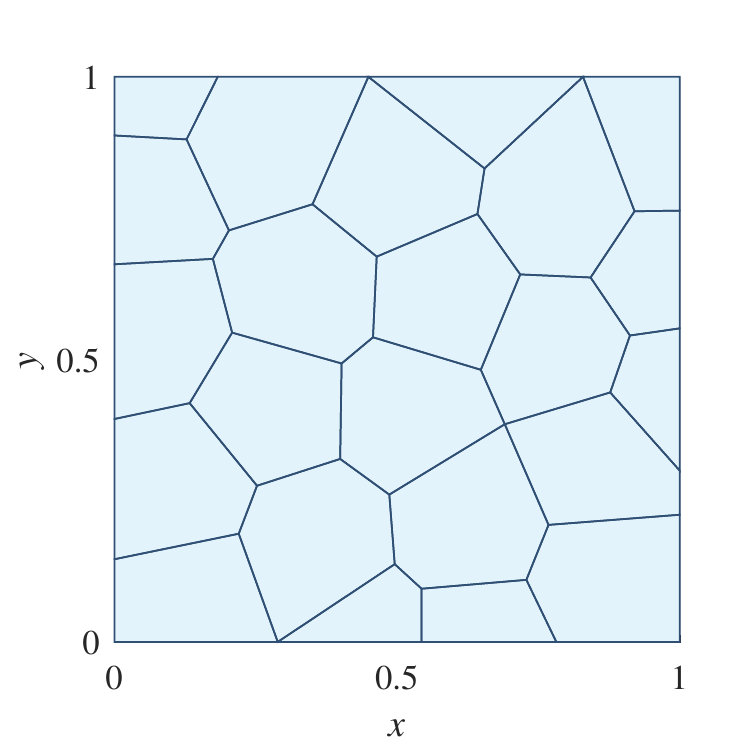}
    \caption{}\label{fig:ex-wachspress-mesh}
  \end{subfigure} \\
  \begin{subfigure}{2in}
    \centering
    \includegraphics[scale=1]{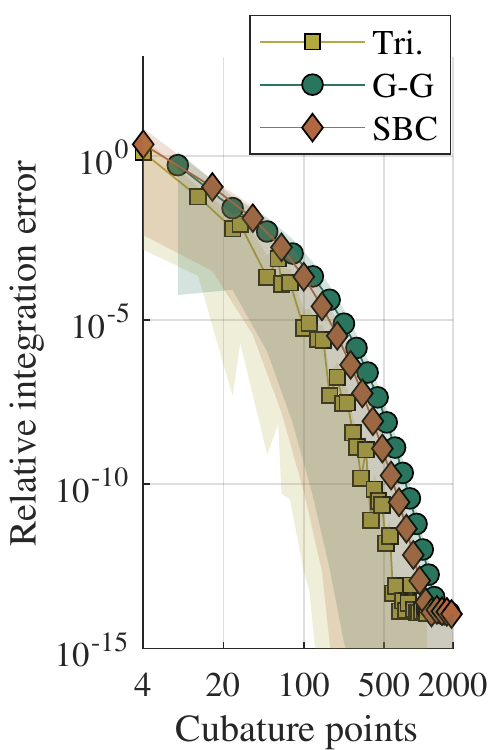}
    \caption{}\label{fig:ex-wachspress-4-sides}
  \end{subfigure}
  \begin{subfigure}{2in}
    \centering
    \includegraphics[scale=1]{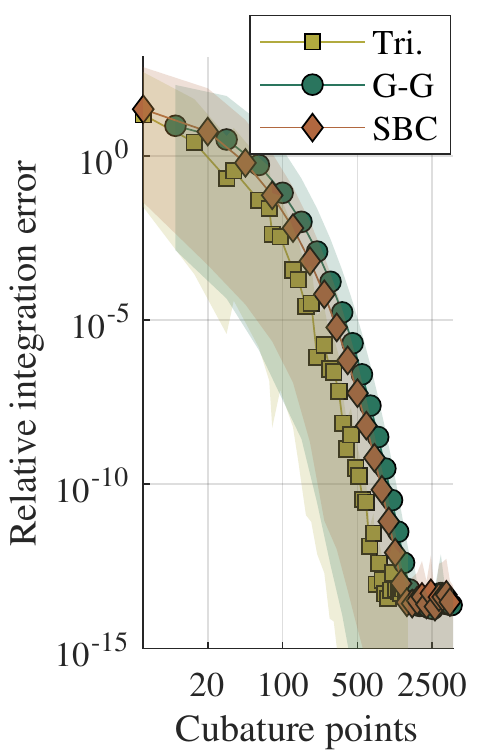}
    \caption{}\label{fig:ex-wachspress-5-sides}
  \end{subfigure}
  \begin{subfigure}{2in}
    \centering
    \includegraphics[scale=1]{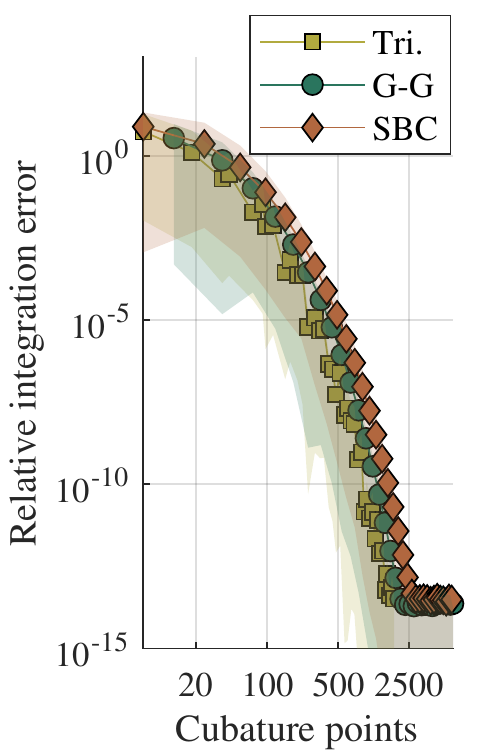}
    \caption{}\label{fig:ex-wachspress-6-sides}
  \end{subfigure}
  \caption{Integrating the products of derivatives of Wachspress basis functions
    using three cubature rules: triangulation (Tri.), G-G cubature, and SBC. (a)
    Polygonal mesh. Integration error versus cubature points over (b)
    quadrilateral elements, (c) pentagonal elements, and (d) hexagonal elements.
    The line denotes the average integration error while the shaded region
    bounds the minimum and maximum integration error.}
  \label{fig:ex-wachspress}
\end{figure}

One potential application of the SBC scheme is in generalized barycentric
coordinates (GBCs)~\cite{Hormann:2017:GBC}, which extends the application of
barycentric coordinates to polytopes.  GBCs have been successfully utilized to
solve boundary\hyp{}value problems under the umbrella of polytopal finite
element methods.  One well\hyp{}known GBC is the Wachspress coordinates (basis
functions)~\cite{Wachspress:1975:ARF}, which uses ideas from projective geometry
to form a rational finite element basis over convex polygonal domains.
Simplified equations to compute Wachspress basis functions are available in
Floater et al.~\cite{Floater:2014:GBW} and they are utilized here.  In this
example, a random polygonal mesh is formed over a square domain through the
Voronoi tesselation of a random point sampling generated using a technique
called maximal Poisson\hyp{}disk sampling.  The mesh used is illustrated
in~\fref{fig:ex-wachspress-mesh}.  Then, the product of derivatives of
Wachspress basis functions are integrated over the polygons using (i) SBC, (ii)
triangulating from the vertex average and using a triangle rule from Xiao and
Gimbutas~\cite{Xiao:2009:ANA}, and (iii) G-G cubature.  The purpose of this
example is to demonstrate the utility of SBC in a realistic scenario.  The
integrand is selected to mimic a typical finite element stiffness matrix
computation.  For SBC cubature, $\vx_0$ is selected as the average of the vertex
coordinates of the polygon.  Similarly, triangulation is performed by connecting
the edges of the polygon with the average of the vertex coordinates of each
polygon.

\smallskip
To investigate convergence properties of the three cubature methods, we
integrate over each polygon in the mesh and compute the relative error of
integration as compared to a reference solution.  Cubature rules of order 1 to
40 are used with each rule.  The reference solution is computed using G-G
cubature with 50 points in each direction per edge.  Average integration
accuracy versus number of cubature points is plotted
in~\fref{fig:ex-wachspress}. Separate plots are presented for quadrilateral
elements (\fref{fig:ex-wachspress-4-sides}), pentagonal elements
(\fref{fig:ex-wachspress-5-sides}), and hexagonal elements
(\fref{fig:ex-wachspress-6-sides}).  In these figures, the line denotes the
average integration error over all of the basis functions integrated.  The
shaded regions give the range of integration error.

\smallskip
The results of \fref{fig:ex-wachspress} indicate integration using a triangle
integration rule is, on average, slightly more accurate than using either SBC or
G-G when integrating with the same number of integration points.  This is
expected since the Xiao\hyp{}Gimbutas triangle rule is a symmetric rule devised
specifically for triangles while both SBC and G-G cubature are
tensor\hyp{}product rules.  When comparing G-G integration to SBC, both methods
exhibit similar error per cubature point.  While the triangle rule is more
accurate on average, the ranges of error overlap in all plots, indicating the
relative error between the three methods is small.  While monotonic convergence
does not generally occur with non\hyp{}polynomial integrands, with SBC and G-G
cubature, error reduces monotonically with additional cubature points.  This
does not occur with the triangle integration rule.  Also, G-G cubature and SBC
rules can be readily adapted to integrate polynomials of arbitrary order; on the
other hand, symmetric triangle rules from moment fitting are time consuming
to generate.

\subsection{\revone{Planar domains} bounded by parametric curves}\label{sec:ex-curved}
\begin{figure}
  \centering
  \includegraphics[scale=1]{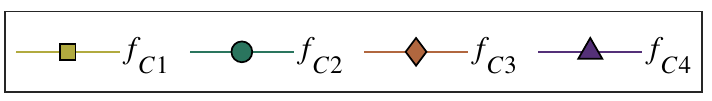}\\
  \begin{subfigure}{2in}
    \centering
    \includegraphics[scale=1]{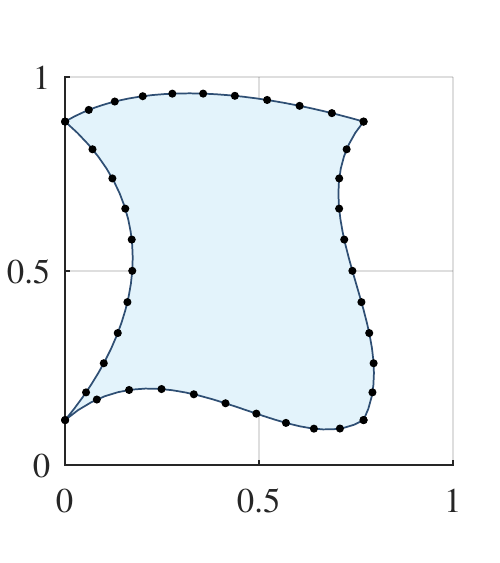}
    \caption{}\label{fig:ex-curve-shape1}
  \end{subfigure}
  \begin{subfigure}{2in}
    \centering
    \includegraphics[scale=1]{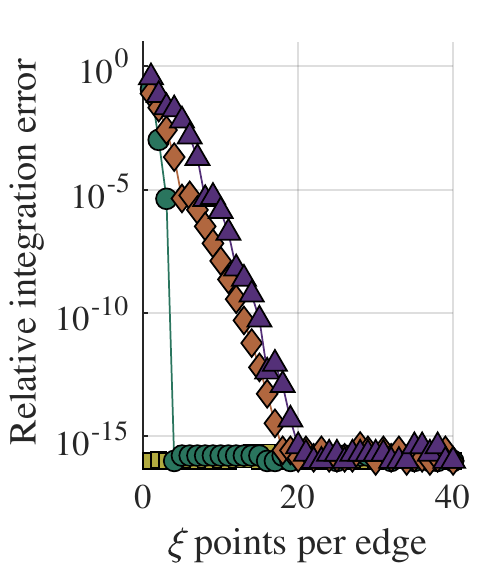}
    \caption{}\label{fig:ex-curve-conv-xi1}
  \end{subfigure}
  \begin{subfigure}{2in}
    \centering
    \includegraphics[scale=1]{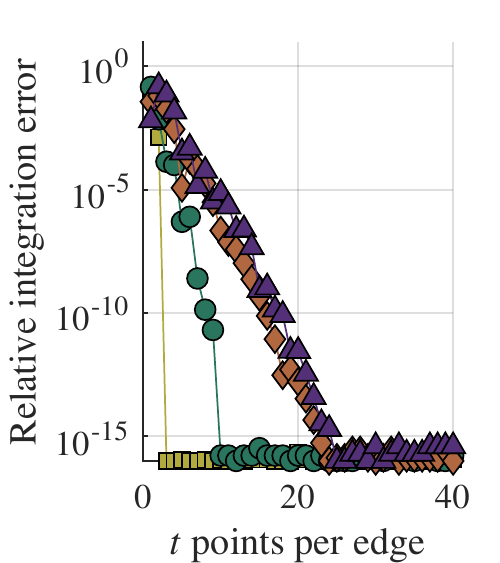}
    \caption{}\label{fig:ex-curve-conv-t1}
  \end{subfigure}
  \begin{subfigure}{2in}
    \centering
    \includegraphics[scale=1]{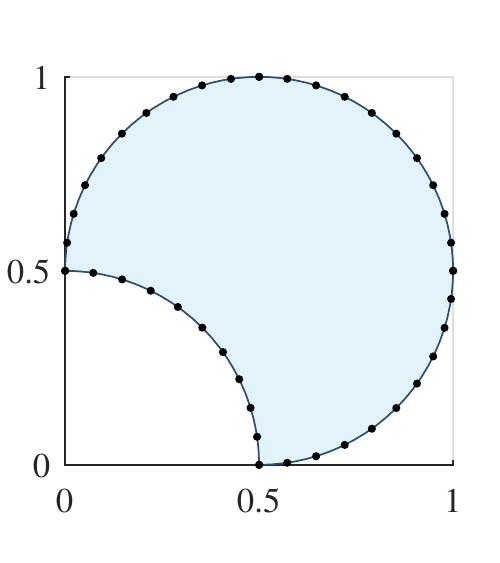}
    \caption{}\label{fig:ex-curve-shape2}
  \end{subfigure}
  \begin{subfigure}{2in}
    \centering
    \includegraphics[scale=1]{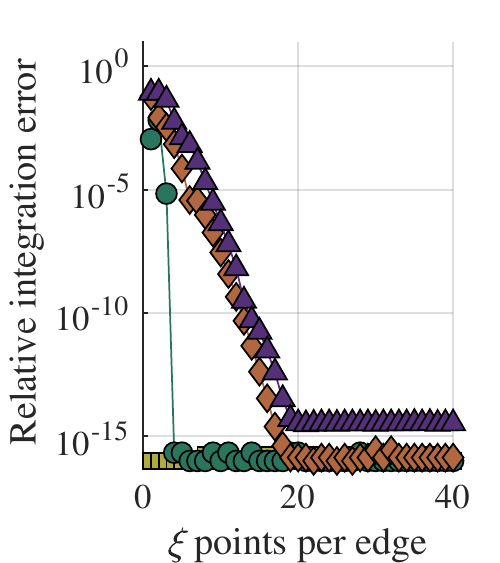}
    \caption{}\label{fig:ex-curve-conv-xi2}
  \end{subfigure}
  \begin{subfigure}{2in}
    \centering
    \includegraphics[scale=1]{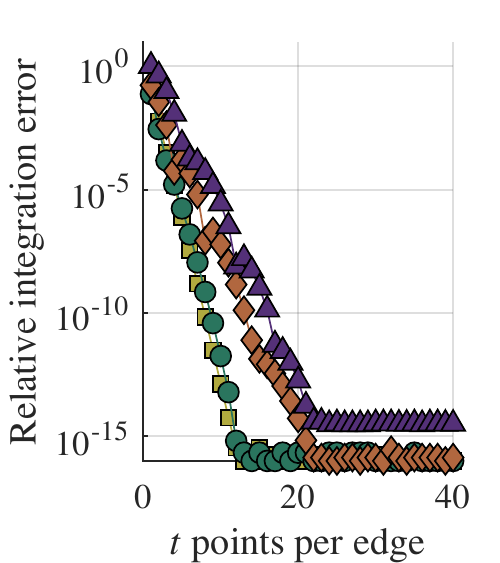}
    \caption{}\label{fig:ex-curve-conv-t2}
  \end{subfigure}
  \begin{subfigure}{2in}
    \centering
    \includegraphics[scale=1]{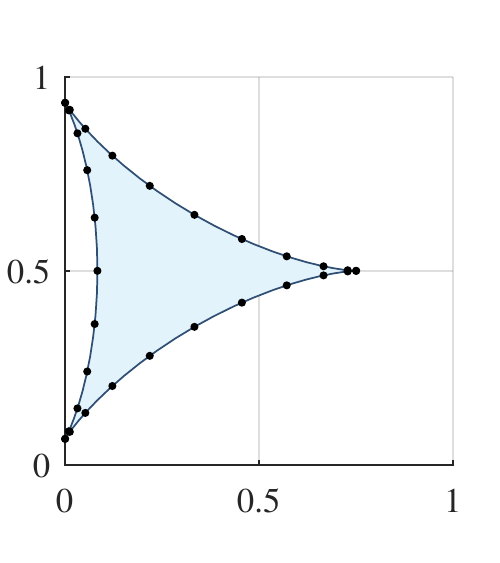}
    \caption{}\label{fig:ex-curve-shape3}
  \end{subfigure}
  \begin{subfigure}{2in}
    \centering
    \includegraphics[scale=1]{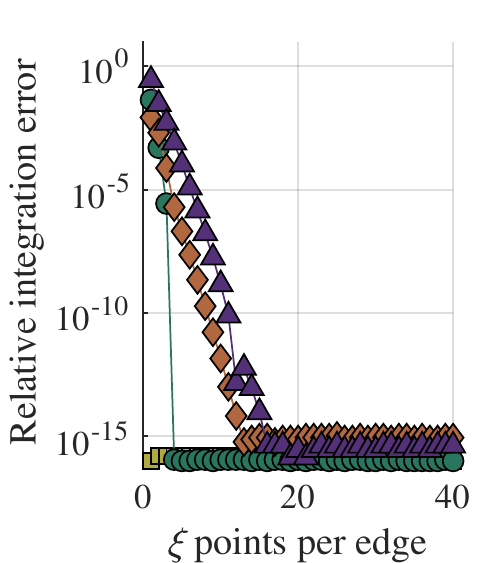}
    \caption{}\label{fig:ex-curve-conv-xi3}
  \end{subfigure}
  \begin{subfigure}{2in}
    \centering
    \includegraphics[scale=1]{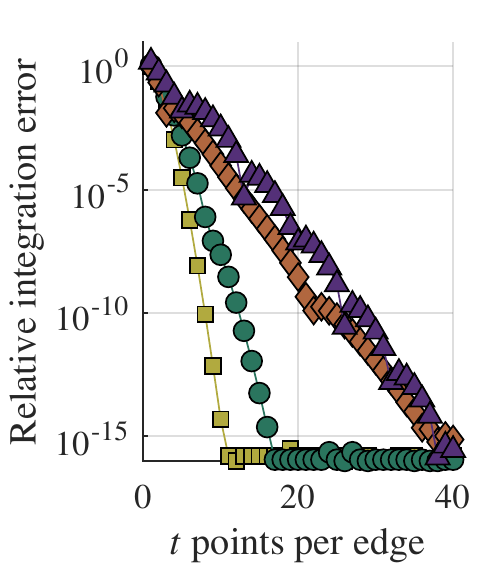}
    \caption{}\label{fig:ex-curve-conv-t3}
  \end{subfigure}
  \caption{Number of quadrature points versus integration accuracy to integrate
  four functions over regions bounded by curves.  (a) Domain bounded by cubic
  B\'{e}zier curves and integration accuracy versus number of (b) $\xi$
  quadrature points and (c) $t$ quadrature points.  (d) Lune (boundary
  represented by rational quadratic B\'{e}zier curves) and integration accuracy
  versus number of (e) $\xi$ quadrature points and (f) $t$ quadrature points.
  (g) Deltoid and integration accuracy versus number of (h) $\xi$ quadrature
  points and (i) $t$ quadrature points.}
  \label{fig:ex-curve}
\end{figure}

In this section, we showcase the capabilities of SBC in computing integrals of
polynomial and non\hyp{}polynomial functions over curved domains bounded by
polynomial and non\hyp{}polynomial parametric curves.  The functions integrated
are:
\begin{subequations}\label{eq:curved-fns}
\begin{align}
  f_{C1} (\vx) & = 1 , \\
  f_{C2} (\vx) & = 10x^5 - 5x^4y - 7x^3y^2 + 6x^2y^3 + 3xy^4 + y^5 - x^4 
      + 2x^3y + 11x^2y^2 - 8xy^3 - 2y^4 - 3x^3 + 9x^2y \nonumber\\
    & + 8xy^2 - 10y^3 - 9x^2 - 6xy + 7y^2 + 5x - 4y + 4 , \\
  f_{C3} (\vx) & = \frac{3}{4} \exp \left( 
    - \frac{(9 x - 2)^2 + (9 y - 2)^2}{4}
  \right) + \frac{3}{4} \exp \left(
    - \frac{(9 x + 1)^2}{49} - \frac{9 y+1}{10}
  \right) \nonumber\\
  & + \frac{1}{2} \exp \left(
    - \frac{(9 x - 7)^2 + (9 y - 3)^2}{4}
  \right) + \frac{1}{5} \exp \left(
    -(9 x - 4)^2 - (9 y - 7)^2
  \right) , \\
  f_{C4} (\vx) & = \exp \left( -\left[ 
      \left( \frac{x - 0.4}{0.3} \right)^2 + \left( \frac{y - 0.5}{0.4} \right)^2
    \right]^2 \right) \, \cos^2 \left( 3x \right) \, \cos^2 \left( 8y \right) .
\end{align}
\end{subequations}
In~\eqref{eq:curved-fns}, the first function is a constant (computes the area of
each curved domain), the second function is a 5th degree polynomial, the third
function is the Franke test function (introduced in~\eqref{eq:franke-function}),
and the fourth function is the product of exponential and trigonometric
functions.  Each of these functions is integrated over three curved \revone{domains} using
the SBC method.  The first domain (pictured in~\fref{fig:ex-curve-shape1}) is
bounded by four polynomial curves that are described using cubic B\'{e}zier
functions.  The control points of the B\'{e}zier curves are provided
in~\tref{tab:bezier1-points}.  The second domain, from Sommariva and
Vianello~\cite{Sommariva:2009:GGC}, is a lune that is illustrated
in~\fref{fig:ex-curve-shape2}.  The boundary of the lune is described using four
rational quadratic B\'{e}zier (non\hyp{}polynomial) curves, which are capable of
exactly representing conic sections.  The final domain is a deltoid that is
shown in~\fref{fig:ex-curve-shape3}.  While a single parametric equation can be
used to describe the entire deltoid boundary, the boundary is subdivided into
three curves to avoid cusps along each curve.  These three trigonometric
parametric equations are given in~\tref{tab:deltoid-eqns}.  For the three
domains, the point $\vx_0$ is selected as the average of the endpoints of each
curve.  With this selection, all cubature points are inside their respective
domains and all weights are positive, despite all three shapes being nonconvex.
Exact integrals are computed in Mathematica 12.0.0 using the functions
\texttt{Integrate[]} and \texttt{NIntegrate[]}.  Beyond applications in
numerical integration, we note the SB parametrization can be used to make
integration over complicated domains, such as the deltoid, tractable using
traditional methods of integration.

{\renewcommand{\arraystretch}{2.1}
\begin{table}[t]
  \centering
  \caption{Control points for the cubic B\'{e}zier curves that form the
  boundary of the domain pictured in~\fref{fig:ex-curve-shape1}.}
  \label{tab:bezier1-points}
  \begin{tabular}{|c|c|c|c|c|}
    \hline
    Curve & Point 1 & Point 2 & Point 3 & Point 4 \\[0.05in]
    \hline
    1 & $\left(0, \dfrac{3}{26}\right)$ & $\left(\dfrac{7}{26}, \dfrac{9}{26}\right)$ 
      & $\left(\dfrac{15}{26}, 0\right)$ & $\left(\dfrac{10}{13}, \dfrac{3}{26}\right)$ \\[0.05in]
    \hline
    2 & $\left(\dfrac{10}{13}, \dfrac{3}{26}\right)$ & $\left(\dfrac{23}{26}, \dfrac{9}{26}\right)$ 
      & $\left(\dfrac{15}{26}, \dfrac{17}{26}\right)$ & $\left(\dfrac{10}{13}, \dfrac{23}{26}\right)$ \\[0.05in]
    \hline
    3 & $\left(\dfrac{10}{13}, \dfrac{23}{26}\right)$ & $\left(\dfrac{1}{2}, \dfrac{25}{26}\right)$ 
      & $\left(\dfrac{5}{26}, 1\right)$ & $\left(0, \dfrac{23}{26}\right)$ \\[0.05in]
    \hline
    4 & $\left(0, \dfrac{23}{26}\right)$ & $\left(\dfrac{7}{26}, \dfrac{17}{26}\right)$ 
      & $\left(\dfrac{5}{26}, \dfrac{9}{26}\right)$ & $\left(0, \dfrac{3}{26}\right)$ \\[0.05in]
    \hline
  \end{tabular}
\end{table}
}

{\renewcommand{\arraystretch}{2.75}
\begin{table}[t]
  \centering
  \caption{Parametric equations for the boundary of the deltoid pictured
  in~\fref{fig:ex-curve-shape3}.}
  \label{tab:deltoid-eqns}
  \begin{tabular}{|c|c|}
    \hline
    $i$ & $\vm{c}_i$ \\
    \hline
    1 & $\begin{bmatrix}
      \dfrac{\left\{1 + 2 \cos \left( \frac{2 \pi}{3} t \right)\right\}^2}{12} &
      \dfrac{1 + 2 \sin \left( \frac{2 \pi}{3} t \right) - 
        4 \sin \left( \frac{4 \pi}{3} t \right)}{6}
    \end{bmatrix}^T$ \\[0.1in]
    \hline
    2 & $\begin{bmatrix}
      \dfrac{\left\{1 + 2 \cos \left( \frac{2 \pi}{3} \left( t + 1 \right) \right)\right\}^2}{12}
        & \dfrac{1 + 2 \sin \left( \frac{2 \pi}{3} \left( t + 1 \right) \right) - 
        4 \sin \left( \frac{4 \pi}{3} \left( t + 1 \right) \right)}{6}
    \end{bmatrix}^T$ \\[0.1in]
    \hline
    3 & $\begin{bmatrix}
      \dfrac{\left\{1 + 2 \cos \left( \frac{2 \pi}{3} \left( t + 2 \right) \right)\right\}^2}{12}
        & \dfrac{1 + 2 \sin \left( \frac{2 \pi}{3} \left( t + 2 \right) \right) - 
        4 \sin \left( \frac{4 \pi}{3} \left( t + 2 \right) \right)}{6}
    \end{bmatrix}^T$ \\[0.1in]
    \hline
  \end{tabular}
\end{table}
}

\smallskip
To investigate cubature convergence in each direction of the tensor\hyp{}product
domain after applying the SB mapping, integration error versus the number of
cubature points in each direction for the four test functions and three shapes
is plotted in~\fref{fig:ex-curve}.  Convergence properties in the
$\xi$ direction (resp.\ $t$ direction) are interrogated by supplying a
surplus of integration points in the $t$ direction (resp.\
$\xi$ direction), which ensures integration error is entirely in the
$\xi$ direction (resp.\ $t$ direction). For the two
non\hyp{}polynomial integrands, approximately $20 \times 20$ cubature points are
needed to realize error close to machine precision over the polynomial curve
domain and the lune. About 40 $t$\hyp{}points and 15 $\xi$\hyp{}points are
needed for the deltoid.  In~Figs.~\ref{fig:ex-curve-shape1},
\ref{fig:ex-curve-shape2}, and~\ref{fig:ex-curve-shape3}, dots are placed on the
boundary at intervals of $\Delta t = 0.1$.  Note the dots are not evenly spaced
on the boundary of the deltoid, since the curve speed $\bigl(\| \vm{c}^\prime
(t) \|\bigr)$ approaches zero around the cusps.  This results in extra weighting
of the parameter space $t$ toward the cusps, and ultimately, more cubature
points near the cusps.  The accuracy per cubature point for integrating
non\hyp{}polynomial functions over the deltoid would likely increase with a
parametrization of the boundary with less variation in curve speed.

\smallskip
While integrating using the SBC method results in a tensor\hyp{}product cubature
rule over each curved triangle, the burden of integration is not shared equally
in each direction.  Besides integrating the transformed integrand, the Jacobian
introduces a factor of $\xi$ and the product of the curve $\vm{c} (t)$ with its
derivative (hodograph).  For polynomial integrands and polynomial boundary
curves, the polynomial order of the integrand in the 
$\xi$- and $t$-directions
can be exactly computed (see~\sref{sec:sbc-curved-solids}). For a degree 0
polynomial bounded by polynomial curves of degree 3, the Gauss rule must be
capable of integrating a polynomial of degree 1 in the 
$\xi$ direction and
degree 5 in the $t$ direction.  This corresponds to 
1 Gauss point in the $\xi$ direction and 3 Gauss points in the $t$ direction for exact integration.  This
matches the results presented for $f_{C1}$ in~\fref{fig:ex-curve-conv-xi1}
and~\fref{fig:ex-curve-conv-t1}. When the polynomial integrand increases to
degree 5 and the bounding curves remain degree 3, the Gauss rule must be capable
of integrating a 6th degree polynomial in $\xi$ and a polynomial of degree 20 in
$t$.  Using a Gauss Legendre rule, 4 points are needed in the 
$\xi$ direction
and 11 points are needed in the $t$ direction. \fref{fig:ex-curve-conv-xi1} and
\fref{fig:ex-curve-conv-t1} reflect this for $f_{C2}$.  Note the number of
quadrature points in the $\xi$ direction is not dependent on the curve
$\vm{c}(t)$.  Therefore, the number of $\xi$\hyp{}points needed in the Gauss
rule is constant for polynomial integrands.  This holds for the lune and deltoid
domains, which are bounded by non\hyp{}polynomial curves
(see~\fref{fig:ex-curve-conv-xi2} and~\fref{fig:ex-curve-conv-xi3}).

\subsection{Integration of weakly singular functions}\label{sec:ex-singular}

\subsubsection{Comparison of $\xi$\hyp{}transformation methods}\label{sec:ex-singular-xi}
\sref{sec:singular-xi} introduced the generalized SB transformation for
cancelling radial singularities in the integrand.  This section compares the
generalized SB transformation and Gauss\hyp{}Jacobi integration to the SB
transformation for integrating different weakly singular functions.
Gauss\hyp{}Jacobi integration rules are constructed by solving for the roots of
the Jacobi polynomial in Mathematica 12.0.0.  Four different functions are
utilized for this comparison:
\begin{subequations}\label{eq:singular-fns-xi}
\begin{align}
  f_{S1} (\vx) & = \frac{4 - 2x + y - x^2 + 2xy - 3y^2 + 3x^3 - 5x^2y + 5xy^2 - 4y^3}
    {(x^2 + y^2)^{1/4}} , \\
  \label{eq:f-s2}
  f_{S2} (\vx) & = \frac{\exp \left( -\left[ \left(\frac{x - 0.25}{0.4} \right)^2
    + \left( \frac{y - 0.2}{0.7} \right)^2 \right]^2 \right)
    \, \cos^2 (5x) \, \cos^2 (5y)}
    {(x^2 + y^2)^{1/4}} , \\
  f_{S3} (\vx) & = \frac{4 - 2x + y - x^2 + 2xy - 3y^2 + 3x^3 - 5x^2y + 5xy^2 - 4y^3}
    {(x^2 + y^2)^{9/10}} , \\
  \label{eq:f-s4}
  f_{S4} (\vx) & = \frac{\exp \left( -\left[ \left(\frac{x - 0.25}{0.4} \right)^2
    + \left( \frac{y - 0.2}{0.7} \right)^2 \right]^2 \right)
    \, \cos^2 (5x) \, \cos^2 (5y)}
    {(x^2 + y^2)^{9/10}} .
\end{align}
\end{subequations}
The first two functions ($f_{S1} (\vx)$ and $f_{S2} (\vx)$) contain $r^{-1/2}$
singularities and the last two functions ($f_{S3} (\vx)$ and $f_{S4} (\vx)$)
contain $r^{9/5}$ in the denominator.  Two of the functions ($f_{S1} (\vx)$ and
$f_{S3} (\vx)$) contain the same cubic polynomial numerator and the remaining
two functions ($f_{S2} (\vx)$ and $f_{S4} (\vx)$) contain the same exponential
numerator.  We note all of these integrands are weakly singular, and contain a
fractional power for the singularity.  For $r^{-1}$ singularities, the SB
transformation is capable of accurate integration in the 
$\xi$ direction without further manipulation.

\begin{figure}
  \centering
  \includegraphics[scale=1]{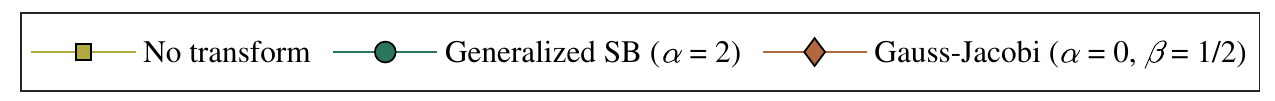}\\[0.1in]
  \begin{subfigure}{2in}
    \caption*{}
  \end{subfigure}
  \begin{subfigure}{2in}
    \centering
    $f_{S1} (\vx)$
    \caption*{}
  \end{subfigure}
  \begin{subfigure}{2in}
    \centering
    $f_{S2} (\vx)$
    \caption*{}
  \end{subfigure}\\[-0.2in]
  \begin{subfigure}{2in}
    \centering
    \includegraphics[scale=1]{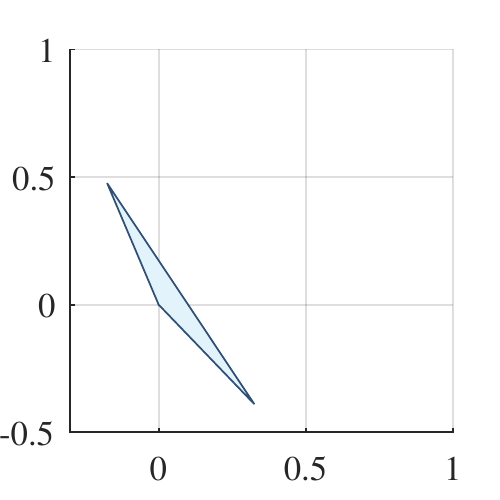}
    \caption{}\label{fig:ex-singular-xi1-shape2}
  \end{subfigure}
  \begin{subfigure}{2in}
    \centering
    \includegraphics[scale=1]{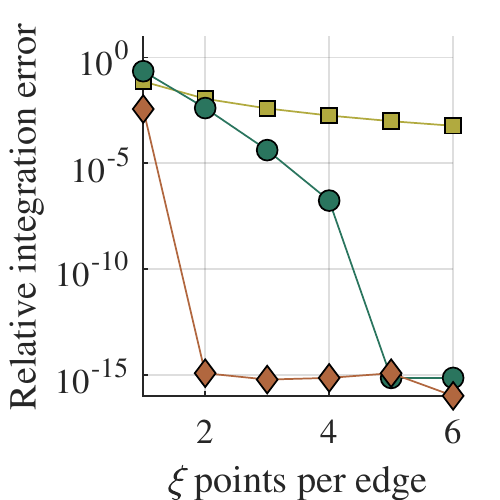}
    \caption{}\label{fig:ex-singular-xi1-conv-1-1}
  \end{subfigure}
  \begin{subfigure}{2in}
    \centering
    \includegraphics[scale=1]{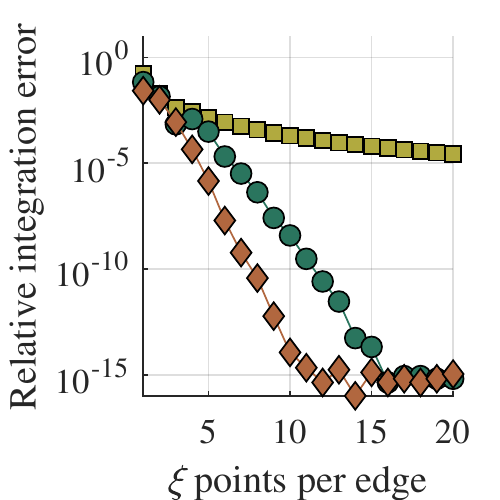}
    \caption{}\label{fig:ex-singular-xi1-conv-1-2}
  \end{subfigure}
  \begin{subfigure}{2in}
    \centering
    \includegraphics[scale=1]{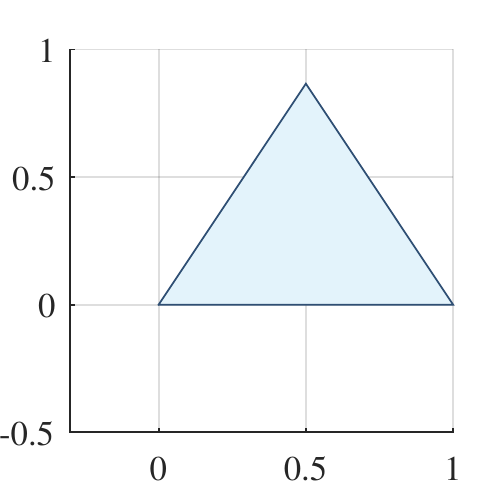}
    \caption{}\label{fig:ex-singular-xi1-shape3}
  \end{subfigure}
  \begin{subfigure}{2in}
    \centering
    \includegraphics[scale=1]{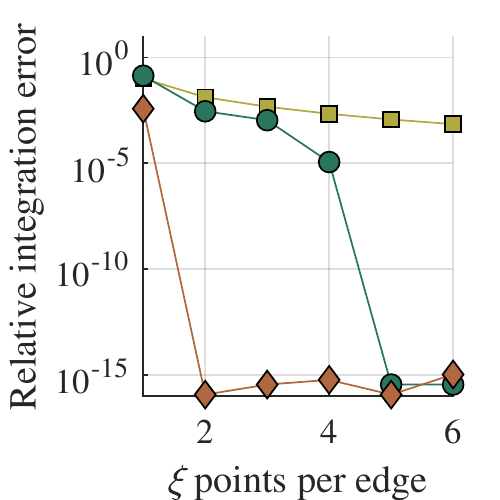}
    \caption{}\label{fig:ex-singular-xi1-conv-2-1}
  \end{subfigure}
  \begin{subfigure}{2in}
    \centering
    \includegraphics[scale=1]{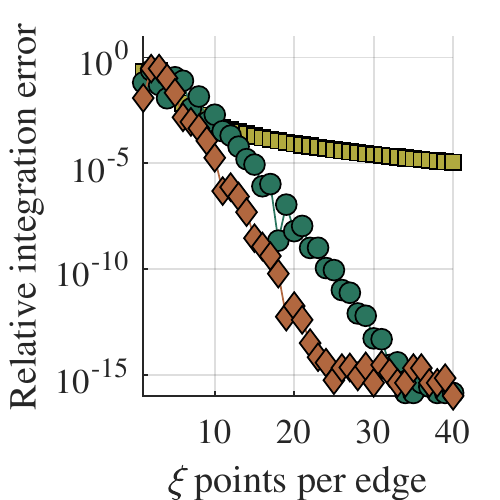}
    \caption{}\label{fig:ex-singular-xi1-conv-2-2}
  \end{subfigure}
  \begin{subfigure}{2in}
    \centering
    \includegraphics[scale=1]{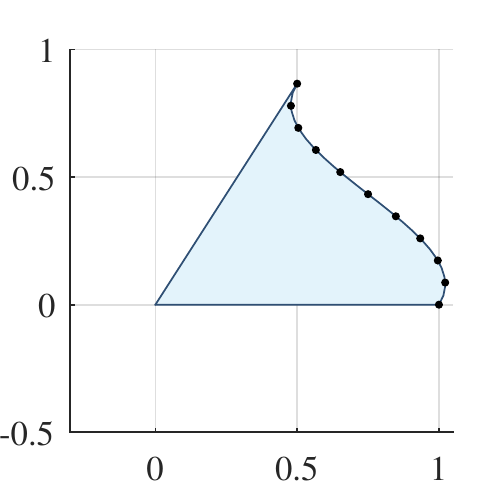}
    \caption{}\label{fig:ex-singular-xi1-shape4}
  \end{subfigure}
  \begin{subfigure}{2in}
    \centering
    \includegraphics[scale=1]{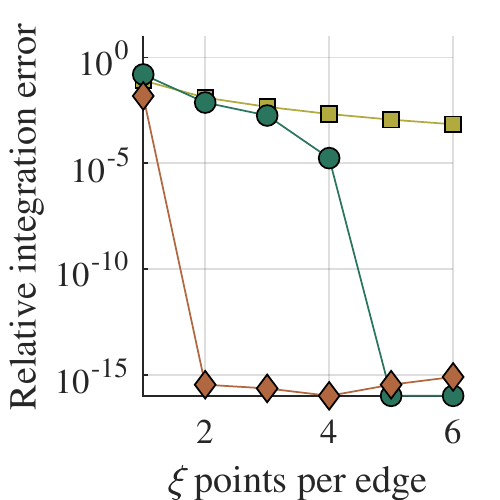}
    \caption{}\label{fig:ex-singular-xi1-conv-3-1}
  \end{subfigure}
  \begin{subfigure}{2in}
    \centering
    \includegraphics[scale=1]{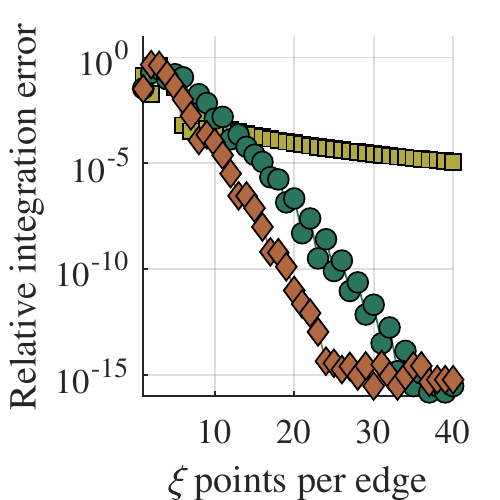}
    \caption{}\label{fig:ex-singular-xi1-conv-3-2}
  \end{subfigure}
  \caption{Number of $\xi$ quadrature points versus integration accuracy to
      integrate $f_{S1} (\vx)$ over (b) $T_2$, (e) $T_3$, and (h) $T_4$ and to
      integrate $f_{S2} (\vx)$ over (c) $T_2$, (f) $T_3$, and (i) $T_4$.  The
      domain of (a) $T_2$, (d) $T_3$, and (g) $T_4$.}
  \label{fig:ex-singular-xi1}
\end{figure}

\begin{figure}
  \centering
  \includegraphics[scale=1]{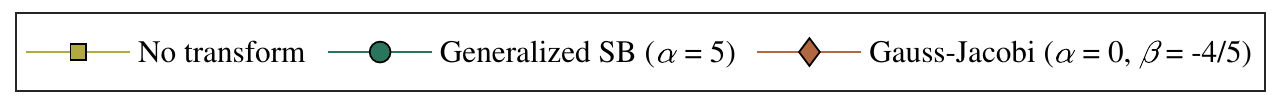}\\[0.1in]
  \begin{subfigure}{2in}
    \caption*{}
  \end{subfigure}
  \begin{subfigure}{2in}
    \centering
    $f_{S3} (\vx)$
    \caption*{}
  \end{subfigure}
  \begin{subfigure}{2in}
    \centering
    $f_{S4} (\vx)$
    \caption*{}
  \end{subfigure}\\[-0.2in]
  \begin{subfigure}{2in}
    \centering
    \includegraphics[scale=1]{ex-singular-shape2}
    \caption{}\label{fig:ex-singular-xi2-shape2}
  \end{subfigure}
  \begin{subfigure}{2in}
    \centering
    \includegraphics[scale=1]{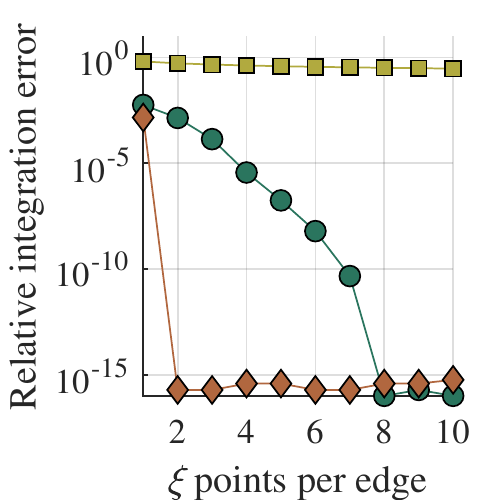}
    \caption{}\label{fig:ex-singular-xi2-conv-1-1}
  \end{subfigure}
  \begin{subfigure}{2in}
    \centering
    \includegraphics[scale=1]{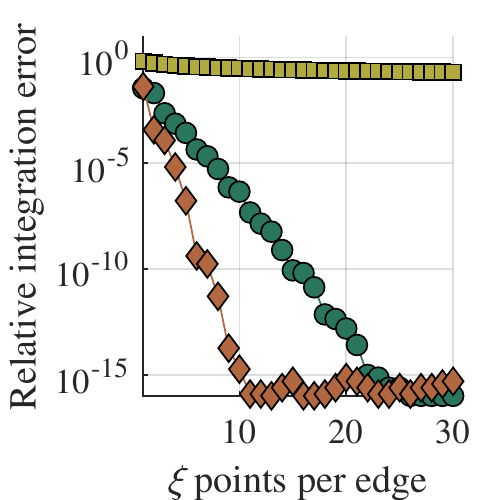}
    \caption{}\label{fig:ex-singular-xi2-conv-1-2}
  \end{subfigure}
  \begin{subfigure}{2in}
    \centering
    \includegraphics[scale=1]{ex-singular-shape3}
    \caption{}\label{fig:ex-singular-xi2-shape3}
  \end{subfigure}
  \begin{subfigure}{2in}
    \centering
    \includegraphics[scale=1]{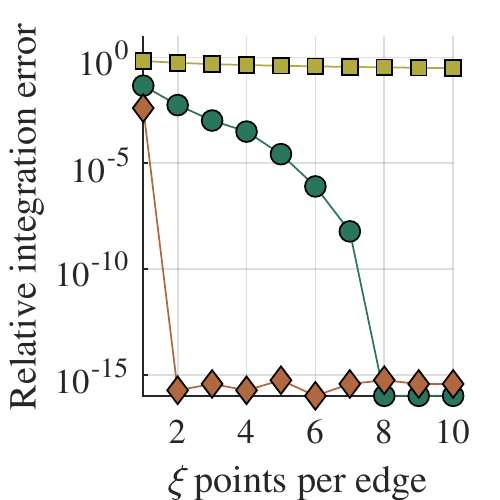}
    \caption{}\label{fig:ex-singular-xi2-conv-2-1}
  \end{subfigure}
  \begin{subfigure}{2in}
    \centering
    \includegraphics[scale=1]{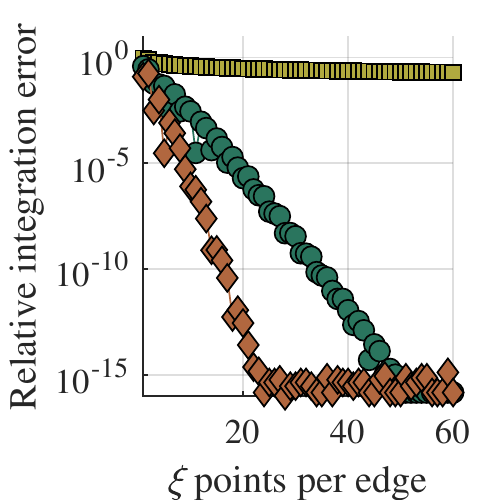}
    \caption{}\label{fig:ex-singular-xi2-conv-2-2}
  \end{subfigure}
  \begin{subfigure}{2in}
    \centering
    \includegraphics[scale=1]{ex-singular-shape4}
    \caption{}\label{fig:ex-singular-xi2-shape4}
  \end{subfigure}
  \begin{subfigure}{2in}
    \centering
    \includegraphics[scale=1]{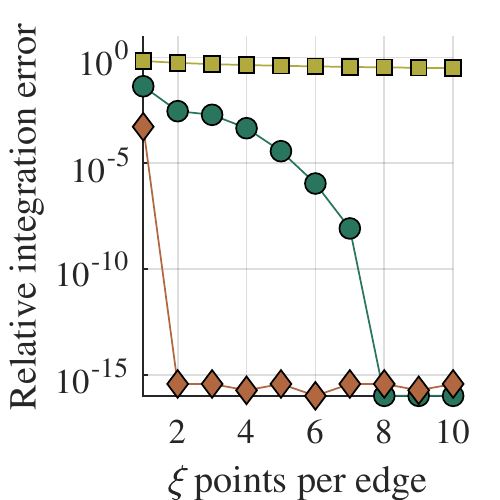}
    \caption{}\label{fig:ex-singular-xi2-conv-3-1}
  \end{subfigure}
  \begin{subfigure}{2in}
    \centering
    \includegraphics[scale=1]{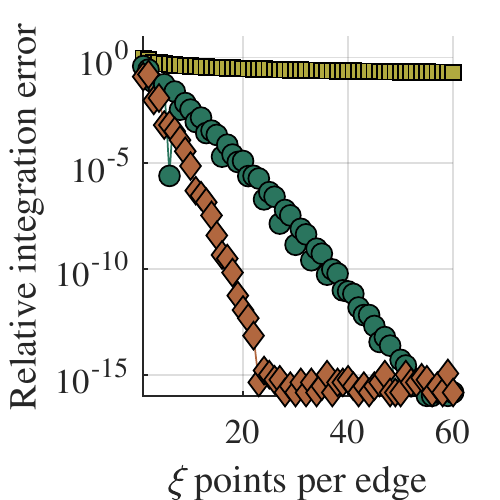}
    \caption{}\label{fig:ex-singular-xi2-conv-3-2}
  \end{subfigure}
  \caption{Number of $\xi$ quadrature points versus integration accuracy to
  integrate $f_{S3} (\vx)$ over (b) $T_2$, (e) $T_3$, and (h) $T_4$ and to
  integrate $f_{S4} (\vx)$ over (c) $T_2$, (f) $T_3$, and (i) $T_4$.  The
  domain of (a) $T_2$, (d) $T_3$, and (g) $T_4$.}
  \label{fig:ex-singular-xi2}
\end{figure}

\smallskip
Each of these functions is integrated over two triangles ($T_2$ and $T_3$) and
one curved triangle ($T_4$) whose boundary is defined by a cubic B\'{e}zier
curve.  The vertex locations of $T_2$ and $T_3$ are provided
in~\tref{tab:ex-singular-triangles}.  $T_2$ represents a triangle with a high
base to height ratio and $T_3$ is a triangle where the base to height ratio is
close to unity.  The control points for the B\'{e}zier curve in $T_4$ are listed
in~\tref{tab:ex-singular-curved}.  The point $\vx_c = \vm{0}$ and the endpoints
of the curve define the two line segments that bound $T_4$.  These three shapes
are illustrated in~\fref{fig:ex-singular-xi1} and \fref{fig:ex-singular-xi2}. To
analyze the effect of the integration techniques in the 
$\xi$ direction, a
surplus of cubature points is provided in the 
$t$ direction ensuring
integration error above machine precision is only present in the
$\xi$ direction.  Exact integrals are computed using Mathematica 12.0.0.

{\renewcommand{\arraystretch}{2.1}
\begin{table}
  \centering
  \caption{Vertices of the triangles over which \revone{weakly} singular integrals 
  are computed.}
  \label{tab:ex-singular-triangles}
  \begin{tabular}{|c|c|c|c|}
    \hline
    Triangle & Vertex 1 & Vertex 2 & Vertex 3 \\[0.05in]
    \hline
    $T_1$ & $(0, 0)$ & $\left( \dfrac{103}{400}, -\dfrac{99 \sqrt{3}}{400} \right)$
      & $\left( -\dfrac{97}{400}, \dfrac{101 \sqrt{3}}{400} \right)$ \\[0.05in]
    \hline
    $T_2$ & $(0, 0)$ & $\left( \dfrac{13}{40}, -\dfrac{9 \sqrt{3}}{40} \right)$
      & $\left( -\dfrac{7}{40}, \dfrac{11 \sqrt{3}}{40} \right)$ \\[0.05in]
    \hline
    $T_3$ & $(0, 0)$ & $\left( 1, 0 \right)$
      & $\left( \dfrac{1}{2}, \dfrac{\sqrt{3}}{2} \right)$ \\[0.05in]
    \hline
  \end{tabular}
\end{table}
}

{\renewcommand{\arraystretch}{2.1}
\begin{table}
  \centering
  \caption{Control points of the curved triangle over which \revone{weakly} singular 
  integrals are computed.}
  \label{tab:ex-singular-curved}
  \begin{tabular}{|c|c|c|c|c|}
    \hline
    Curved Triangle & Point 1 & Point 2 & Point 3 & Point 4 \\[0.05in]
    \hline
    $T_4$ & $\left(1, 0\right)$ & $\left(\dfrac{7}{6}, \dfrac{\sqrt{3}}{6}\right)$ 
      & $\left(\dfrac{1}{3}, \dfrac{\sqrt{3}}{3} \right)$ 
      & $\left(\dfrac{1}{2}, \dfrac{\sqrt{3}}{2}\right)$ \\[0.05in]
    \hline
  \end{tabular}
\end{table}
}

\smallskip
Results for the two functions with the $r^{-1/2}$ singularity are presented
in~\fref{fig:ex-singular-xi1}.  For both the function with polynomial numerator
($f_{S1} (\vx)$) and the \revone{one with}
non\hyp{}polynomial numerator ($f_{S2} (\vx)$), a
Gauss\hyp{}Jacobi rule with $\alpha = 0$ and $\beta = 1/2$ provides the most
accurate integration per quadrature point.  To integrate $f_{S1} (\vx)$ to about
15 digits of precision, two quadrature points are needed in the
$\xi$ direction over all three shapes.  For $\mathcal{O}(10^{-15})$
accurate integration of $f_{S2} (\vx)$, 12 quadrature points are needed in the
$\xi$ direction for $T_2$ and approximately 25 quadrature points are needed
for $T_3$ and $T_4$.  The generalized SB transformation with $\alpha = 2$ also
enables accurate cubature of $f_{S1} (\vx)$ and $f_{S2} (\vx)$, though not as
efficiently as the Gauss\hyp{}Jacobi rule.  Machine precision integration (about
15 digits of precision) of $f_{S1} (\vx)$ over the three domains requires 5
points in the $\xi$ direction.  For machine precision integration of the
non\hyp{}polynomial integrand ($f_{S2} (\vx)$), 16 quadrature points are needed
over $T_2$ and about 35 quadrature points are needed over $T_3$ and $T_4$.  When
SBC is used with no transformation in the $\xi$ direction, no more than six
digits of precision can be obtained with 64 $\xi$\hyp{}quadrature points for
either $f_{S1} (\vx)$ or $f_{S2} (\vx)$ over any of the three domains.

\smallskip
Integration error versus number of quadrature points in the 
$\xi$ direction
for integrating $f_{S3} (\vx)$ and $f_{S4} (\vx)$ are illustrated
in~\fref{fig:ex-singular-xi2}.  Qualitatively, the results of integrating the
two functions with a $r^{-9/5}$ singularity mirror those for the $r^{-1/2}$
singularity.  For the Gauss\hyp{}Jacobi rule with $\alpha = 0$ and $\beta =
-4/5$, accuracy per quadrature point is almost identical to the results
presented in~\fref{fig:ex-singular-xi1}, though slightly fewer quadrature points
are needed to integrate $f_{S4} (\vx)$ to machine precision as compared to
$f_{S2} (\vx)$.  For these integrands, the generalized SB transformation with
$\alpha = 5$ is less efficient than the $\alpha = 2$ transformation used to
integrate $f_{S1} (\vx)$ and $f_{S2} (\vx)$.  Eight quadrature points are needed
in the $\xi$ direction to integrate $f_{S3} (\vx)$ to machine precision and
24 (resp.\ 50) quadrature points are needed to integrate $f_{S4} (\vx)$ to 15
digits of accuracy on $T_2$ (resp.\ on $T_3$ and $T_4$).

\subsubsection{Comparison of $t$\hyp{}transformation methods}\label{sec:ex-singular-t}
This section compares the three integral transformations introduced
in~\sref{sec:singular-t} that are designed to smooth near\hyp{}singularities in
the $t$ direction.  Functions with $r^{-1/2}$, $r^{-1}$, and $r^{-9/5}$
singularities are integrated.  The functions with the $r^{-1/2}$ and $r^{-9/5}$
singularities are given in~\eqref{eq:singular-fns-xi} and the functions with the
$r^{-1}$ singularity are
\begin{subequations}
\begin{align}\label{eq:f-s5}
  f_{S5} (\vx) & = \frac{2x^2 + 2y \left( y + \sqrt{x^2 + y^2} \right)
    + x \left( y + 2\sqrt{x^2 + y^2} \right)}{\left( x^2 + y^2 \right)^{3/2}}, \\
  \label{eq:f-s6}
  f_{S6} (\vx) & = \frac{\exp \left( -\left[ \left(\frac{x - 0.25}{0.4} \right)^2
  + \left( \frac{y - 0.2}{0.7} \right)^2 \right]^2 \right)
  \, \cos^2 (5x) \, \cos^2 (5y)}
  {\sqrt{x^2 + y^2}} .
\end{align}
\end{subequations}
Equation~\eqref{eq:f-s5} is the sum of the terms that appear in the kernel for
linear elasticity and elasto\hyp{}plasticity in the
BEM~\cite{Nagarajan:1993:AMM}.  Equation~\eqref{eq:f-s6} is analogous
to~\eqref{eq:f-s2} and~\eqref{eq:f-s4}, though the singularity is of the form
$r^{-1}$.  The accuracy of each integral transformation is tested by integrating
in the $\xi$ direction to $\mathcal{O}(10^{-15})$ accuracy, then increasing
the number of cubature points in the $t$ direction.  
The exact integral values are computed using Mathematica 12.0.0.

\smallskip
These six functions are integrated over three different triangles: $T_1$, $T_2$,
and $T_3$, which correspond to triangles with very high, high, and near unitary
base to height ratios, respectively.  The vertices of these triangles are listed
in~\tref{tab:ex-singular-triangles}.  No triangles bounded by non\hyp{}affine
curves are included since the integral transformations in the 
$t$ direction are
formulated for domains bounded by line segments.  Since polygons are decomposed
into triangles by the SBC method, the results in this section (integration over
triangles) can be generalized to polygonal domains of integration.

\begin{figure}
  \centering
  \includegraphics[scale=1]{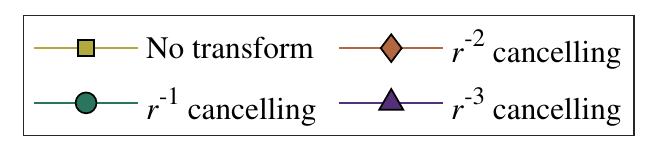}\\[0.1in]
  \begin{subfigure}{2in}
    \caption*{}
  \end{subfigure}
  \begin{subfigure}{2in}
    \centering
    $f_{S1} (\vx)$
    \caption*{}
  \end{subfigure}
  \begin{subfigure}{2in}
    \centering
    $f_{S2} (\vx)$
    \caption*{}
  \end{subfigure}\\[-0.2in]
  \begin{subfigure}{2in}
    \centering
    \includegraphics[scale=1]{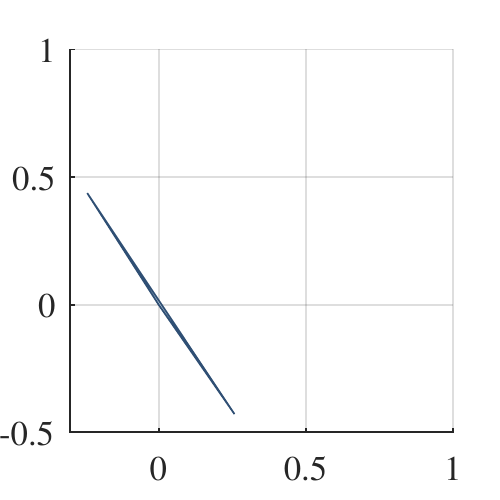}
    \caption{}\label{fig:ex-singular-t1-shape1}
  \end{subfigure}
  \begin{subfigure}{2in}
    \centering
    \includegraphics[scale=1]{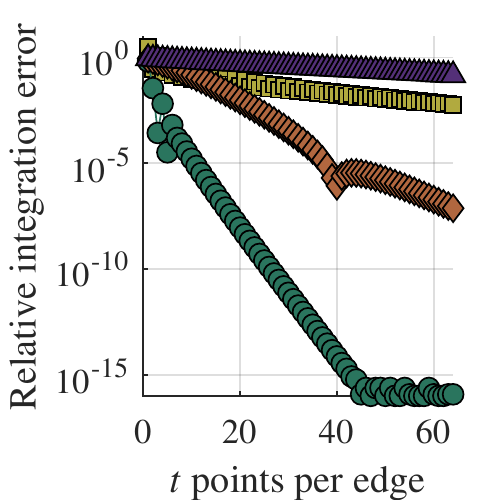}
    \caption{}\label{fig:ex-singular-t1-conv-1-1}
  \end{subfigure}
  \begin{subfigure}{2in}
    \centering
    \includegraphics[scale=1]{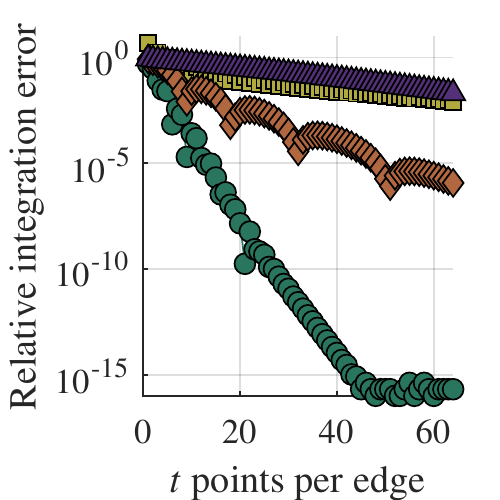}
    \caption{}\label{fig:ex-singular-t1-conv-1-2}
  \end{subfigure}
  \begin{subfigure}{2in}
    \centering
    \includegraphics[scale=1]{ex-singular-shape2}
    \caption{}\label{fig:ex-singular-t1-shape2}
  \end{subfigure}
  \begin{subfigure}{2in}
    \centering
    \includegraphics[scale=1]{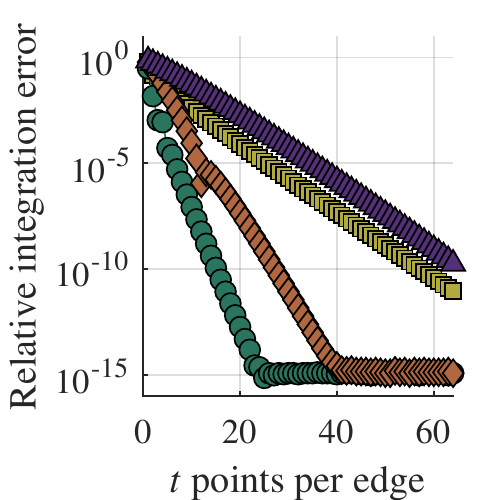}
    \caption{}\label{fig:ex-singular-t1-conv-2-1}
  \end{subfigure}
  \begin{subfigure}{2in}
    \centering
    \includegraphics[scale=1]{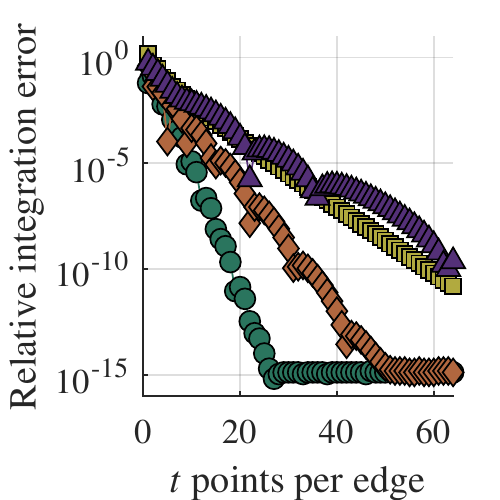}
    \caption{}\label{fig:ex-singular-t1-conv-2-2}
  \end{subfigure}
  \begin{subfigure}{2in}
    \centering
    \includegraphics[scale=1]{ex-singular-shape3}
    \caption{}\label{fig:ex-singular-t1-shape3}
  \end{subfigure}
  \begin{subfigure}{2in}
    \centering
    \includegraphics[scale=1]{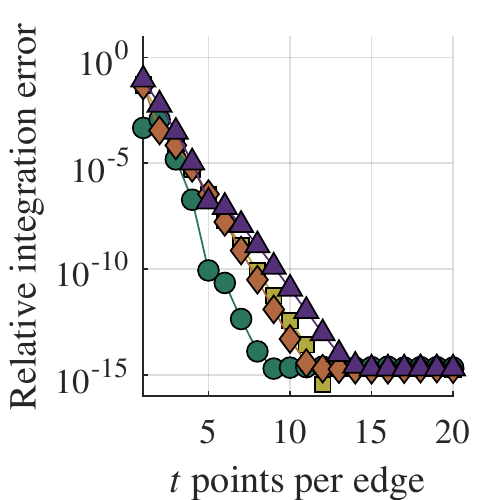}
    \caption{}\label{fig:ex-singular-t1-conv-3-1}
  \end{subfigure}
  \begin{subfigure}{2in}
    \centering
    \includegraphics[scale=1]{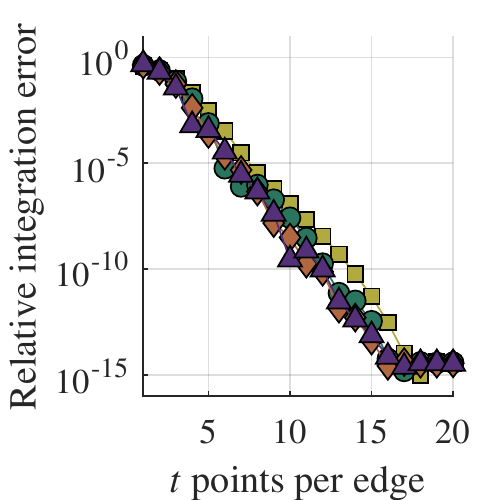}
    \caption{}\label{fig:ex-singular-t1-conv-3-2}
  \end{subfigure}
  \caption{Number of $t$ quadrature points versus integration accuracy to
  integrate $f_{S1} (\vx)$ over (b) $T_1$, (e) $T_2$, and (h) $T_3$ and to
  integrate $f_{S2} (\vx)$ over (c) $T_1$, (f) $T_2$, and (i) $T_3$.  The
  domain of (a) $T_1$, (d) $T_2$, and (g) $T_3$.}
  \label{fig:ex-singular-t1}
\end{figure}

\begin{figure}
  \centering
  \includegraphics[scale=1]{ex-singular-t-legend}\\[0.1in]
  \begin{subfigure}{2in}
    \caption*{}
  \end{subfigure}
  \begin{subfigure}{2in}
    \centering
    $f_{S5} (\vx)$
    \caption*{}
  \end{subfigure}
  \begin{subfigure}{2in}
    \centering
    $f_{S6} (\vx)$
    \caption*{}
  \end{subfigure}\\[-0.2in]
  \begin{subfigure}{2in}
    \centering
    \includegraphics[scale=1]{ex-singular-shape1}
    \caption{}\label{fig:ex-singular-t2-shape1}
  \end{subfigure}
  \begin{subfigure}{2in}
    \centering
    \includegraphics[scale=1]{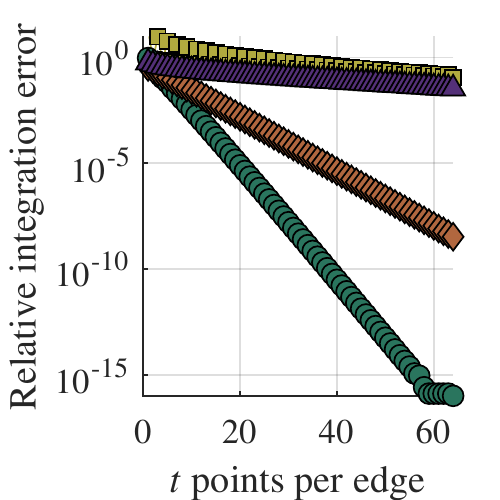}
    \caption{}\label{fig:ex-singular-t2-conv-1-1}
  \end{subfigure}
  \begin{subfigure}{2in}
    \centering
    \includegraphics[scale=1]{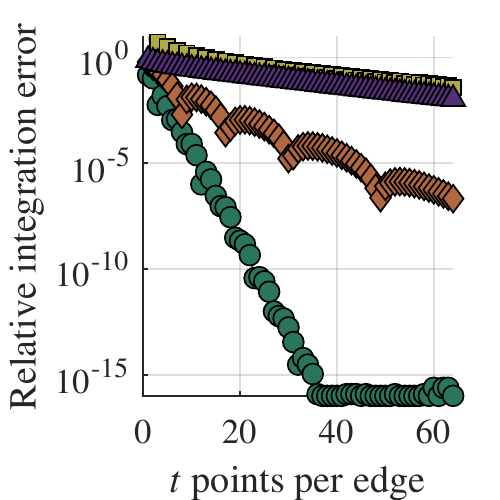}
    \caption{}\label{fig:ex-singular-t2-conv-1-2}
  \end{subfigure}
  \begin{subfigure}{2in}
    \centering
    \includegraphics[scale=1]{ex-singular-shape2}
    \caption{}\label{fig:ex-singular-t2-shape2}
  \end{subfigure}
  \begin{subfigure}{2in}
    \centering
    \includegraphics[scale=1]{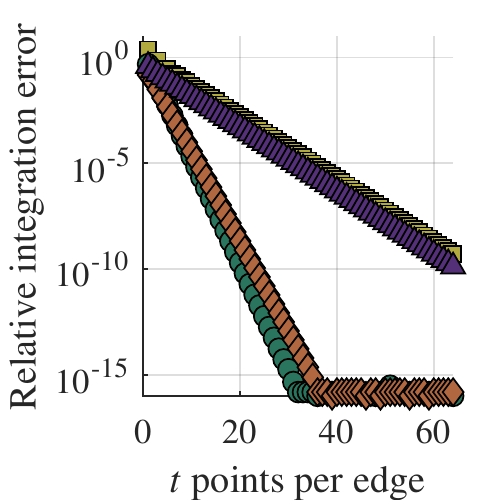}
    \caption{}\label{fig:ex-singular-t2-conv-2-1}
  \end{subfigure}
  \begin{subfigure}{2in}
    \centering
    \includegraphics[scale=1]{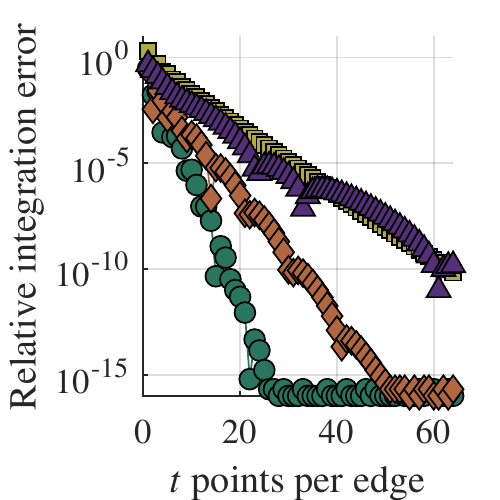}
    \caption{}\label{fig:ex-singular-t2-conv-2-2}
  \end{subfigure}
  \begin{subfigure}{2in}
    \centering
    \includegraphics[scale=1]{ex-singular-shape3}
    \caption{}\label{fig:ex-singular-t2-shape3}
  \end{subfigure}
  \begin{subfigure}{2in}
    \centering
    \includegraphics[scale=1]{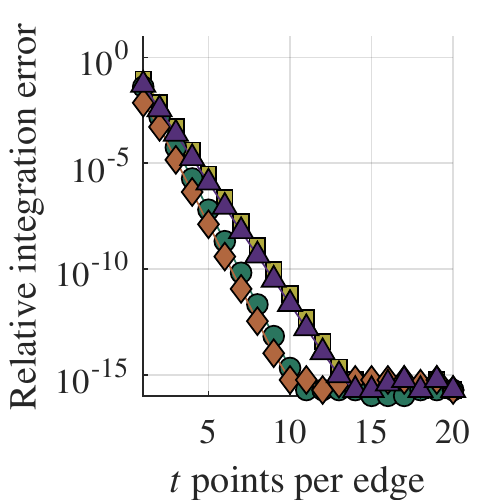}
    \caption{}\label{fig:ex-singular-t2-conv-3-1}
  \end{subfigure}
  \begin{subfigure}{2in}
    \centering
    \includegraphics[scale=1]{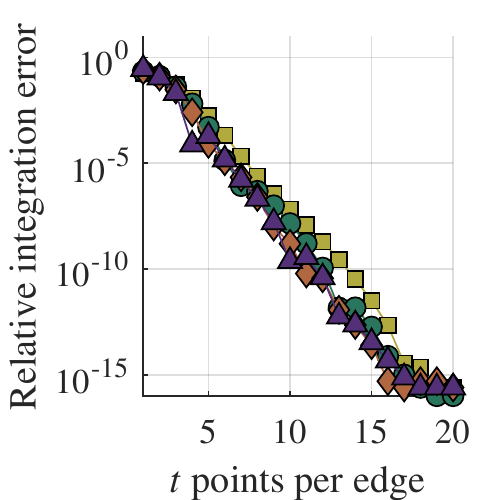}
    \caption{}\label{fig:ex-singular-t2-conv-3-2}
  \end{subfigure}
  \caption{Number of $t$ quadrature points versus integration accuracy to
  integrate $f_{S5} (\vx)$ over (b) $T_1$, (e) $T_2$, and (h) $T_3$ and to
  integrate $f_{S6} (\vx)$ over (c) $T_1$, (f) $T_2$, and (i) $T_3$.  The
  domain of (a) $T_1$, (d) $T_2$, and (g) $T_3$.}
  \label{fig:ex-singular-t2}
\end{figure}

\begin{figure}
  \centering
  \includegraphics[scale=1]{ex-singular-t-legend}\\[0.1in]
  \begin{subfigure}{2in}
    \caption*{}
  \end{subfigure}
  \begin{subfigure}{2in}
    \centering
    $f_{S3} (\vx)$
    \caption*{}
  \end{subfigure}
  \begin{subfigure}{2in}
    \centering
    $f_{S4} (\vx)$
    \caption*{}
  \end{subfigure}\\[-0.2in]
  \begin{subfigure}{2in}
    \centering
    \includegraphics[scale=1]{ex-singular-shape1}
    \caption{}\label{fig:ex-singular-t3-shape1}
  \end{subfigure}
  \begin{subfigure}{2in}
    \centering
    \includegraphics[scale=1]{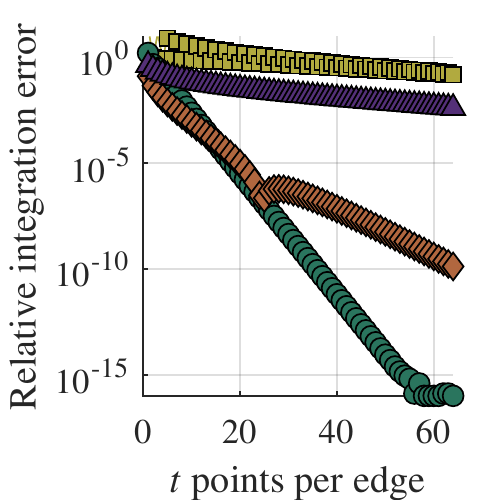}
    \caption{}\label{fig:ex-singular-t3-conv-1-1}
  \end{subfigure}
  \begin{subfigure}{2in}
    \centering
    \includegraphics[scale=1]{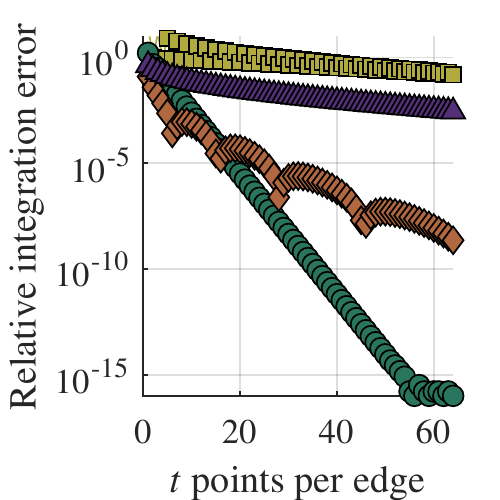}
    \caption{}\label{fig:ex-singular-t3-conv-1-2}
  \end{subfigure}
  \begin{subfigure}{2in}
    \centering
    \includegraphics[scale=1]{ex-singular-shape2}
    \caption{}\label{fig:ex-singular-t3-shape2}
  \end{subfigure}
  \begin{subfigure}{2in}
    \centering
    \includegraphics[scale=1]{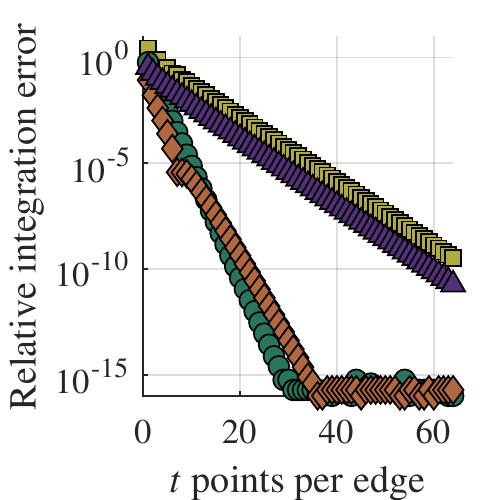}
    \caption{}\label{fig:ex-singular-t3-conv-2-1}
  \end{subfigure}
  \begin{subfigure}{2in}
    \centering
    \includegraphics[scale=1]{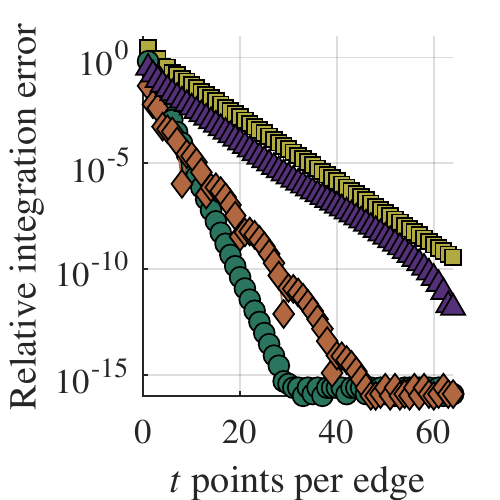}
    \caption{}\label{fig:ex-singular-t3-conv-2-2}
  \end{subfigure}
  \begin{subfigure}{2in}
    \centering
    \includegraphics[scale=1]{ex-singular-shape3}
    \caption{}\label{fig:ex-singular-t3-shape3}
  \end{subfigure}
  \begin{subfigure}{2in}
    \centering
    \includegraphics[scale=1]{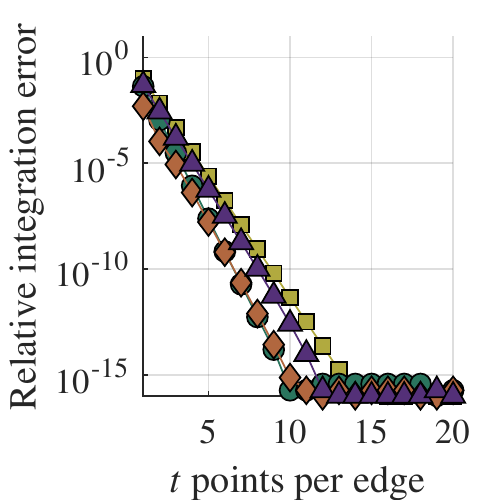}
    \caption{}\label{fig:ex-singular-t3-conv-3-1}
  \end{subfigure}
  \begin{subfigure}{2in}
    \centering
    \includegraphics[scale=1]{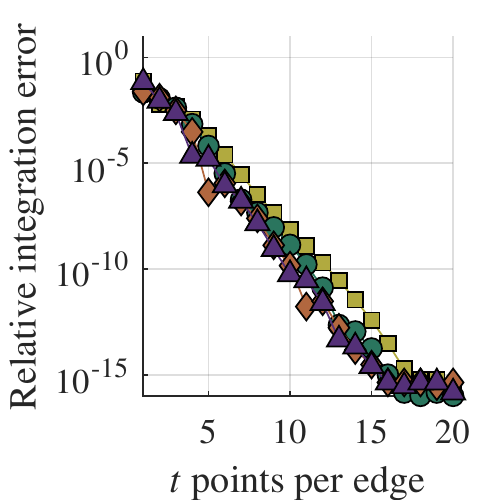}
    \caption{}\label{fig:ex-singular-t3-conv-3-2}
  \end{subfigure}
  \caption{Number of $t$ quadrature points versus integration accuracy to
  integrate $f_{S3} (\vx)$ over (b) $T_1$, (e) $T_2$, and (h) $T_3$ and to
  integrate $f_{S4} (\vx)$ over (c) $T_1$, (f) $T_2$, and (i) $T_3$.  The
  domain of (a) $T_1$, (d) $T_2$, and (g) $T_3$.}
  \label{fig:ex-singular-t3}
\end{figure}

\smallskip
Integration accuracy versus the number of $t$\hyp{}direction points is plotted
in Figs.~\ref{fig:ex-singular-t1}, \ref{fig:ex-singular-t2},
and~\ref{fig:ex-singular-t3} for functions containing $r^{-1/2}$, $r^{-1}$, and
$r^{-9/5}$ singularities, respectively.  Only the $r^{-1}$\hyp{}cancelling
transformation is able to consistently provide about 15 digits of precision in
integration with fewer than 60 quadrature points in the $t$ direction.  As the
strength of the singularity increases, the performance of the
$r^{-2}$\hyp{}cancelling transformation improves.  This mirrors results
presented in Chin~\cite{Chin:2019:NIH}, where the $r^{-2}$\hyp{}cancelling
transformation outperforms the $r^{-1}$\hyp{}cancelling transformation as the
power of the singularity becomes more negative.  For integration of the
weakly singular integrands presented in this section, the lowest accuracy
per quadrature point is consistently obtained using either no transformation or
the $r^{-3}$\hyp{}cancelling transformation.  Stronger singularities are likely
needed for the $r^{-3}$ transformation to be competitive versus the $r^{-1}$ and
$r^{-2}$ transformations.

\smallskip
As the base to height ratio of the triangular integration domain approaches
unity, the efficacy of the $r^{-1}$\hyp{} and $r^{-2}$\hyp{}cancelling
transformations wanes.  When integrating the non\hyp{}polynomial integrands
($f_{S2} (\vx)$, $f_{S6} (\vx)$, and $f_{S4} (\vx)$) over $T_3$, the benefit of
the transformation is minimal---as compared to SBC with no
$t$\hyp{}transformation, only two fewer quadrature points are needed to
integrate the non\hyp{}polynomial integrands to 15 digits of precision when the
$r^{-1}$\hyp{}cancelling transformation is applied. As the base to height ratio
of the triangular integration domain increases, the near\hyp{}singularity in the
integrand is less pronounced, and consequently, the number of quadrature points
needed for accurate integration in the $t$ direction decreases.  This
suggests the number of quadrature points in the $t$ direction can be scaled
based on the base to height ratio of a triangle for a more efficient cubature
rule.

\subsubsection{Integrals in the extended finite element method}
Modeling cracks in linear elastic continua using the X-FEM requires integration
of weakly singular functions.  In the X-FEM, one component of the stiffness
matrix calculations is of the form:
\begin{equation}\label{eq:xfem-k}
  K_{IJ} = \int_{\Omega} \left[ \frac{\partial 
    \Bigl(F(\vx) N_I \bigl(\vm{\xi}(\vx)\bigr)\Bigr)}{\partial x}
    \frac{\partial N_J \bigl(\vm{\xi}(\vx)\bigr)}{\partial x} + \frac{\partial 
    \Bigl(F(\vx) N_I \bigl(\vm{\xi}(\vx)\bigr)\Bigr)}{\partial y}
    \frac{\partial N_J \bigl(\vm{\xi}(\vx)\bigr)}{\partial y} \right] d\vx ,
\end{equation}
where $F(\vx)$ is an enrichment function used to capture singular behavior in
the vicinity of a crack tip and $N_I (\vm{\xi})$ is a finite element basis
function.  Typically, four crack\hyp{}tip enrichment functions are used, but 
in this study, only the function
\begin{subequations}
\begin{equation}
  F(\vx) = \sqrt{r} \sin \left( \frac{1}{2} \theta \right) ,
\end{equation}
where
\begin{equation}
  r = \sqrt{x^2 + y^2} \qquad \text{and} \qquad
  \theta = \arctan \frac{y}{x} \quad \forall \theta \in [-\pi, \pi] ,
\end{equation}
\end{subequations}
is considered.  Note the derivative of $F(\vx)$ is a weakly singular
function with a $r^{-1/2}$ singularity, and furthermore, $F(\vx)$ contains a
discontinuity at $\theta = \pm \pi$.  In the X-FEM, the location of the
discontinuity coincides with the location of a crack in the domain. Bilinear
quadrilaterals with four nodes are used to form the finite element basis.  In
the isoparametric domain $\Xi = [-1, 1]^2$, these shape functions are
\begin{equation}
  N_1(\vm{\xi}) = \frac{(1-\xi)(1-\eta)}{4}, \quad 
  N_2(\vm{\xi}) = \frac{(1+\xi)(1-\eta)}{4}, \quad 
  N_3(\vm{\xi}) = \frac{(1+\xi)(1+\eta)}{4}, \quad 
  N_4(\vm{\xi}) = \frac{(1-\xi)(1+\eta)}{4} .
\end{equation}
Mapping between an arbitrary quadrilateral and $\Xi$ is performed using the
isoparametric mapping
\begin{equation}
  \vx(\vm{\xi}) = \sum_{I=1}^4 N_I (\vm{\xi}) \vx_I ,
\end{equation}
where $\vx_I$ is the nodal coordinate associated with $N_I (\vm{\xi})$.

\begin{figure}
  \centering
  \begin{subfigure}{3in}
    \centering
    \begin{tikzpicture}
      \node at (0in, 0in) {\includegraphics[scale=1]{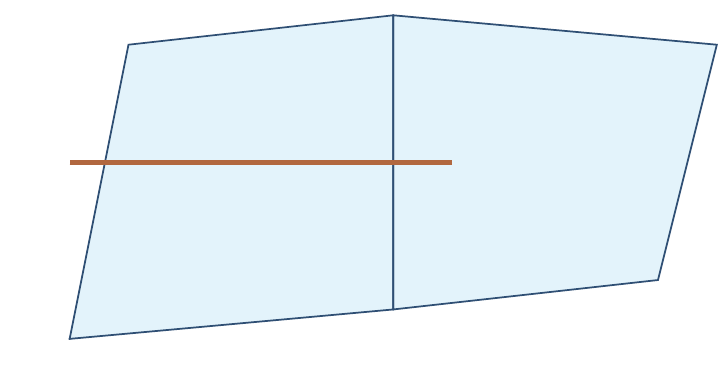}};
      \node at (0.36in, 0.18in) {$\vx_c$};
      \node at (0.24in, -0.05in) {$\Delta x$};
      \draw (0.355in, 0.1in) -- (0.355in, 0in);
      \draw (0.12in, 0.1in) -- (0.12in, 0in);
      \draw (0.12in, 0.03in) -- (0.355in, 0.03in);
      \node at (-0.5in, -0.1in) {$\Omega_1$};
      \node at (0.75in, 0.075in) {$\Omega_2$};
      \node at (-1in, -0.7in) {$(-1.1, -0.1)$};
      \node at (0.1in, -0.6in) {$(0, 0)$};
      \node at (1.05in, -0.5in) {$(0.9, 0.1)$};
      \node at (-0.95in, 0.675in) {$(-0.9, 0.9)$};
      \node at (0.1in, 0.775in) {$(0, 1)$};
      \node at (1.2in, 0.675in) {$(1.1, 0.9)$};
      \node at (0, 1in) {};
      \node at (0, -0.92in) {};
    \end{tikzpicture}
    \caption{}\label{fig:ex-xfem-domain}
  \end{subfigure}
  \begin{subfigure}{3in}
    \centering
    \begin{tikzpicture}
      \draw[pal58, fill=pal29] (0,0) rectangle (0.2in,0.2in);
      \draw[pal58, pattern=vertical lines, pattern color=black!40!white] (0,0) rectangle (0.2in,0.2in);
      \node[anchor=west] at (0.25in,0.09in) {Positive integral};
      \draw[pal58, fill=pal29] (1.5in,0) rectangle (1.7in,0.2in);
      \draw[pal58, pattern=horizontal lines, pattern color=black!40!white] (1.5in,0) rectangle (1.7in,0.2in);
      \node[anchor=west] at (1.75in,0.09in) {Negative integral};
      \draw (-0.1in,0.23in) rectangle (2.8in,-0.03in);
    \end{tikzpicture}\\
    \begin{subfigure}{1.48in}
      \centering
      \includegraphics[scale=1]{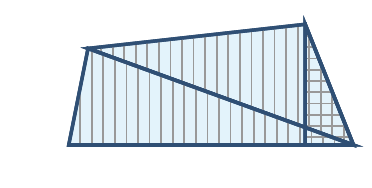}
      \includegraphics[scale=1]{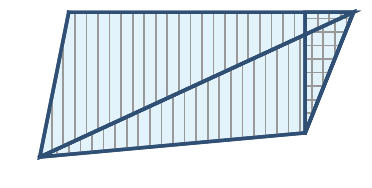}
      \caption*{}
    \end{subfigure}
    \begin{subfigure}{1.48in}
      \centering
      \includegraphics[scale=1]{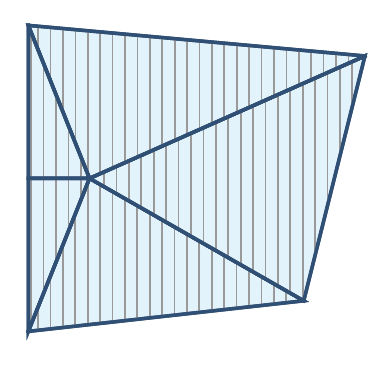}
      \caption*{}
    \end{subfigure}
    \caption{}\label{fig:ex-xfem-decomp}
  \end{subfigure}
  \caption{Problem setup for numerically integrating functions that appear in
  the X-FEM.  (a) Problem domain with an integrand discontinuity (thick line)
  and (b) domain partitioning after applying SBC to regions with continuous
  integrands.}
  \label{fig:ex-xfem-problem}
\end{figure}

\smallskip
For the bilinear quadrilateral element, \eqref{eq:xfem-k} defines 16 different
integrals.  In this example, these 16 integrals are computed over two element
domains: $\Omega_1$ and $\Omega_2$ pictured in~\fref{fig:ex-xfem-domain}.  In
the element domains, the nodal locations are perturbed such that the integrands
are non\hyp{}polynomial functions.  This example is presented to demonstrate a
realistic application of singularity\hyp{}cancelling methods applicable to the
SBC method described in~\sref{sec:singular}.  Three different crack\hyp{}tip
locations are tested with varying proximity to the element interface.  The
distances to the interface in the three different crack\hyp{}tip locations are
$\Delta x = 0.001$, $\Delta x = 0.01$, and $\Delta x = 0.1$.  The length $\Delta
x$ is defined in~\fref{fig:ex-xfem-domain}.  The relative distance of a point
singularity to the boundary of the integration domain is shown to affect
integration accuracy in the examples in~\sref{sec:ex-singular}.  Selecting
different values of $\Delta x$ for this problem is designed to highlight this
effect.  Reference integrals for this example are computed using the SBC method
with many cubature points.

\smallskip
To enable singularity cancelling transformations with the SBC method, the
location of the crack tip ($r = 0$ in $F(\vx)$) is selected as $\vx_0 = \vx_c$.
Since $F(\vx)$ is discontinuous, each integration domain must be split across
the discontinuity to ensure smooth, continuous integrands over each triangular
partition.  The location of the discontinuity for this example is denoted by a
thick line in~\fref{fig:ex-xfem-domain}. Selecting the point $r = 0$ as $\vx_0$
and splitting the domains across the discontinuity results in the triangular
decomposition of the element domains pictured in~\fref{fig:ex-xfem-decomp}.
Since \revone{$\vx_0 = \vx_c \notin \bar{\Omega}_1$}, some cubature
points will lie outside the domain of integration and some cubature weights will
be negative. These points and weights coincide with the horizontally hatched
triangles in~\fref{fig:ex-xfem-decomp}.

\begin{figure}[t]
  \centering
  \includegraphics[scale=1]{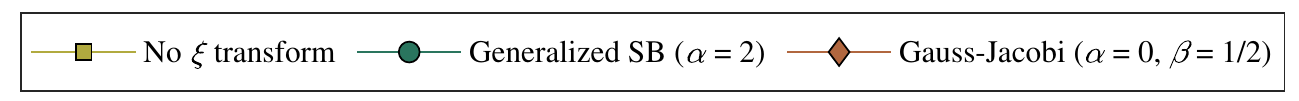}\\[0.1in]
  \begin{subfigure}{0.57in}
    \caption*{}
  \end{subfigure}
  \begin{subfigure}{1.81in}
    \centering
    $\Delta x = 0.001$
    \caption*{}
  \end{subfigure}
  \begin{subfigure}{1.81in}
    \centering
    $\Delta x = 0.01$
    \caption*{}
  \end{subfigure}
  \begin{subfigure}{1.81in}
    \centering
    $\Delta x = 0.1$
    \caption*{}
  \end{subfigure}\\[-0.2in]
  \begin{subfigure}{0.57in}
    \centering
    $\Omega_1$
    \caption*{}
  \end{subfigure}
  \begin{subfigure}{1.81in}
    \centering
    \includegraphics[scale=1]{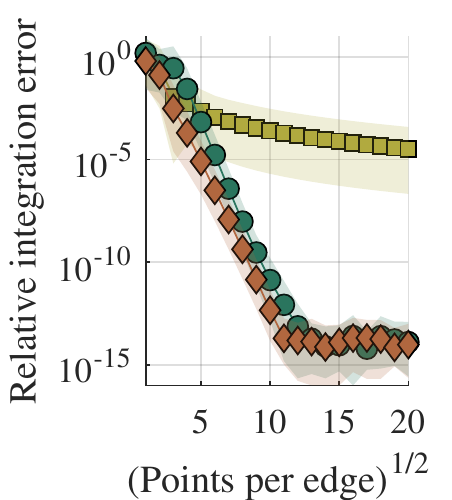}
    \caption{}\label{fig:ex-xfem-conv-11}
  \end{subfigure}
  \begin{subfigure}{1.81in}
    \centering
    \includegraphics[scale=1]{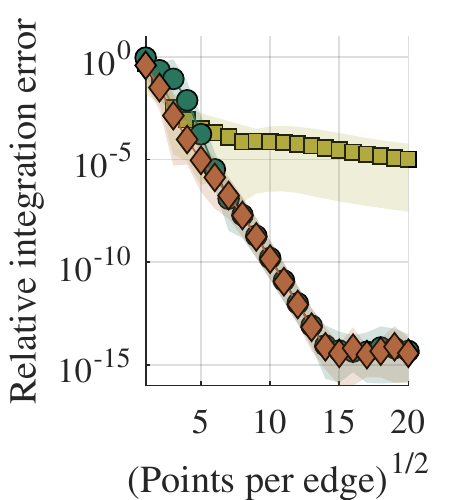}
    \caption{}\label{fig:ex-xfem-conv-21}
  \end{subfigure}
  \begin{subfigure}{1.81in}
    \centering
    \includegraphics[scale=1]{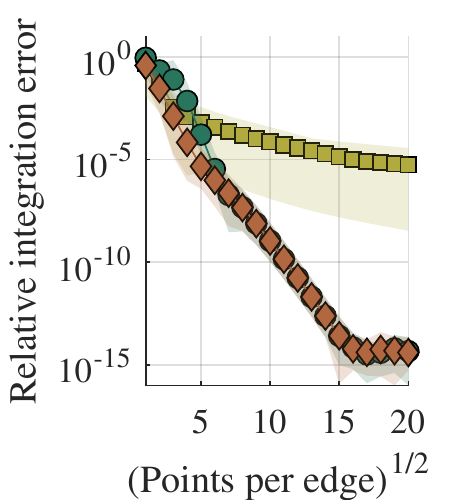}
    \caption{}\label{fig:ex-xfem-conv-31}
  \end{subfigure}
  \begin{subfigure}{0.57in}
    \centering
    $\Omega_2$
    \caption*{}
  \end{subfigure}
  \begin{subfigure}{1.81in}
    \centering
    \includegraphics[scale=1]{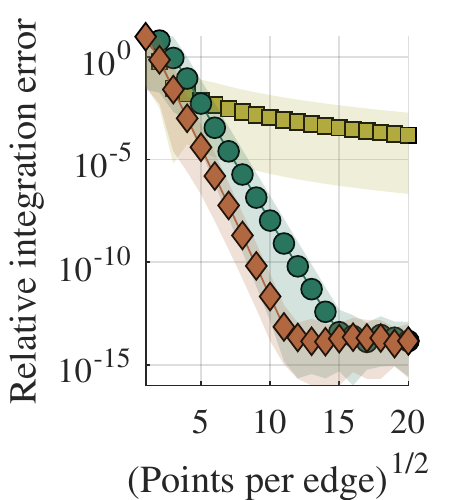}
    \caption{}\label{fig:ex-xfem-conv-12}
  \end{subfigure}
  \begin{subfigure}{1.81in}
    \centering
    \includegraphics[scale=1]{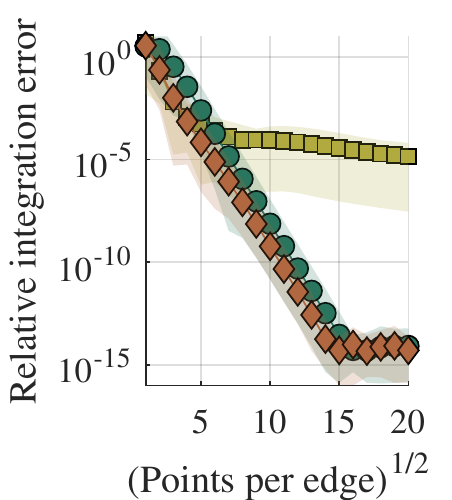}
    \caption{}\label{fig:ex-xfem-conv-22}
  \end{subfigure}
  \begin{subfigure}{1.81in}
    \centering
    \includegraphics[scale=1]{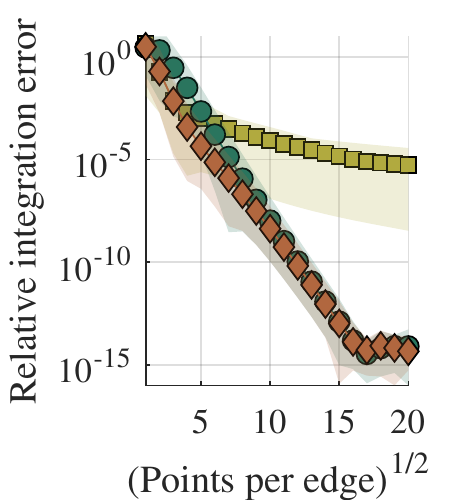}
    \caption{}\label{fig:ex-xfem-conv-32}
  \end{subfigure}
  \caption{Cubature points per edge versus integration error for a \revone{weakly} 
  singular extended finite element enriched basis function that is integrated 
    over $\Omega_1$ with (a) $\Delta x = 0.001$, (b)
    $\Delta x = 0.01$, and (c) $\Delta x = 0.1$ and over $\Omega_2$ with (d)
    $\Delta x = 0.001$, (e) $\Delta x = 0.01$, and (f) $\Delta x = 0.1$. }
  \label{fig:ex-xfem-conv}
\end{figure}

\smallskip
Integration accuracy versus the number of cubature points is illustrated
in~\fref{fig:ex-xfem-conv}.  In the $t$ direction, the $r^{-1}$ singularity
cancelling transformation is applied to all results.  Based on the example
presented in~\sref{sec:ex-singular-t}, this is the transformation that provides
the most accurate integration for an integrand with a $r^{-1/2}$ singularity. In
the $\xi$ direction, three integration methods are tested: no
transformation in the $\xi$ direction, the generalized SB transformation
with $\alpha = 2$, and Gauss\hyp{}Jacobi quadrature with $\alpha = 0$ and $\beta
= 1/2$.  An equal number of cubature points are used in the $t$ and
$\xi$ directions on each triangular partition.  In the plots presented
in~\fref{fig:ex-xfem-conv}, the average error over the 16 integrals is presented
as a solid line.  Shaded regions give the range of error over the 16 integrals.
Compared to the results in~\sref{sec:ex-singular-xi}, the generalized SB
transformation and the Gauss\hyp{}Jacobi rule are closely aligned.  The only
case where the Gauss\hyp{}Jacobi rule provided noticeably better integration is
over $\Omega_2$ with $\Delta x = 0.001$.  Regardless, either method is able to
provide approximately 14 digits of integration accuracy with a $16 \times 16$
rule over each triangular partition.  Without any $\xi$\hyp{}transformation,
less than five digits of integration accuracy are obtained on average with a $20
\times 20$ rule over each triangle.  For this problem, the accuracy of
integration did not seem to be affected by the presence of points outside the
element domain and negative cubature weights, since the results for $\Omega_1$
and $\Omega_2$ are comparable.  The ability to easily generate rules with points
outside the domain with the SBC scheme enables economical integration of
singular integrands \revone{when} the singularity is not present within the domain.
Computing the same integral using, for example, very high\hyp{}order
tensor\hyp{}product Gauss cubature provides unsatisfactory
results~\cite{Chin:2017:MCD}.

\subsection{Transfinite mean value interpolation}\label{sec:ex-tmvi}
\begin{figure}
  \centering
  \begin{tikzpicture}
    \node at (0in,0in) {\includegraphics[scale=1]{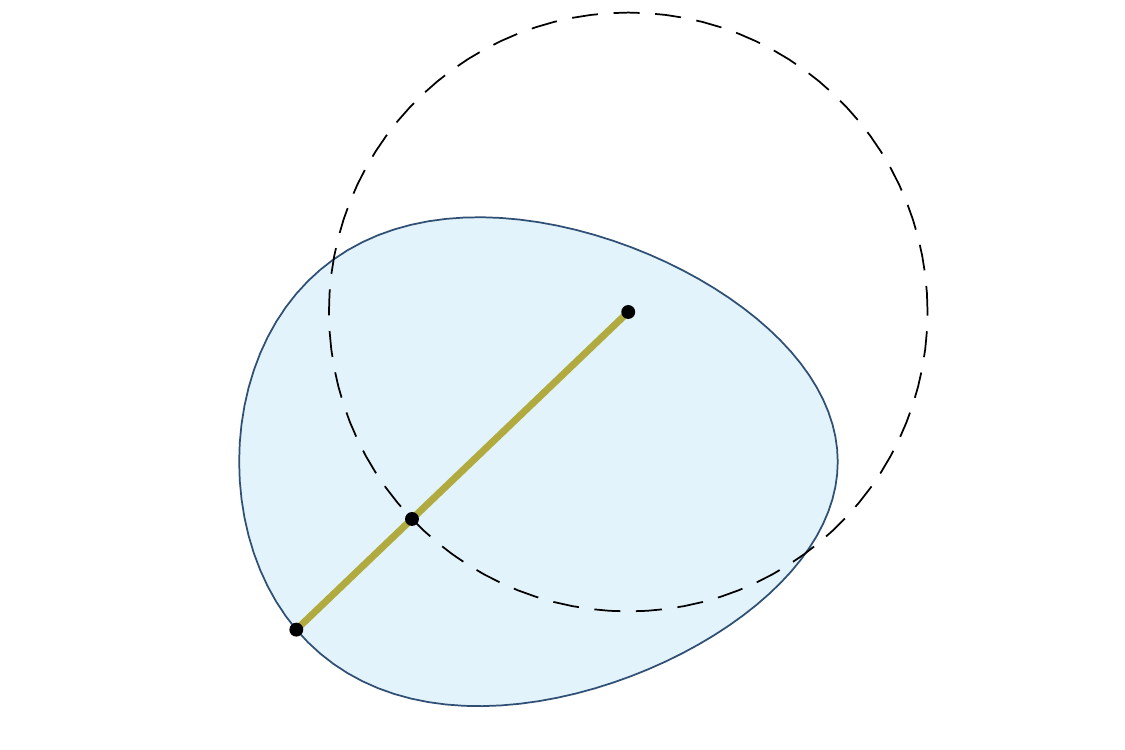}};
    \node at (0.35in,0.20in) {$\vx$};
    \node at (0.1in,-0.4in) {$\Omega$};
    \node at (-1.35in,0in) {$\Gamma$};
    \node at (-1.29in,-1.1in) {$\vm{y}(\vx, \vm{v})$};
    \node at (-0.72in,-0.57in) {$\vm{v}$};
    \node at (1.35in,1in) {$S_v$};
  \end{tikzpicture}
  \caption{Egg\hyp{}shaped convex domain, with illustration of variables that
            appear in~\eqref{eq:tmvi}.}
  \label{fig:ex-tmvi-egg}
\end{figure}

Transfinite barycentric interpolation over domains bounded by curves is the
continuous counterpart of generalized barycentric coordinates over
polygons~\cite{Belyaev:2017:TBC}.  Consider an open, bounded convex domain
$\Omega$ with boundary $\Gamma = \partial \Omega$.  Generalized barycentric
coordinates~\cite{Floater:2013:GBC,Anisimov:2017:BCP}, $\vm{\phi} : \Omega
\rightarrow \Re_+^m$, on an $m$\hyp{}gon are nonnegative and satisfy the
partition of unity and linear precision properties:
\begin{equation*}
  \vx = \sum_{i=1}^m \phi_i(\vx) \vx_i, \quad
  \phi_i(\vx) = \frac{w_i(\vx)}{W(\vx)}, \quad
  W(\vx) = \sum_{j=1}^m w_j(\vx),
\end{equation*}
where the weight function $w_i : \bar{\Omega} \rightarrow \Re_+$ is nonnegative.
The weights $\vm{w} = (w_1, w_2, \dotsc, w_m)$ are also known as homogeneous
coordinates.  Similarly, given a function $g : \Gamma \rightarrow \Re$ that is
prescribed on a curved boundary, the transfinite mean value interpolant $u :
\Omega \rightarrow \Re$ is defined as~\cite{Ju:2005:MVC,Dyken:2009:TMV}
\begin{equation}\label{eq:tmvi}
  u(\vx) = \frac{\int_{S_v} g \bigl( \vm{y} (\vx, \vm{v}) \bigr) K(\vx, \vm{y}) \, dS_v}
    {W(\vx)} , \quad
  W(\vx) = \int_{S_v} K(\vx, \vm{v}) \, dS_v , \quad
  K(\vx, \vm{v}) = \frac{1}{\|\vm{v} - \vx\|} ,
\end{equation}
where $\vx \in \Omega \backslash S_v$ and $\vm{v} \in S_v$, $S_v$ is the unit
circle that is centered at $\vx$, the ray from $\vx$ that passes through
$\vm{v}$ intersects the boundary $\Gamma$ at $\vm{y}$, and $K(\vx, \vm{v})$ is a
singular kernel function~\cite{Ju:2005:MVC,Dyken:2009:TMV}.  See
\fref{fig:ex-tmvi-egg} for an illustration.  
For this study, we let $g(\vx)$ be a function that is defined over
$\bar{\Omega}$ whose mean value interpolant is sought. We have $\lim_{x
\rightarrow y} u(\vx) = g(\vm{y})$, and in addition, if $g$ is an affine
function over $\bar{\Omega}$, then the transfinite interpolant $u = g$.  So
similar to GBCs on polygons, the TMVI can exactly reproduce any linear
polynomial.  The function $\psi(\vx) = 1 / W(\vx)$ behaves like an approximate
distance function to the boundary and its normal derivative on the boundary
$\Gamma$ is $1/2$~\cite{Dyken:2009:TMV}.  Dyken and
Floater~\cite{Dyken:2009:TMV} use the circumferential mean value theorem to
derive the transfinite mean value interpolant in~\eqref{eq:tmvi}, and its
representation is analogous in the discrete case to the Shepard interpolant
where inverse-distance singular nodal weight functions are used for scattered
data interpolation~\cite{Shepard:1968:ATD}.

\smallskip
On projecting the boundary curve $\Gamma$ onto the unit circle, the integrals
in~\eqref{eq:tmvi} can be written for a convex domain in terms of the curve
parameter $t \in [0,1]$ as~\cite{Ju:2005:MVC}
\begin{equation}\label{eq:tmvi-t}
  u(\vx) = \frac{\int_0^1 g \bigl( \vm{c}(t) \bigr) K(\vx, t) \, dt}{W(\vx)}, \quad
  W(\vx) = \int_0^1 K(\vx, t) \, dt, \quad
  K(\vx, t) = \frac{\bigl( \vm{c}(t) - \vx \bigr) \cdot \vm{c}^{\prime\perp} (t)}
    {\| \vm{c}(t) - \vx \|^3} ,
\end{equation}
which is convenient for numerical computations.  In~\eqref{eq:tmvi-t},
$\vm{c}^{\prime\perp}(t) := \text{rot} \bigl( \vm{c}^\prime (t) \bigr)$ is
obtained by rotating $\vm{c}^\prime (t)$ through $90\degree$ in the clockwise
direction.  For a nonconvex domain, a ray from $\vx$ can intersect the boundary
$C$ at one or more (odd number) points; see Dyken and
Floater~\cite{Dyken:2009:TMV} for the general form of~\eqref{eq:tmvi-t}.

\smallskip
As noted in~\cite{Dyken:2009:TMV}, the TMVI is derived by
Lee~\cite{Lee:2007:MVR} using flux integrals and Euler's homogeneous function
theorem, with generalizations of mean value coordinates on polygons and
polyhedra appearing in Lee~\cite{Lee:2009:VFM}.  Here we present the
connections.  Consider a convex domain $\Omega$ with boundary $\Gamma$.  Since
the Laplace operator is rotationally invariant, the kernel $K(\vx, \vm{y}) :=
K(\vm{y} - \vx)$ is only dependent on the radial coordinate.  Let $\vm{\rho} =
\vm{y} - \vx$ and $\rho := \| \vm{\rho} \|$, where $\vm{y} \in C$.  On noting
that $\vm{n}_\rho = \frac{\vm{c}^{\prime\perp} (t)}{\| \vm{c}^\prime (t) \|}$ is
the unit outward normal vector to $C$ and $ds_\rho = \| \vm{c}^\prime (t) \| \, dt$
is the arc\hyp{}length differential on $C$, we can rewrite~\eqref{eq:tmvi-t} as
\begin{equation}\label{eq:tmvi-flux}
  u(\vx) = \frac{\int_\Gamma g(\vx + \vm{\rho}) K(\vx, \vm{\rho}) \, ds_\rho}
    {\int_\Gamma K(\vx, \vm{\rho}) \, ds_\rho} , \quad
  K(\vx, \vm{\rho}) = \frac{\vm{\rho}}{\rho^3} \cdot \vm{n}_\rho .
\end{equation}

\smallskip
Now, \eqref{eq:tmvi-flux} is in the form of a flux integral.  So, on noting the
correspondence, we consider $K(\vx, \vm{\rho}) := \vm{F}(\vm{\rho}) \cdot
\vm{n}_\rho$ to be the ansatz for the kernel function, where $\vm{F}(\vm{\rho})
:= h(\rho) \vm{\rho}$ is a vector field with $h(\rho)$ being a
scalar\hyp{}valued function.  If $g$ is an affine function, then the TMVI
interpolant $u = g$.  So, let $g(\vx) = \vm{a} \cdot \vx$, where $\vm{a} \in
\Re^2$ is an arbitrary vector.  Then, substituting $u(\vx) = \vm{a} \cdot \vx$
on the left\hyp{}hand side and $g(\vx + \vm{\rho}) = \vm{a} \cdot (\vx +
\vm{\rho})$ on the right\hyp{}hand side, we obtain
\begin{equation}
  \vm{a} \cdot \vx = \vm{a} \cdot
    \frac{\int_\Gamma (\vx + \vm{\rho}) 
      \left[ \vm{F} (\vm{\rho}) \cdot \vm{n}_\rho \right] ds_\rho}
    {\int_\Gamma \left[ \vm{F}(\vm{\rho}) \cdot \vm{n}_\rho \right] ds_\rho} =
  \vm{a} \cdot \vx + \vm{a} \cdot
    \frac{\int_\Gamma \vm{\rho}
      \left[ \vm{F} (\vm{\rho}) \cdot \vm{n}_\rho \right] ds_\rho}
    {\int_\Gamma \left[ \vm{F}(\vm{\rho}) \cdot \vm{n}_\rho \right] ds_\rho} ,
\end{equation}
and therefore, since $\vm{a}$ is arbitrary,
\begin{equation}
  \vm{I} := \frac{\int_\Gamma \vm{\rho}
    \left[ \vm{F} (\vm{\rho}) \cdot \vm{n}_\rho \right] ds_\rho}
    {\int_\Gamma \left[ \vm{F}(\vm{\rho}) \cdot \vm{n}_\rho \right] ds_\rho} = \vm{0}
\end{equation}
must hold in $\Re^2 \backslash \{ \vm{0} \}$.  For the above to be true, the
numerator must be identically equal to $\vm{0}$ and the denominator must not
vanish.  On substituting $\vm{F}(\vm{\rho}) = h(\rho) \vm{\rho}$ in the above
equation, we obtain
\begin{equation}
  \vm{I} = \frac{\int_\Gamma
    \left[ h(\rho) \vm{\rho} \otimes \vm{\rho} \right] \cdot \vm{n}_\rho ds_\rho}
    {\int_\Gamma \left[ h(\rho) \vm{\rho} \cdot \vm{n}_\rho \right] ds_\rho} ,
\end{equation}
and on invoking the divergence theorem, we have
\begin{equation}
  \vm{I} = \lim_{\epsilon \rightarrow 0^+} \frac{\int_{\Omega \backslash S_\epsilon}
    \left[ h(\rho) (\vm{I} \cdot \vm{\rho}) + h(\rho) (\nabla \cdot \vm{\rho}) \vm{\rho}
      + \vm{\rho} (\nabla h (\rho) \cdot \vm{\rho}) \right] dA}
    {\int_{\Omega \backslash S_\epsilon} \left[ \nabla h(\rho) \cdot \vm{\rho} 
      + h(\rho) (\nabla \cdot \vm{\rho}) \right] dA} ,
\end{equation}
where $S_\epsilon$ is a circle of radius $\epsilon$ that is centered at $\vx$.
Now, let $h(\rho) := \rho^q$ be a $q$\hyp{}homogeneous function.  We have
$\nabla \cdot \vm{\rho} = 2$ and by Euler's homogeneous function theorem $\nabla
h \cdot \vm{\rho} = q h$.  Then, the above equation becomes
\begin{equation}
  \vm{I} = \frac{\int_{\Omega \backslash S_\epsilon} (q+3) h(\rho) \vm{\rho} \, dA}
    {\int_{\Omega \backslash S_\epsilon} (q+2) h(\rho) \, dA} .
\end{equation}
So if $q = -3$, then the numerator vanishes for $\vm{\rho} \in \Re^2 \backslash
\{ \vm{0} \}$ and $\vm{I} = \vm{0}$ a.e., since the integrand in the
denominator, $-h(\rho)$, is strictly negative.  Therefore, by ensuring linear
precision, we find that the vector field $\vm{F} = \frac{\vm{\rho}}{\rho^3}$
matches that used in the transfinite mean value interpolant
in~\eqref{eq:tmvi-t}.  This provides the flux integral representation to
construct transfinite interpolants~\cite{Lee:2007:MVR,Lee:2009:VFM}.

\begin{figure}[t]
  \centering
  \begin{subfigure}{3in}
    \centering
    \includegraphics[scale=1]{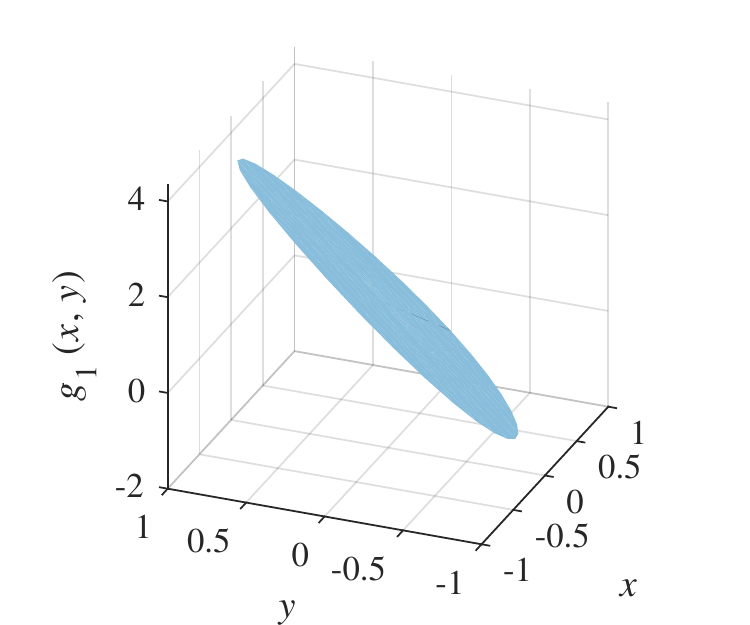}
    \caption{}\label{fig:ex-tmvi-g1}
  \end{subfigure}
  \begin{subfigure}{3in}
    \centering
    \includegraphics[scale=1]{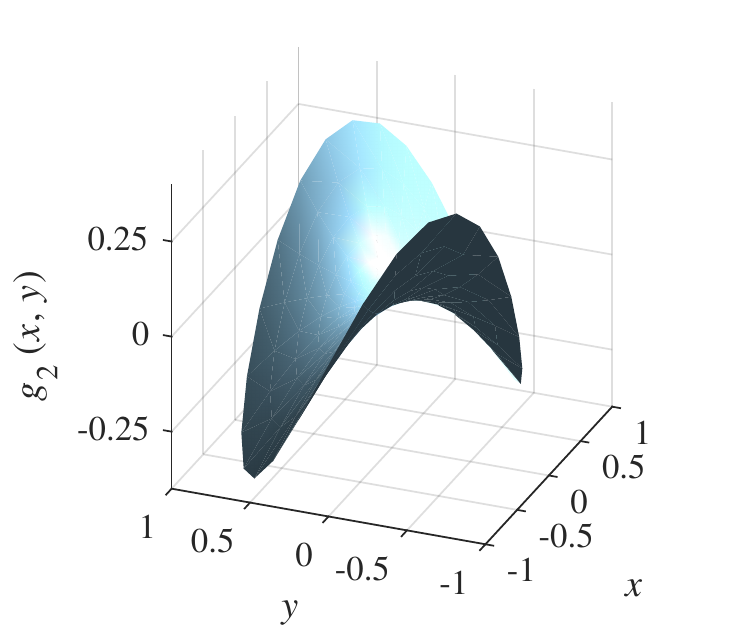}
    \caption{}\label{fig:ex-tmvi-g2}
  \end{subfigure}
  \begin{subfigure}{3in}
    \centering
    \includegraphics[scale=1]{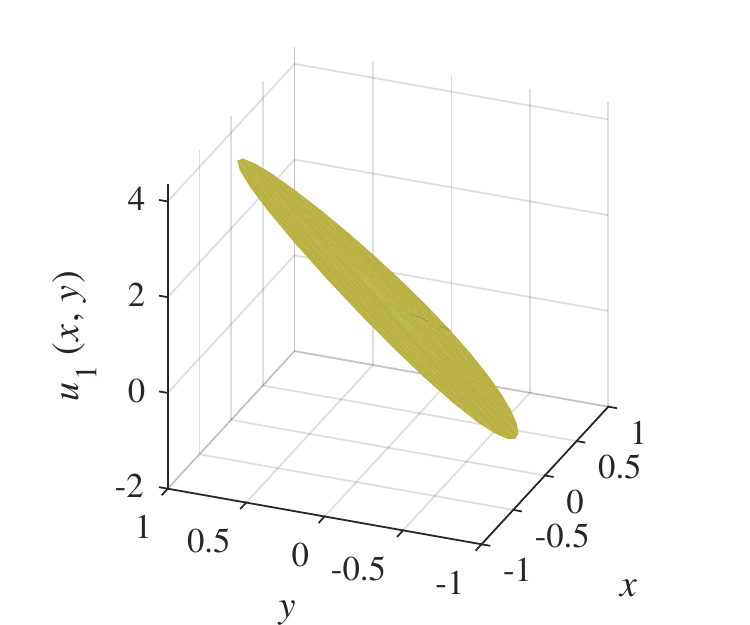}
    \caption{}\label{fig:ex-tmvi-u1}
  \end{subfigure}
  \begin{subfigure}{3in}
    \centering
    \includegraphics[scale=1]{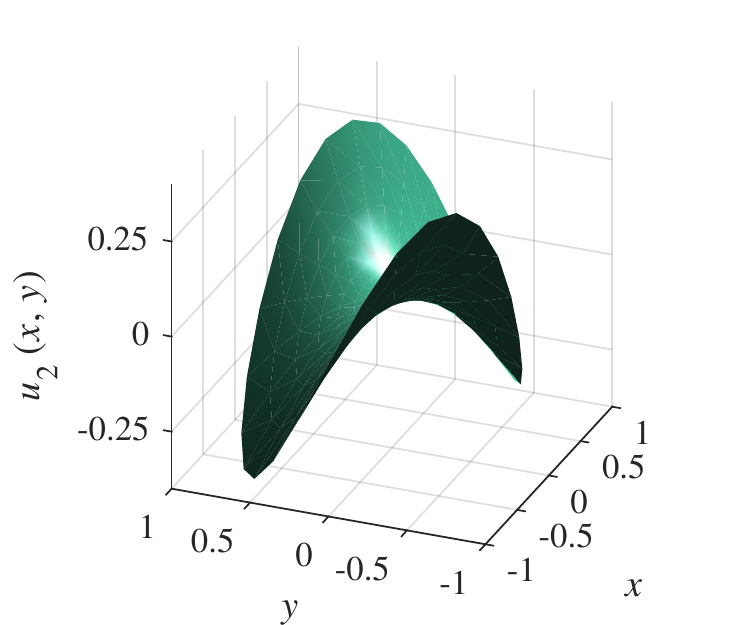}
    \caption{}\label{fig:ex-tmvi-u2}
  \end{subfigure}
  \caption{Interpolating functions over the domain pictured
  in~\fref{fig:ex-tmvi-egg} using the TMVI.  (a) $g_1(\vx)$; (b) $g_2 (\vx)$;
  (c) $u_1(\vx)$, the TMVI of $g_1(\vx)$; and (d) $u_2(\vx)$, the TMVI of
  $g_2(\vx)$.}
  \label{fig:ex-tmvi}
\end{figure}

\smallskip
We now apply the TMVI to two functions over a curved domain: a linear polynomial
and a trigonometric function.  Since linear polynomials are exactly reproduced
by the TMVI, the linear polynomial should be reproduced without any error.  Both
functions are interpolated over an open, bounded convex domain $\Omega$ with
boundary $\Gamma = \partial \Omega$ (see \fref{fig:ex-tmvi-egg}).  The curved
boundary $\Gamma : [0,1] \rightarrow \Re^2$ shown in~\fref{fig:ex-tmvi-egg} has
the parametric representation:
\begin{equation}\label{eq:egg-bdry}
  \vm{c} \bigl(\theta(t)\bigr) = r \cos \theta \, \vm{e}_1 + \frac{a r \sin \theta}
    {b + r \cos \theta} \, \vm{e}_2 ,
  \quad
  \vm{c}^\prime \bigl( \theta(t) \bigr) = -2 \pi r \sin \theta \, \vm{e}_1 + 
    \frac{2 \pi a r (r + b \cos \theta)}{(b + r \cos \theta)^2} \, \vm{e}_2 ,
\end{equation}
where $\theta := \theta(t) = 2 \pi t$ and $(\cdot)^\prime$ denotes
differentiation with respect to $t$.  In these numerical computations, we choose
$a = 4$, $b = 5$, and $r = 1$.  The functions we choose to interpolate over
$\bar{\Omega}$ are:
\begin{equation}
  g_1(\vx) = 1 - 2x + 3y , \qquad g_2(\vx) = \sin x \, \sin y .
\end{equation}
These functions are plotted in Figs.~\ref{fig:ex-tmvi-g1} and
\ref{fig:ex-tmvi-g2}.  The TMVI is computed using~\eqref{eq:tmvi-t} and the
resulting interpolants are illustrated in Figs.~\ref{fig:ex-tmvi-u1} and
\ref{fig:ex-tmvi-u2}.  The relative $L_2$ error over $\Omega$ is defined as:
\begin{equation}
{\cal E} = \frac{\| u - g \|_{2,\Omega}} {\| u \|_{2,\Omega}},
\end{equation}
which is computed for $g_1(\vx)$ and $g_2(\vx)$ using the SBC scheme.  As
expected, interpolation of $g_1(\vx)$ yields machine precision relative $L_2$
error of $2.2 \times 10^{-16}$. Interpolation of the trigonometric function
using TMVI is not exact away from $\Gamma$.  Though the surfaces displayed
in~\fref{fig:ex-tmvi-g2} and \fref{fig:ex-tmvi-u2} appear nearly identical, the
relative $L_2$ error is $1.24 \times 10^{-2}$.

\begin{figure}[t]
  \centering
  \begin{subfigure}{3in}
    \centering
    \includegraphics[scale=1]{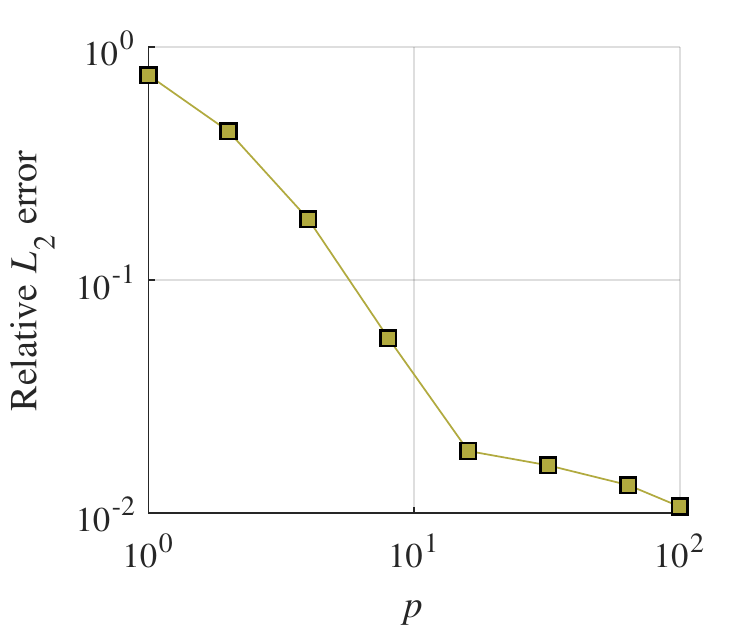}
    \caption{}\label{fig:ex-tmvi-dist-conv}
  \end{subfigure}
  \begin{subfigure}{3in}
    \centering
    \includegraphics[scale=1]{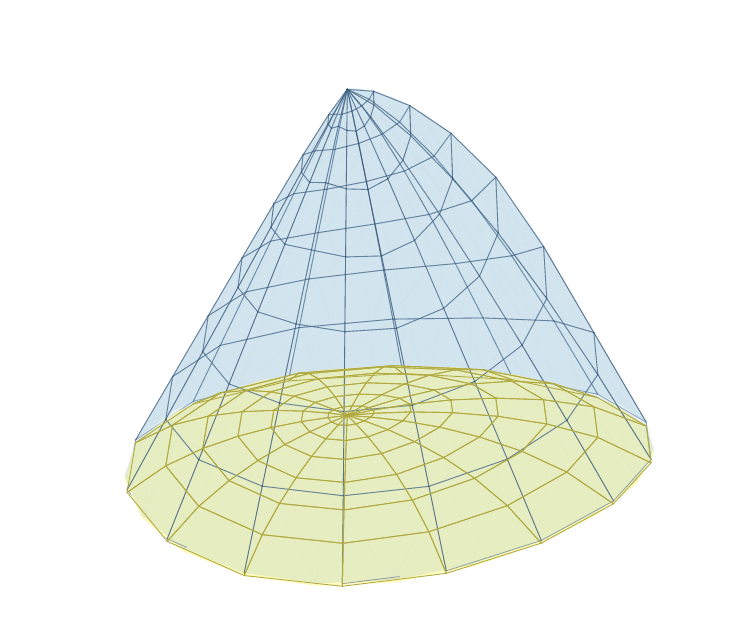}
    \caption{}\label{fig:ex-tmvi-dist-p1}
  \end{subfigure}
  \begin{subfigure}{3in}
    \centering
    \includegraphics[scale=1]{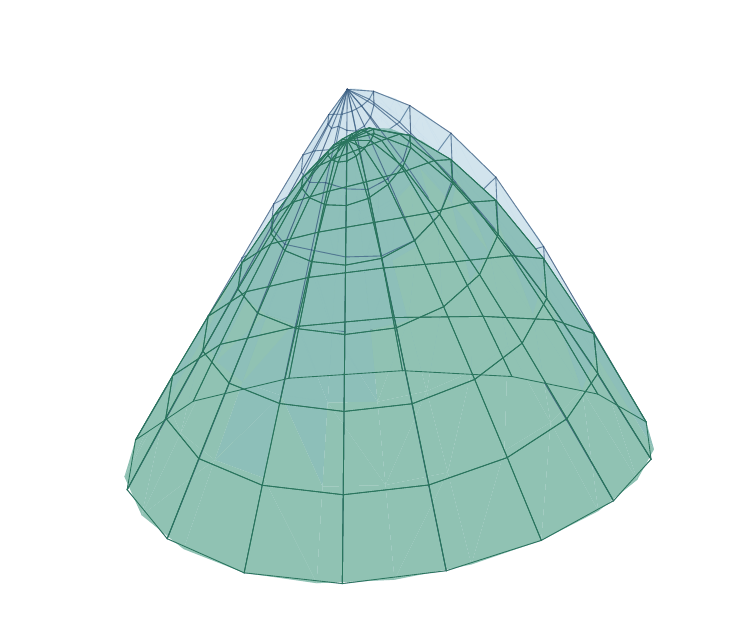}
    \caption{}\label{fig:ex-tmvi-dist-p10}
  \end{subfigure}
  \begin{subfigure}{3in}
    \centering
    \includegraphics[scale=1]{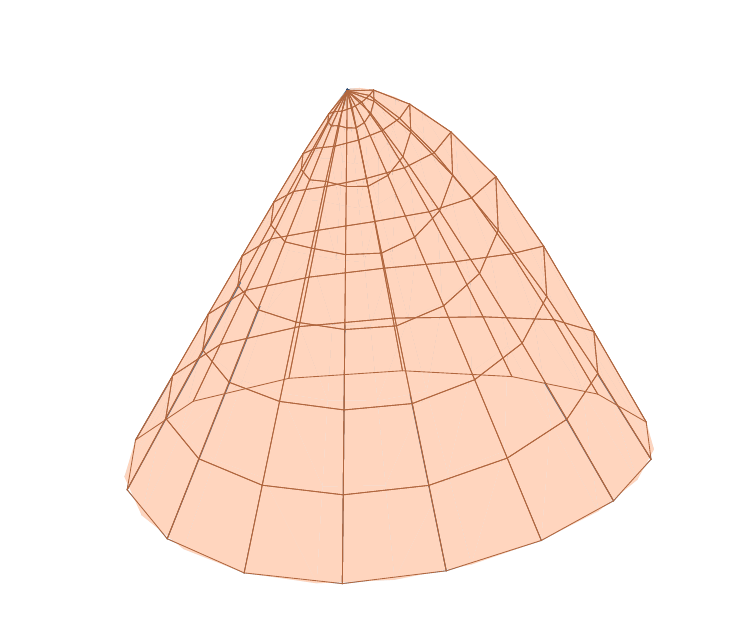}
    \caption{}\label{fig:ex-tmvi-dist-p100}
  \end{subfigure}
  \caption{Approximating the distance function over the domain pictured
  in~\fref{fig:ex-tmvi-egg} using $L_p$\hyp{}distance fields.  (a) $L_2$ error
  versus $p$, $L_p$\hyp{}distance field with (b) $p = 1$, (c) $p = 10$, and (d)
  $p = 100$.  In (b)--(d), the distance function is plotted as a
  semi\hyp{}transparent overlay.}
  \label{fig:ex-tmvi-dist}
\end{figure}

\smallskip
As a second, related example, we use a modified inverse of the weight function
in the TMVI to provide a smooth approximation of the distance function.  Belyaev
et al.~\cite{Belyaev:2013:SLD} introduced $L_p$\hyp{}distance fields ($p \geq
1$):
\begin{equation}\label{eq:tmvi-p}
  \psi (\vx) = \left( \frac{1}{W_p (\vx)} \right)^{1/p}, \qquad
  W_p(\vx) = \int_0^1 \frac{ \bigl( \vm{c}(t) - \vx \bigr) \cdot 
    \vm{c}^{\prime\perp} (t) }{ \| \vm{c}(t) - \vx \|^{2+p} } ,
\end{equation}
which approximates the exact distance function.  In~\eqref{eq:tmvi-p}, as $p
\rightarrow \infty$, $\psi (\vx)$ approaches the exact distance function.
Selecting $p = 1$ recovers the weight function in the TMVI.  
Consider the egg\hyp{}shaped domain $\Omega$ (\fref{fig:ex-tmvi-egg}) with
boundary $\Gamma$ whose parametric representation is given
in~\eqref{eq:egg-bdry}. In~\fref{fig:ex-tmvi-dist-conv}, we plot the relative
$L_2$ error for the distance function versus $p$ ($p \in [1,100]$).  With $p =
1$, the relative $L_2$ error is 76 percent and with $p = 100$, relative $L_2$
error reduces to $1.1$ percent.  The approximate distance functions for $p = 1$,
$10$, and $100$ are plotted in Figs.~\ref{fig:ex-tmvi-dist-p1},
\ref{fig:ex-tmvi-dist-p10}, and~\ref{fig:ex-tmvi-dist-p100}, respectively.  The
exact distance function,
\begin{equation}
  d(\vx) = \min_{t\in[0,1]} \| \vx - \vm{c}(t) \| \qquad \forall \vx \in \Omega ,
\end{equation}
is calculated by minimizing the square of the distance function using Newton's
method.

\section{Conclusions}\label{sec:conclusion} 
In this paper, we proposed the scaled boundary cubature (SBC) method for
numerical integration over polygons and domains that are bounded by parametric
curves. This method of integration applies the SB parametrization introduced in
the context of the SB-FEM~\cite{Song:1997:TSB,Wolf:2001:TSB} to transform
integration from a region bounded by a curve and two line segments to the unit
square.  As explained in~\sref{sec:hni}, the method is very closely related to
the homogeneous numerical integration scheme
(HNI)~\cite{Chin:2015:NIH,Chin:2019:MCI}, but is not solely limited to the
integration of homogeneous functions.  A broad suite of numerical examples were
considered in~\sref{sec:examples} that revealed the promise of using the SBC
scheme in the polygonal finite element method, extended finite element method,
and transfinite interpolation of functions over curved domains using the
transfinite mean value interpolant (TMVI). In these examples, integration was
performed over a number of convex and nonconvex shapes, demonstrating the
flexibility of the method. Furthermore, in~\sref{sec:ex-tmvi} the close
connection between the TMVI and homogeneous functions was noted, paralleling the
link between the SBC and HNI schemes. The simplicity of the SBC scheme enabled
fast cubature rule generation when compared to other methods of integration,
such as Gauss\hyp{}Green cubature~\cite{Sommariva:2007:PGC,Sommariva:2009:GGC}
and triangulation.  In the examples presented in~\sref{sec:examples}, accuracy
per cubature point for the SBC scheme was comparable to triangulation and the
Gauss\hyp{}Green scheme.  For star\hyp{}convex domains, it is possible to place
all cubature points inside the domain of integration and to have all positive
cubature weights with the SBC scheme.

\smallskip
The SBC method was also suitable for integrating functions with point
singularities. On subdividing the domain of integration into curved triangles
that share a point, the SBC scheme simplifies the formulation of integral
transformations that smooth singularities at that point.
In~\sref{sec:singular-xi}, we demonstrated the SBC scheme is equivalent to the
Duffy transformation~\cite{Duffy:1982:QOP} when the domain of integration is the
standard triangle.  Therefore, the SBC scheme provides efficient integration of
functions containing a $r^{-1}$ radial singularity without modification. We also
introduced a modified mapping in~\sref{sec:singular-xi} that we named the
generalized SB transformation due to its similarity to the generalized Duffy
transformation~\cite{Mousavi:2010:GDT} when applied to the standard triangle.
This led to efficient integration of functions with fractional radial
singularities.  Additionally, we pointed out that Gauss\hyp{}Jacobi quadrature
also provides efficient integration of functions with fractional power point
singularities.  To smooth near\hyp{}singularities along the edges of polygons,
we described several $t$\hyp{}direction integral transformations
in~\sref{sec:singular-t}.  When the transformations in~\sref{sec:singular-t} are
combined with the transformations in~\sref{sec:singular-xi}, accurate
integration of weakly singular functions with minimal cubature points is
realized.  These capabilities were explored in examples that appeared
in~\sref{sec:ex-singular}. As these examples revealed, the SBC scheme with
singularity cancelling transformations can be readily applied to the X-FEM
and the BEM.

\smallskip
The fundamental need to integrate functions over polytopes and curved solids
appears in numerical methods and also within many higher\hyp{}order embedded
domain methods, where \textit{exact} representation of geometry is desirable. In
this paper, we have proposed a new method for numerical integration of smooth
functions and weakly singular functions over elements with curved boundaries. We
showed that the method is accuracy and efficient, while underscoring the
simplicity in its implementation. The performance of SBC reveals that it can
find use in polygonal finite element methods, boundary element methods, enriched
partition-of-unity methods, embedded domain methods, and isogeometric analysis.
Since the SBC scheme enables tensor\hyp{}product Gauss integration rules to be
utilized over complex element shapes, there is greater flexibility in the choice
of approximation spaces that one can adopt in the underlying computational
method.  The main ideas and contributions in this work readily extend to
three\hyp{}dimensional regions that are bounded by affine and curved parametric
surfaces, which \revone{is the subject of ongoing research}.

\appendix
\section{Connection to the Poincar\'{e} Lemma}\label{sec:poincare}
In this appendix, we use the Poincar\'{e} Lemma to provide an alternate
derivation for the SBC scheme.  The Poincar\'{e} Lemma for vector fields
is~\cite{Terrell:2009:FTC}
\begin{equation}\label{eq:poincare}
  h(\vx) = \nabla \cdot \left[ \int_0^1 h(\xi \vx) \, \xi \vx
    \, d\xi \right] .
\end{equation}
We introduce a curved triangle $\mathcal{T}$ that is bounded by two affine
lines, $\mathcal{L}_1$ and $\mathcal{L}_2$, and a curved boundary $\mathcal{C} :
[0,1] \rightarrow \Re^2$ that is parametrized by $\vm{c}(t)$.  Let the vertex at
the intersection of $\mathcal{L}_1$ and $\mathcal{L}_2$ lie at the origin.  We
integrate $h(\vx)$ over $\mathcal{T}$ and apply the divergence theorem:
\begin{equation}
  \int_\mathcal{T} h(\vx) \, d\vx = \int_\mathcal{T} \nabla \cdot \left[
    \int_0^1 h(\xi \vx) \, \xi \vx \, d \xi \right] d\vx = 
    \int_{\partial \mathcal{T}} \left[ \int_0^1 h(\xi \vx) \,
      \xi \vx \, d\xi \right] \cdot \vm{n} \, ds.
\end{equation}
On noting that coordinates on the segments $\mathcal{L}_1 \backslash \{ \vm{0}
\}$ and $\mathcal{L}_2 \backslash \{ \vm{0} \}$ are perpendicular to the outward
normal (see~\eqref{eq:hni-curvedtri2}), we can write
\begin{equation}\label{eq:poincare-sbc}
  \int_\mathcal{T} h(\vx) \, d\vx = \int_\mathcal{C} \left[
    \int_0^1 h(\xi \vx) \, \xi \, d \xi \right] 
      (\vx \cdot \vm{n}) \, ds =
    \int_0^1 \int_0^1 h(\xi \vm{c}(t)) \, \xi \, d \xi \,
      \bigl( \vm{c}(t) \cdot \vm{n}(t) \bigr) \, \| \vm{c}^\prime (t) \| \, dt .
\end{equation}
Upon substituting the expression for the unit normal given
in~\eqref{eq:curve-norm} into~\eqref{eq:poincare-sbc} and noting $h(\xi
\vm{c}(t)) = (h \circ \vm{\psi})(\xi \rho(t), t)$ (see~\eqref{eq:polar-xform}),
\eqref{eq:hni-xi} is recovered.

If~\eqref{eq:poincare} is modified to incorporate a constant shift $\vx_0$:
\begin{equation}
  h(\vx) = \nabla \cdot \left[ \int_0^1 h(\vm{\varphi}) \, \xi (\vx - \vx_0)
    \, d\xi \right] ,
\end{equation}
where $\vm{\varphi} := \vm{\varphi}(\xi, t) = \vx_0 + \xi(\vx - \vx_0)$ is the
SB parametrization introduced in~\eqref{eq:sb-xform},
then~\eqref{eq:2d-curved-subdomain-integral} is the result after integrating
over $\mathcal{T}$.  This reveals that the SBC method can be derived using the
Poincar\'{e} Lemma.

\section*{Declaration of competing interest}
The authors declare that they have no known competing financial interests or
personal relationships that could have appeared to influence the work reported
in this paper.

\section*{Acknowledgements}
This work was performed under the auspices of the U.S. Department of Energy by
Lawrence Livermore National Laboratory under Contract DE-AC52-07NA27344.

\section*{Disclaimer}
{\small
This document was prepared as an account of work sponsored by an agency of the
United States government. Neither the United States government nor Lawrence
Livermore National Security, LLC, nor any of their employees makes any warranty,
expressed or implied, or assumes any legal liability or responsibility for the
accuracy, completeness, or usefulness of any information, apparatus, product, or
process disclosed, or represents that its use would not infringe privately owned
rights. Reference herein to any specific commercial product, process, or service
by trade name, trademark, manufacturer, or otherwise does not necessarily
constitute or imply its endorsement, recommendation, or favoring by the United
States government or Lawrence Livermore National Security, LLC. The views and
opinions of authors expressed herein do not necessarily state or reflect those
of the United States government or Lawrence Livermore National Security, LLC,
and shall not be used for advertising or product endorsement purposes.
}

\bibliography{journals,references}

\begin{thebibliography}{10}
\expandafter\ifx\csname url\endcsname\relax
  \def\url#1{\texttt{#1}}\fi
\expandafter\ifx\csname urlprefix\endcsname\relax\def\urlprefix{URL }\fi
\expandafter\ifx\csname href\endcsname\relax
  \def\href#1#2{#2} \def\path#1{#1}\fi

\bibitem{Lasserre:1998:ICP}
J.~B. Lasserre, Integration on a convex polytope, Proceedings of the American
  Mathematical Society 126~(8) (1998) 2433--2441.

\bibitem{Chin:2015:NIH}
E.~B. Chin, J.~B. Lasserre, N.~Sukumar, Numerical integration of homogeneous
  functions on convex and nonconvex polygons and polyhedra, Computational
  Mechanics 56~(6) (2015) 967--981.

\bibitem{Chin:2019:MCI}
E.~B. Chin, N.~Sukumar, Modeling curved interfaces without element-partitioning
  in the extended finite element method, International Journal for Numerical
  Methods in Engineering 120~(5) (2019) 607--649.

\bibitem{Song:1997:TSB}
C.~Song, J.~P. Wolf, The scaled boundary finite-element method---alias
  consistent infinitesimal finite-element cell method---for elastodynamics,
  Computer Methods in Applied Mechanics and Engineering 147~(3--4) (1997)
  329--355.

\bibitem{Wolf:2001:TSB}
J.~P. Wolf, C.~Song, The scaled boundary finite-element method---a fundamental
  solution-less boundary-element method, Computer Methods in Applied Mechanics
  and Engineering 190 (2001) 5551--5568.

\bibitem{Moes:1999:AFE}
N.~Mo\"{e}s, J.~Dolbow, T.~Belytschko, A finite element method for crack growth
  without remeshing, International Journal for Numerical Methods in Engineering
  46~(1) (1999) 131--150.

\bibitem{Sevilla:2008:NEF}
R.~Sevilla, S.~Fern\'{a}ndez-M\'{e}ndez, A.~Huerta, {NURBS}-enhanced finite
  element method ({NEFEM}), International Journal for Numerical Methods in
  Engineering 76 (2008) 56--83.

\bibitem{Fries:2017:HOM}
T.~P. Fries, S.~Omerovi\'{c}, D.~Sch\"{o}llhammer, J.~Steidl, Higher-order
  meshing of implicit geometries---{P}art {I}: {I}ntegration and interpolation
  in cut elements, Computer Methods in Applied Mechanics and Engineering 313
  (2017) 759--784.

\bibitem{Artioli:2020:ACP}
E.~Artioli, A.~Sommariva, M.~Vianello, Algebraic cubature on polygonal elements
  with a circular edge, Computers \& Mathematics with Applications 79~(7)
  (2020) 2057--2066.

\bibitem{Mousavi:2010:GGQ}
S.~E. Mousavi, H.~Xiao, N.~Sukumar, Generalized {G}aussian quadrature rules on
  arbitrary polygons, International Journal for Numerical Methods in
  Engineering 82~(1) (2010) 99--113.

\bibitem{Saye:2015:HOQ}
R.~I. Saye, High-order quadrature methods for implicitly defined surfaces and
  volumes in hyperrectangles, SIAM Journal on Scientific Computing 37~(2)
  (2015) A993--A1019.

\bibitem{Scholz:2019:NIT}
F.~Scholz, B.~J\"{u}ttler, Numerical integration on trimmed three-dimensional
  domains with implicitly defined trimming surfaces, Computer Methods in
  Applied Mechanics and Engineering 357 (2019) 112577.

\bibitem{Scholz:2020:HOQ}
F.~Scholz, B.~J\"{u}ttler, High-order quadrature on planar domains based on
  transport theorems for implicitly defined moving curves, Tech. Rep.~89,
  Nationales Forschungsnetzwerk, Geometry + Simulation (May 2020).

\bibitem{Sommariva:2007:PGC}
A.~Sommariva, M.~Vianello, Product {G}auss cubature over polygons based on
  {G}reen's integration formula, BIT Numerical Mathematics 47 (2007) 441--453.

\bibitem{Sommariva:2009:GGC}
A.~Sommariva, M.~Vianello, {Gauss-Green} cubature and moment computation over
  arbitrary geometries, Journal of Computational and Applied Mathematics 231
  (2009) 886--896.

\bibitem{Gunderman:2021:SMF}
D.~Gunderman, K.~Weiss, J.~A. Evans, Spectral mesh-free quadrature for planar
  regions bounded by rational parametric curves, Computer-Aided Design 130
  (2021) 102944.

\bibitem{Lasserre:1999:IHF}
J.~B. Lasserre, Integration and homogeneous functions, Proceedings of the
  American Mathematical Society 127~(3) (1999) 813--818.

\bibitem{Chin:2020:AEM}
E.~B. Chin, N.~Sukumar, An efficient method to integrate polynomials over
  polytopes and curved solids, Computer Aided Geometric Design 82 (2020)
  101914.

\bibitem{Chin:2017:MCD}
E.~B. Chin, J.~B. Lasserre, N.~Sukumar, Modeling crack discontinuities without
  element-partitioning in the extended finite element method, International
  Journal for Numerical Methods in Engineering 110~(11) (2017) 1021--1048.

\bibitem{Chen:2016:ANB}
L.~Chen, B.~Simeon, S.~Klinkel, A {NURBS} based {G}alerkin approach for the
  analysis of solids in boundary representation, Computer Methods in Applied
  Mechanics and Engineering 305 (2016) 777--805.

\bibitem{Klinkel:2019:AFE}
S.~Klinkel, R.~Reichel, A finite element formulation in boundary representation
  for the analysis of nonlinear problems in solid mechanics, Computer Methods
  in Applied Mechanics and Engineering 347 (2019) 295--315.

\bibitem{Sukumar:2015:EFE}
N.~Sukumar, J.~E. Dolbow, N.~Mo\"{e}s, Extended finite element method in
  computational fracture mechanics: a retrospective examination, International
  Journal of Fracture 196 (2015) 189--206.

\bibitem{Hormann:2017:GBC}
K.~Hormann, N.~Sukumar (Eds.), Generalized Barycentric Coordinates in Computer
  Graphics and Computational Mechanics, CRC Press, New York, NY, 2017.

\bibitem{Hughes:2005:IAC}
T.~J.~R. Hughes, J.~A. Cottrell, Y.~Bazilevs, Isogeometric analysis: {CAD},
  finite elements, {NURBS}, exact geometry and mesh refinement, Computer
  Methods in Applied Mechanics and Engineering 194~(39--41) (2005) 4135--4195.

\bibitem{Marussig:2018:ICI}
B.~Marussig, T.~J.~R. Hughes, A review of trimming in isogeometric analysis:
  challenges, data exchange, and simulation aspects, Archives of Computational
  Methods in Engineering 25 (2018) 1059--1127.

\bibitem{Bishop:2003:RSA}
J.~E. Bishop, Rapid stress analysis of geometrically complex domains using
  implicit meshing, Computational Mechanics 30~(5) (2003) 460--478.

\bibitem{Burman:2014:CDG}
E.~Burman, S.~Claus, P.~Hansbo, M.~G. Larson, A.~Massing, {CutFEM}:
  {D}iscretizing geometry and partial differential equations, International
  Journal for Numerical Methods in Engineering 104~(7) (2014) 472--501.

\bibitem{Schillinger:2015:TFC}
D.~Schillinger, M.~Ruess, The finite cell method: {A} review in the context of
  higher-order structural analysis of {CAD} and image-based geometric models,
  Archives of Computational Methods in Engineering 22 (2015) 391--455.

\bibitem{BeiraodaVeiga:2013:BPV}
L.~{Beir\~{a}o da Veiga}, F.~Brezzi, A.~Cangiani, G.~Manzini, L.~D. Marini,
  A.~Russo, Basic principles of virtual element methods, Mathematical Models
  and Methods in Applied Sciences 23~(1) (2013) 199--214.

\bibitem{BeiraodaVeiga:2019:TVE}
L.~{Beir\~{a}o da Veiga}, A.~Russo, G.~Vacca, The virtual element method with
  curved edges, ESAIM: Mathematical Modelling and Numerical Analysis 53 (2019)
  375--404.

\bibitem{Benvenuti:2019:EVE}
E.~Benvenuti, A.~Chiozzi, G.~Manzini, N.~Sukumar, Extended virtual element
  method for the {L}aplace problem with singularities and discontinuities,
  Computer Methods in Applied Mechanics and Engineering 256 (2019) 571--597.

\bibitem{Lew:2008:ADG}
A.~J. Lew, G.~C. Buscaglia, A {discontinuous-Galerkin-based} immersed boundary
  method, International Journal for Numerical Methods in Engineering 76 (2008)
  427--454.

\bibitem{Cangiani:2014:HVD}
A.~Cangiani, E.~H. Georgoulis, P.~Houston, $hp$-version discontinuous
  {G}alerkin methods on polygonal and polyhedral meshes, Mathematical Models
  and Methods in Applied Sciences 24~(10) (2014) 2009--2041.

\bibitem{Puso:2004:AMS}
M.~A. Puso, T.~A. Laursen, A mortar segment-to-segment contact method for large
  deformation solid mechanics, Computer Methods in Applied Mechanics and
  Engineering 193~(6--8) (2004) 601--629.

\bibitem{Hesch:2011:TTD}
C.~Hesch, P.~Betsch, Transient three-dimensional contact problems: mortar
  method. {M}ixed methods and conserving integration, Computational Mechanics
  48 (2011) 461--475.

\bibitem{Guendelman:2003:NRB}
E.~Guendelman, R.~Bridson, R.~Fedkiw, Nonconvex rigid bodies with stacking, ACM
  Transactions on Graphics 22~(3) (2003) 871--878.

\bibitem{Krishnamurthy:2011:AGA}
A.~Krishnamurthy, S.~McMains, Accurate {GPU}-accelerated surface integrals for
  moment computation, Computer-Aided Design 43 (2011) 1284--1295.

\bibitem{Mousavi:2010:GDT}
S.~E. Mousavi, N.~Sukumar, Generalized {D}uffy transformation for integrating
  vertex singularities, Computational Mechanics 45 (2010) 127--140.

\bibitem{Franke:1979:ACC}
R.~Franke, A critical comparison of some methods for interpolation of scattered
  data, Tech. Rep. NPS-53-79-003, Naval Postgraduate School (1979).

\bibitem{Ju:2005:MVC}
T.~Ju, S.~Schaefer, J.~Warren, Mean value coordinates for closed triangular
  meshes, ACM Transactions on Graphics 24~(3) (2005) 561--566.

\bibitem{Dyken:2009:TMV}
C.~Dyken, M.~S. Floater, Transfinite mean value interpolation, Computer Aided
  Geometric Design 26~(1) (2009) 117--134.

\bibitem{Terrell:2009:FTC}
R.~E. Terrell,
  \href{http://citeseerx.ist.psu.edu/viewdoc/download?doi=10.1.1.159.3286&rep=rep1&type=pdf}{The
  fundamental theorem of calculus and the {Poincar{\'e}} lemma}, {C}ornell
  University, Ithaca, NY 14853 (2009).
\newline\urlprefix\url{http://citeseerx.ist.psu.edu/viewdoc/download?doi=10.1.1.159.3286&rep=rep1&type=pdf}

\bibitem{Arioli:2019:SBP}
C.~Arioli, A.~Shamanskiy, S.~Klinkel, B.~Simeon, Scaled boundary
  parametrizations in isogeometric analysis, Computer Methods in Applied
  Mechanics and Engineering 349 (2019) 576--594.

\bibitem{Duffy:1982:QOP}
M.~G. Duffy, Quadrature over a pyramid or cube of integrands with a singularity
  at a vertex, SIAM Journal on Numerical Analysis 19~(6) (1982) 1260--1262.

\bibitem{Chernov:2012:ECG}
A.~Chernov, C.~Schwab, Exponential convergence of {G}auss-{J}acobi quadratures
  for singular integrals over simplices in arbitrary dimension, SIAM Journal on
  Numerical Analysis 50~(3) (2012) 1433--1455.

\bibitem{Ma:2002:DTN}
H.~Ma, N.~Kamiya, Distance transformation for the numerical evaluation of near
  singular boundary integrals with various kernels in boundary element method,
  Engineering Analysis with Boundary Elements 26 (2002) 329--339.

\bibitem{Lv:2019:ASD}
J.-H. Lv, Y.-Y. Jiao, X.-T. Feng, P.~Wriggers, X.-Y. Zhuang, T.~Rabczuk, A
  series of {D}uffy-distance transformation for integrating 2{D} and 3{D}
  vertex singularities, International Journal for Numerical Methods in
  Engineering 118~(1) (2019) 38--60.

\bibitem{Chin:2019:NIH}
E.~B. Chin, Numerical integration of homogeneous functions with applications in
  the extended finite element method, Ph.D. thesis, University of California,
  Davis (2019).

\bibitem{Dunavant:1985:HDE}
D.~A. Dunavant, High degree efficient symmetrical {G}aussian quadrature rules
  for the triangle, International Journal for Numerical Methods in Engineering
  21 (1985) 1129--1148.

\bibitem{Wachspress:1975:ARF}
E.~L. Wachspress, A Rational Finite Element Basis, Academic Press, New York,
  1975.

\bibitem{Floater:2014:GBW}
M.~Floater, A.~Gillette, N.~Sukumar, Gradient bounds for {W}achspress
  coordinates on polytopes, SIAM Journal on Numerical Analysis 52~(1) (2014)
  515--532.

\bibitem{Xiao:2009:ANA}
H.~Xiao, Z.~Gimbutas, A numerical algorithm for the construction of efficient
  quadrature rules in two and higher dimensions, Computers \& Mathematics with
  Applications 59 (2009) 663--676.

\bibitem{Nagarajan:1993:AMM}
A.~Nagarajan, S.~Mukherjee, A mapping method for numerical evaluation of
  two-dimensional integrals with $1/r$ singularity, Computational Mechanics 12
  (1993) 19--26.

\bibitem{Belyaev:2017:TBC}
A.~G. Belyaev, P.-A. Fayolle, Transfinite barycentric coordinates, in: Hormann
  and Sukumar  \cite{Hormann:2017:GBC}, pp. 43--62.

\bibitem{Floater:2013:GBC}
M.~S. Floater, Generalized barycentric coordinates and applications, Acta
  Numerica 24 (2015) 161--214.

\bibitem{Anisimov:2017:BCP}
D.~Anisimov, Barycentric coordinates and their properties, in: Hormann and
  Sukumar  \cite{Hormann:2017:GBC}, pp. 3--22.

\bibitem{Shepard:1968:ATD}
D.~Shepard, A two-dimensional interpolation function for irregularly-spaced
  data, in: Proceedings of the 23rd ACM national conference, Association for
  Computing Machinery, New York, New York, 1968, pp. 517--524.

\bibitem{Lee:2007:MVR}
S.~L. Lee,
  \href{https://citeseerx.ist.psu.edu/viewdoc/download?doi=10.1.1.389.951\&rep=rep1\&type=pdf}{Mean
  value representations and curvatures of compact convex hypersurfaces},
  {P}reprint (2007).
\newline\urlprefix\url{https://citeseerx.ist.psu.edu/viewdoc/download?doi=10.1.1.389.951\&rep=rep1\&type=pdf}

\bibitem{Lee:2009:VFM}
S.~L. Lee, Vector fields for mean value coordinates, SIAM Journal on
  Mathematical Analysis 40~(6) (2009) 2437--2450.

\bibitem{Belyaev:2013:SLD}
A.~Belyaev, P.-A. Fayolle, A.~Pasko, Signed {$L_p$}-distance fields,
  Computer-Aided Design 45~(2) (2013) 523--528.

\end{thebibliography}

\end{document}